\documentclass[11pt,twoside,reqno]{amsart}

\usepackage{microtype}
\usepackage{cite}
\usepackage[OT1]{fontenc}
\usepackage{lmodern}
\usepackage{type1cm}
\usepackage{mathtools}
\usepackage{amssymb}
\usepackage[shortlabels]{enumitem}
\usepackage{comment}
\usepackage{xcolor}
\usepackage{graphicx}
\usepackage{dsfont}
\usepackage{aliascnt}
\usepackage{tikz}
\usetikzlibrary{calc,patterns,angles,quotes,arrows}

\usepackage{geometry}
\usepackage{hyperref}
\hypersetup{
  colorlinks=true,
  linkcolor=black,
  anchorcolor=black,
  citecolor=black,
  filecolor=black,      
  menucolor=red,
  runcolor=black,
  urlcolor=black,
}
\usepackage[capitalize,noabbrev,nameinlink]{cleveref}

\numberwithin{equation}{section}
\crefformat{equation}{(#2#1#3)}
\Crefformat{equation}{(#2#1#3)}
\crefrangeformat{equation}{(#3#1#4)--(#5#2#6)}
\Crefrangeformat{equation}{(#3#1#4)--(#5#2#6)}
\crefmultiformat{equation}{(#2#1#3)}{ and~(#2#1#3)}{, (#2#1#3)}{, and~(#2#1#3)}
\Crefmultiformat{equation}{(#2#1#3)}{ and~(#2#1#3)}{, (#2#1#3)}{, and~(#2#1#3)}
\crefrangemultiformat{equation}{(#3#1#4)--(#5#2#6)}{ and~(#3#1#4)--(#5#2#6)}%
  {, (#3#1#4)--(#5#2#6)}{, and~(#3#1#4)--(#5#2#6)}
\Crefrangemultiformat{equation}{(#3#1#4)--(#5#2#6)}{ and~(#3#1#4)--(#5#2#6)}%
  {, (#3#1#4)--(#5#2#6)}{, and~(#3#1#4)--(#5#2#6)}
\crefformat{figure}{#2Figure #1#3}
\Crefformat{figure}{#2Figure #1#3}
\crefformat{enumi}{(#2#1#3)}
\Crefformat{enumi}{(#2#1#3)}
\crefrangeformat{enumi}{(#3#1#4)--(#5#2#6)}
\Crefrangeformat{enumi}{(#3#1#4)--(#5#2#6)}
\crefmultiformat{enumi}{(#2#1#3)}{ and~(#2#1#3)}{, (#2#1#3)}{, and~(#2#1#3)}
\Crefmultiformat{enumi}{(#2#1#3)}{ and~(#2#1#3)}{, (#2#1#3)}{, and~(#2#1#3)}
\crefrangemultiformat{enumi}{(#3#1#4)--(#5#2#6)}{ and~(#3#1#4)--(#5#2#6)}%
  {, (#3#1#4)--(#5#2#6)}{, and~(#3#1#4)--(#5#2#6)}
\Crefrangemultiformat{enumi}{(#3#1#4)--(#5#2#6)}{ and~(#3#1#4)--(#5#2#6)}%
  {, (#3#1#4)--(#5#2#6)}{, and~(#3#1#4)--(#5#2#6)}
\crefformat{theorem}{#2Theorem~#1#3}
\Crefformat{theorem}{#2Theorem~#1#3}
\crefname{maintheorem}{Theorem}{Theorems}
\Crefname{maintheorem}{Theorem}{Theorems}
\crefformat{maintheorem}{#2Theorem~#1#3}
\Crefformat{maintheorem}{#2Theorem~#1#3}
\crefformat{lemma}{#2Lemma~#1#3}
\Crefformat{lemma}{#2Lemma~#1#3}
\crefformat{proposition}{#2Proposition~#1#3}
\Crefformat{proposition}{#2Proposition~#1#3}
\crefformat{corollary}{#2Corollary~#1#3}
\Crefformat{corollary}{#2Corollary~#1#3}
\crefformat{claim}{#2Claim~#1#3}
\Crefformat{claim}{#2Claim~#1#3}
\crefformat{section}{#2\S#1#3}
\Crefformat{section}{#2\S#1#3}
\crefformat{subsection}{#2\S#1#3}
\Crefformat{subsection}{#2\S#1#3}

\theoremstyle{plain}
\newtheorem{theorem}{Theorem}[section]
\newtheorem{maintheorem}{Theorem}

\newaliascnt{corollary}{theorem}

\aliascntresetthe{corollary}
\newaliascnt{proposition}{theorem}
\newtheorem{proposition}[proposition]{Proposition}
\aliascntresetthe{proposition}

\newaliascnt{lemma}{theorem}
\newtheorem{lemma}[lemma]{Lemma}
\aliascntresetthe{lemma}
\newaliascnt{conjecture}{theorem}

\aliascntresetthe{conjecture}

\theoremstyle{remark}
\newaliascnt{remark}{theorem}
\newtheorem{remark}[remark]{Remark}
\aliascntresetthe{remark}
\newaliascnt{remarks}{theorem}

\aliascntresetthe{remarks}
\newaliascnt{example}{theorem}
\newtheorem{example}[example]{Example}
\aliascntresetthe{example}

\newtheorem*{acknowledgements}{Acknowledgements}

\newtheorem*{declaration}{Declaration}
\newtheorem*{leandevelopment}{Lean formalisation}
\newtheorem*{financial}{Funding information}

\theoremstyle{definition}
\newaliascnt{definition}{theorem}

\aliascntresetthe{definition}
\newaliascnt{question}{theorem}
\newtheorem{question}[question]{Question}
\aliascntresetthe{question}

\crefname{theorem}{Theorem}{Theorems}
\Crefname{theorem}{Theorem}{Theorems}
\crefname{lemma}{Lemma}{Lemmas}
\Crefname{lemma}{Lemma}{Lemmas}
\crefname{proposition}{Proposition}{Propositions}
\Crefname{proposition}{Proposition}{Propositions}
\crefname{corollary}{Corollary}{Corollaries}
\Crefname{corollary}{Corollary}{Corollaries}
\crefname{conjecture}{Conjecture}{Conjectures}
\Crefname{conjecture}{Conjecture}{Conjectures}
\crefname{remark}{Remark}{Remarks}
\Crefname{remark}{Remark}{Remarks}
\crefname{remarks}{Remarks}{Remarks}
\Crefname{remarks}{Remarks}{Remarks}
\crefname{example}{Example}{Examples}
\Crefname{example}{Example}{Examples}
\crefname{definition}{Definition}{Definitions}
\Crefname{definition}{Definition}{Definitions}
\crefname{question}{Question}{Questions}
\Crefname{question}{Question}{Questions}

\DeclareFontFamily{U}{mathx}{}
\DeclareFontShape{U}{mathx}{m}{n}{<-> mathx10}{}
\DeclareSymbolFont{mathx}{U}{mathx}{m}{n}
\DeclareMathAccent{\widehat}{0}{mathx}{"70}
\DeclareMathAccent{\widecheck}{0}{mathx}{"71}

\newcommand{\BB}{\mathcal{B}}

\newcommand{\LL}{\mathcal{L}}
\newcommand{\cP}{\mathcal{P}}
\newcommand{\cN}{\mathcal{N}}
\newcommand{\cI}{\mathcal{I}}
\newcommand{\cQ}{\mathcal{Q}}
\newcommand{\MM}{\mathcal{M}}

\newcommand{\cV}{\mathcal{V}}
\newcommand{\cA}{\mathcal{A}}
\newcommand{\R}{\mathbb{R}}
\newcommand{\RP}{\mathbb{RP}^1}

\newcommand{\Q}{\mathbb{Q}}
\newcommand{\Z}{\mathbb{Z}}
\newcommand{\N}{\mathbb{N}}
\newcommand{\E}{\mathbb{E}}

\newcommand{\iii}{\mathtt{i}}
\newcommand{\jjj}{\mathtt{j}}
\newcommand{\kkk}{\mathtt{k}}
\renewcommand{\lll}{\mathtt{l}}
\newcommand{\eps}{\varepsilon}
\newcommand{\fii}{\varphi}

\newcommand{\A}{\mathsf{A}}

\newcommand{\dd}{\,\mathrm{d}}

\renewcommand{\ge}{\geqslant}
\renewcommand{\le}{\leqslant}
\renewcommand{\geq}{\geqslant}
\renewcommand{\leq}{\leqslant}

\DeclareMathOperator{\dimloc}{dim_{loc}}
\DeclareMathOperator{\udimloc}{\overline{dim}_{loc}}
\DeclareMathOperator{\ldimloc}{\underline{dim}_{loc}}

\DeclareMathOperator{\dimm}{dim_M}
\DeclareMathOperator{\udimm}{\overline{dim}_M}
\DeclareMathOperator{\ldimm}{\underline{dim}_M}

\DeclareMathOperator{\dimh}{dim_H}
\DeclareMathOperator{\udimh}{\overline{dim}_H}
\DeclareMathOperator{\ldimh}{\underline{dim}_H}

\DeclareMathOperator{\dimp}{dim_p}
\DeclareMathOperator{\udimp}{\overline{dim}_p}

\DeclareMathOperator{\diml}{dim_L}
\DeclareMathOperator{\dimaff}{dim_{aff}}
\newcommand{\dimmaff}{\mathop{\widehat{\mathrm{dim}}}\nolimits_{\mathrm{aff}}}

\DeclareMathOperator*{\essinf}{ess\,inf}
\DeclareMathOperator*{\esssup}{ess\,sup}

\DeclareMathOperator{\GL}{GL_2(\mathbb{R})}

\DeclareMathOperator{\linspan}{span}
\DeclareMathOperator{\dist}{dist}
\DeclareMathOperator{\diam}{diam}
\DeclareMathOperator{\proj}{proj}

\renewcommand{\atop}[2]{\genfrac{}{}{0pt}{}{#1}{#2}}

\begin{document}

\title[Planar box-like self-affine sets with rotations]{Hausdorff dimension of planar box-like
self-affine sets with rotations}

\author[A.~K\"aenm\"aki]{Antti K\"aenm\"aki}
\address[Antti K\"aenm\"aki]
{University of Eastern Finland \\ 
  Department of Physics and Mathematics \\
  P.O.\ Box 111 \\ 
  FI-80101 Joensuu \\ 
Finland}
\email{antti@kaenmaki.net}

\author[R.~D.~Prokaj]{R. D\'aniel Prokaj}
\address[R.~D.~Prokaj]
{University of North Texas \\ 
  Department of Mathematics \\ 
  1155 Union Circle \\ 
  Denton, Texas 76203-5017 \\
USA}
\email{rudolf.prokaj@unt.edu}

\author[K.~Simon]{K\'aroly Simon}
\address[K.~Simon]
{Department of Stochastics \\
  Institute of Mathematics \\
  Budapest University of Technology and Economics \\
  M\H{u}egyetem rkp.\ 3 \\
  H-1111 Budapest \\
  Hungary}
\address[K.~Simon]
{HUN-REN--BME Stochastics Research Group \\
  M\H{u}egyetem rkp.\ 3 \\
  H-1111 Budapest \\
  Hungary}
\email{simonk@math.bme.hu}

\author[S.~Troscheit]{Sascha Troscheit}
\address[S.~Troscheit]
{Uppsala University \\
Department of Mathematics \\
Box 524\\
751 20 Uppsala\\
Sweden}
\email{sascha@troscheit.eu}

\subjclass[2020]{Primary 28A80; Secondary 37C45, 37D35}
\keywords{Self-affine set, Hausdorff dimension, affinity dimension, equilibrium state,
Ledrappier--Young formula, exponential separation}

\begin{abstract}
We study the Hausdorff dimension of proximal planar self-affine sets whose defining maps have
diagonal or anti-diagonal linear parts, with both types present. Such systems are irreducible but
not strongly irreducible. Under a projective exponential separation condition, we show that the
Hausdorff dimension equals the minimum of $2$ and the affinity dimension, thereby extending the
dimension formula of Hochman--Rapaport to all proximal irreducible planar systems. When the affinity
dimension is at most an explicit threshold determined by the Lyapunov exponents, the formula holds
assuming only that there are no exact overlaps. Above the threshold, it holds under the strong open
set condition. These dimension formulae are proved by approximation with finite diagonal subsystems.
Assuming projective exponential separation, we also determine the pointwise dimensions of the
coordinate projections of the equilibrium state. Below the threshold, the projections are not
exact-dimensional: the pointwise dimension of each exists almost everywhere and takes exactly two
values, each on a set of positive measure. At or above the threshold, both projections are
exact-dimensional with dimension one.
\end{abstract}

\maketitle

\tableofcontents

\section{Introduction}

The dimension theory of self-affine sets occupies a central place in fractal geometry. Given a
finite tuple $\Phi = (\fii_1,\ldots,\fii_N)$ of contracting affine maps of the plane, a classical
result of Hutchinson \cite{Hutchinson1981} guarantees that there is a unique non-empty compact set
$X$ satisfying $X = \bigcup_{i=1}^N \fii_i(X)$, the \emph{self-affine set} associated to $\Phi$.
Throughout the paper, we assume that $X$ is not a singleton or, equivalently, that the maps
$\fii_i$ have no common fixed point. Falconer \cite{Falconer1988} introduced the \emph{affinity
dimension}, a number determined by the linear parts of the maps alone, which always bounds the
Hausdorff dimension of $X$ from above and whose minimum with $2$, for fixed and sufficiently
contracting linear parts, equals the Hausdorff dimension for Lebesgue-almost every choice of the
translation parts. A driving problem of the field has been to replace this almost-sure statement by
verifiable conditions which guarantee the equality for explicit systems.

The answer is sensitive to the algebraic structure of the linear parts: for the classical
carpets of Bedford~\cite{Bedford1984} and McMullen~\cite{McMullen1984}, the equality typically
fails because the systems induced on the coordinate axes have exact overlaps.
B\'ar\'any--Hochman--Rapaport~\cite{BaranyHochmanRapaport2019} and
Hochman--Rapaport~\cite{HochmanRapaport2022} established the equality for proximal strongly
irreducible planar systems under separation assumptions, while
Rapaport~\cite{Rapaport2023-preprint} treated proximal diagonal systems, see
\cref{thm:hochman_rapaport,thm:rapaport}. This paper treats the intermediate case of proximal
systems that are irreducible but not strongly irreducible.
Moreover, we exhibit non-exact-dimensional projections in this concrete and deterministic
setting, a phenomenon independently observed by Morris--Sert~\cite{MorrisSert2025-preprint}
and Feng--Xie~\cite{FengXie2025-preprint} in the almost every translation setting.

\subsection{Setting and previous results} \label{sec:intro_setting}

A \emph{box-like self-affine set} is a planar self-affine set for which the image of the square
$[-1,1]^2$ under any composition of the defining maps is a rectangle with sides parallel to the
coordinate axes. Without rotations, this class
contains the Bedford--McMullen carpets \cite{Bedford1984,McMullen1984} and the
generalisations thereof considered by Gatzouras--Lalley \cite{LalleyGatzouras1992}, Bara\'nski
\cite{Baranski2007}, and Feng--Wang \cite{FengWang2005}, as well as their overlapping versions
studied by Fraser--Shmerkin \cite{FraserShmerkin2016}, B\'ar\'any--Rams--Simon
\cite{BaranyRamsSimon2016}, and Pardo-Sim\'on \cite{Pardo-Simon2019}. With rotations, these sets
have been studied much less and we refer to Fraser \cite{Fraser2012,Fraser2016},
Fraser--Jordan--Jurga \cite{FraserJordanJurga2020}, and Qiu--Wang \cite{QiuWang2023-preprint} for
examples in the literature. The term ``box-like'' was introduced by Fraser \cite{Fraser2012} and
describes self-affine sets defined by mappings $\fii_i \colon \R^2 \to \R^2$ of the form
\begin{equation*}
  \fii_i(x) = A_ix + v_i,
\end{equation*}
where $A_i \in \GL$ is either diagonal or anti-diagonal and $v_i \in \R^2$. We focus on the case
where both diagonal and anti-diagonal linear parts are present and refer to such sets as box-like
self-affine sets \emph{with rotations}, even though the isometries involved may also be
reflections.

We call a matrix tuple $\A = (A_1,\ldots,A_N) \in \GL^N$, and the associated self-affine set $X$,
\emph{proximal} if the semigroup generated by $\A$ contains a matrix with two real eigenvalues of
different absolute values, see \cref{sec:irreducibility}. As we will show in \cref{thm:A_classification}, a proximal
tuple is either \emph{strongly irreducible}, meaning that no finite union of lines through the origin
is invariant under all the matrices $A_i$; or irreducible but not strongly irreducible; or
\emph{reducible}, meaning that some line through the origin is invariant under all of them. In the
strongly irreducible case, B\'ar\'any--Hochman--Rapaport \cite{BaranyHochmanRapaport2019} determined the
Hausdorff dimension under the strong open set condition, and Hochman--Rapaport
\cite{HochmanRapaport2022} extended their result to systems with overlaps.
Most recently, Li--Pan--Zhang~\cite{LiPanZhang2026-preprint} proved equality of Hausdorff and affinity dimensions
in arbitrary dimension under Zariski density of the linear parts and the strong open set condition.

Fix a left-invariant Riemannian metric $d$ on the group of invertible affine maps of the plane.
For a word $\iii=i_0\cdots i_{n-1}$, write
$\fii_\iii=\fii_{i_0}\circ\cdots\circ\fii_{i_{n-1}}$.
Following Hochman--Rapaport~\cite[\S1.1]{HochmanRapaport2022}, we call $\Phi$
\emph{exponentially separated} if there is $c>0$ such that, for infinitely many $n$,
\begin{equation*}
  \min_{\atop{\iii,\jjj\in\{1,\ldots,N\}^n}{\iii\ne\jjj}}
  d(\fii_\iii,\fii_\jjj) \ge c^n.
\end{equation*}

\begin{theorem}[Hochman--Rapaport {\cite[Theorem~1.1]{HochmanRapaport2022}}] \label{thm:hochman_rapaport}
  If $X$ is a proximal strongly irreducible planar self-affine set satisfying the exponential separation condition, then
  \begin{equation*}
    \dimh(X) = \min\{2,\dimaff(X)\}.
  \end{equation*}
\end{theorem}

The reducible case contains all box-like self-affine sets without rotations, that is, all systems
with diagonal linear parts, including the carpets listed above. For a diagonal system
$\fii_i(x,y) = (a_ix+v_i^{(1)},b_iy+v_i^{(2)})$, proximality means that $|a_i| \ne |b_i|$ for some
$i$. Its projection onto the principal axes induces a self-similar system.
For a system $(g_i)_{i=1}^N$ of similarities of the line, write
$g_\iii=g_{i_0}\circ\cdots\circ g_{i_{n-1}}$ for a word
$\iii=i_0\cdots i_{n-1}\in\{1,\ldots,N\}^n$. We call this system \emph{exponentially separated}
if there is $c>0$ such that, for infinitely many $n$,
\begin{equation*}
  \min_{\atop{\iii,\jjj\in\{1,\ldots,N\}^n}{\iii\ne\jjj,\; g_\iii'=g_\jjj'}}
  |g_\iii(0)-g_\jjj(0)| \ge c^n.
\end{equation*}
B\'ar\'any--Rams--Simon~\cite[Theorems~A and~B]{BaranyRamsSimon2016} treated, under exponential
separation, the case in which one of the induced systems has similarity dimension at most one.

\begin{theorem}[Rapaport {\cite[Theorem~1.3]{Rapaport2023-preprint}}] \label{thm:rapaport}
  If $X$ is a proximal box-like self-affine set without rotations inducing exponentially separated
  systems on both coordinate axes, then
  \begin{equation*}
    \dimh(X) = \min\{2,\dimaff(X)\}.
  \end{equation*}
\end{theorem}

Z.~Feng~\cite[Theorem~1.3]{FengZ2025-preprint} proved the corresponding result for self-affine
measures. We call an affine system \emph{homogeneous} if all its maps have the same linear part.
For proximal homogeneous systems with positive diagonal entries and rational translations,
Rapaport--Ren~\cite[Theorem~1.2]{RapaportRen2024-preprint} showed that the absence of exact overlaps
in the induced systems suffices in place of exponential separation. Bedford--McMullen carpets,
whose induced systems have exact overlaps, show that this assumption cannot be dropped.

Beyond diagonal systems, the reducible case contains the upper-triangular tuples that are not
simultaneously diagonalisable. For these,
B\'ar\'any--Hochman--Rapaport~\cite[Proposition~6.6]{BaranyHochmanRapaport2019} and
Hochman--Rapaport~\cite[Theorem~1.7]{HochmanRapaport2022} cover, under additional assumptions, the
cases in which the invariant direction is contracted more weakly or more strongly, respectively,
while overlapping triangular carpets, such as the triangular Gatzouras--Lalley carpets of
Kolossv\'ary--Simon~\cite{KolossvarySimon2019}, are studied with carpet-type methods.

Finally, consider the proximal tuples that are irreducible but not strongly irreducible. By
\cref{thm:A_classification}, after a change of basis, their matrices are simultaneously diagonal or
anti-diagonal, with at least one anti-diagonal generator. If all generators are anti-diagonal, the
system of all pairwise compositions is diagonal and has the same attractor. Products of an odd
number of anti-diagonal matrices are anti-diagonal and hence not proximal, so this diagonal system
is again proximal, and \cref{thm:rapaport} applies to it whenever its induced systems are
exponentially separated. We study the remaining case, in which both types are present, and the
associated attractor is a box-like self-affine set with rotations. The union of the coordinate axes
is invariant under every $A_i$, but no single line through the origin is invariant under all of
them.

For box-like self-affine sets, with or without rotations, Fraser \cite{Fraser2012} calculated the
packing and Minkowski (or box-counting) dimensions under the rectangular open set condition, defined
in \cref{sec:box_like_self_affine}. We denote by $\dimh, \dimp, \dimm, \udimm, \dimmaff, \dimaff$
the Hausdorff dimension, the packing dimension, the Minkowski dimension, the upper Minkowski
dimension, Fraser's modified affinity dimension, and the affinity dimension, respectively, and by
$\proj_x$ and $\proj_y$ the orthogonal projections onto the coordinate axes, see
\cref{sec:affinity_dim}. The modified affinity dimension is bounded above by the affinity dimension,
$\dimmaff(X) \le \dimaff(X)$.

\begin{theorem}[Fraser {\cite[Theorem~2.3 and Corollary~2.5]{Fraser2012}}] \label{thm:fraser}
  If $X$ is a box-like self-affine set, then
  \begin{equation*}
    \dimp(X) = \udimm(X) \le \dimmaff(X).
  \end{equation*}
  Furthermore, if $X$ satisfies the rectangular open set condition, then
  \begin{equation*}
    \dimp(X) = \dimm(X) = \dimmaff(X),
  \end{equation*}
  and, if in addition
  $\dimm(\proj_x(X)) = 1 = \dimm(\proj_y(X))$, then
  \begin{equation*}
    \dimmaff(X) = \dimaff(X).
  \end{equation*}
\end{theorem}

\subsection{Main results} \label{sec:intro_main}

Passing from the packing and Minkowski dimensions to the Hausdorff dimension requires a finer
analysis, which we carry out by investigating the unique $s$-equilibrium state $\mu$ for
$s=\dimaff(X)$, see \cref{sec:affinity_dim}. This measure is defined on the symbolic space.
We write $\pi_*\mu$ for its image under the canonical projection $\pi$, which maps a sequence
of symbols to the point of $X$ that it codes, see \cref{sec:two_sided_shift_space}.
Under the strong separation condition and for linear parts with non-negative entries,
Fraser--Jordan--Jurga \cite{FraserJordanJurga2020} proved that $\pi_*\mu$ is
\emph{exact-dimensional}, that is, its pointwise dimension exists and is constant almost
everywhere (see \cref{sec:ledrappier_young_formula}), and satisfies a Ledrappier--Young formula.

The coordinate projections of $\pi_*\mu$, however, need not be exact-dimensional. We show that,
under the projective exponential separation condition described below, they are
\emph{double-dimensional} below an explicit threshold: their pointwise dimension exists almost
everywhere and takes exactly two values, each on a set of positive measure. In the
Ledrappier--Young formula, which writes $\dim(\pi_*\mu)$ as the sum of a projection dimension and a
slice dimension, the projection dimension is the larger of these two values.

The projective exponential separation condition is the analogue of the separation in
\cref{thm:rapaport} for the induced iterated function system of the coordinate projections:
distinct compositions of the projected maps with a common contraction ratio must stay
exponentially far apart along a subsequence of scales. It is introduced in
\cref{sec:exponential_separation} and holds, for instance, whenever all parameters are
algebraic and the projected maps generate a free semigroup. We write
$\lambda_2(\A,\mu) \le \lambda_1(\A,\mu) < 0$ for the Lyapunov exponents of $\mu$,
as defined in \cref{sec:lyapunov_exponents}.
Our first result is the following.

\begin{maintheorem} \label{thm:main1}
  Let $X$ be a proximal box-like self-affine set with rotations satisfying the projective
  exponential separation condition, and let $\mu$ be the unique $s$-equilibrium state for $s =
  \dimaff(X)$. Then the coordinate projections of $\pi_*\mu$ are double-dimensional if
  \begin{equation*}
    s < 2 - \frac{\lambda_1(\A,\mu)}{\lambda_2(\A,\mu)},
  \end{equation*}
  and exact-dimensional of dimension one otherwise.
\end{maintheorem}

\cref{thm:main1} is proved in \cref{sec:dimension_of_projections}, where the two values of the
pointwise dimension are also identified explicitly, see \cref{thm:eq_state_proj_dim}. We observed
the failure of the coordinate projections to be exact-dimensional independently of
Morris--Sert~\cite[Theorem~C]{MorrisSert2025-preprint} and Feng--Xie~\cite[Theorem~1.2(iv) and
Example~8.1]{FengXie2025-preprint}, who found it simultaneously. Morris--Sert consider the
equilibrium state when the affinity dimension is at most one, and Feng--Xie consider a Markov
measure on a system of anti-diagonal maps. Both results are proved only for Lebesgue-almost every
choice of translation parameters in the sense of Falconer. The Morris--Sert construction also
persists over an open family of linear parts within their prescribed algebraic class. By contrast,
\cref{thm:main1} exhibits the phenomenon for explicit translation parameters satisfying the
projective exponential separation condition and establishes its precise double-dimensional
structure.

For planar \emph{diagonal} self-affine measures, by contrast, the coordinate projections are
self-similar and hence exact-dimensional by Feng--Hu~\cite[Theorem~2.8]{FengHu2009}, and every
non-principal orthogonal projection is exact-dimensional under the nonconformality and
irrationality hypotheses of Py\"or\"al\"a~\cite[Theorem~1.1]{Pyorala2024-preprint}.
Py\"or\"al\"a's theorem builds on the earlier projection theorem of
B\'ar\'any--K\"aenm\"aki--Py\"or\"al\"a--Wu~\cite[Theorem~1.6]{BaranyKaenmakiPyoralaWu2023-preprint}
and removes the rectangular strong separation condition assumed there.

To determine the dimension of the attractor, we approximate with finite diagonal subsystems.
In the following theorem, we assume that the projective system has \emph{no exact overlaps},
meaning that distinct finite words in its generators define distinct similarities. Equivalently, the
projected maps generate a free semigroup. Up to the threshold in \cref{thm:main1}, this
assumption alone suffices for the dimension formula below.

\begin{maintheorem} \label{thm:main2}
  Let $X$ be a proximal box-like self-affine set with rotations whose projective system has no exact
  overlaps, and let $\mu$ be the unique $s$-equilibrium state for $s = \dimaff(X)$. Then
  \begin{equation*}
    \dimh(X) = \min\{2,\dimaff(X)\}
  \end{equation*}
  whenever at least one of the following conditions holds:
  \begin{enumerate}
    \item\label{it:main2_threshold} $\dimaff(X) \le 2-\lambda_1(\A,\mu)/\lambda_2(\A,\mu)$,
    \item\label{it:main2_sosc} $X$ satisfies the strong open set condition,
    \item\label{it:main2_pesc} $X$ satisfies the projective exponential separation condition.
  \end{enumerate}
  Furthermore, if \cref{it:main2_pesc} holds together with \cref{it:main2_threshold} or
  \cref{it:main2_sosc}, then $\pi_*\mu$ is exact-dimensional with $\dim(\pi_*\mu) = \dimaff(X) \le
  2$.
\end{maintheorem}

The strong open set condition is defined in \cref{sec:box_like_self_affine}. When $s<2$, the
exponents are distinct and the threshold of condition \cref{it:main2_threshold} lies strictly
between $1$ and $2$. For $s\ge2$, the exponents coincide and the threshold is one. Below the
threshold, no strong open set condition is required. In particular, the absence of projective
exact overlaps suffices whenever $\dimaff(X) \le 1$. Conditions \cref{it:main2_threshold} and
\cref{it:main2_sosc} force $s\le2$, and the case $s=2$ follows from the strong open set condition
alone.

Under condition \cref{it:main2_pesc}, the affinity dimension may take any value, and neither the
threshold nor the strong open set condition is needed. Projective exponential separation implies
that the projective system has no exact overlaps, as we show in \cref{sec:exponential_separation},
so under condition \cref{it:main2_pesc} this assumption of \cref{thm:main2} holds automatically.
\cref{thm:main2} is proved in \cref{sec:dimension_result_proof}. It also covers the example in
\cref{prop:ssc_no_overlap_not_pesc}, which satisfies strong separation but fails projective
exponential separation.

Under projective exponential separation, \cref{thm:main2} establishes for box-like sets with
rotations the same dimension formula as the earlier theorems of Hochman--Rapaport
(\cref{thm:hochman_rapaport}) and Rapaport (\cref{thm:rapaport}). As in the diagonal case,
exponential separation is imposed on the systems induced on the coordinate axes, since these
lines are permuted by the
linear parts, and exact overlaps of the induced systems force a dimension drop even under the
strong separation condition: the carpets show this for diagonal systems, and an example of
Fraser discussed in \cref{sec:intro_fraser} shows it for box-like sets with rotations.
Hochman--Rapaport's theorem covers strongly irreducible systems under exponential separation
of the planar maps. For the other irreducible systems, we use the basis given by
\cref{thm:A_classification}. If all generators are anti-diagonal, Rapaport's theorem applies
to the system of all pairwise compositions, provided its two induced coordinate systems are
exponentially separated. Our \cref{thm:main2} treats the remaining case, in which both matrix
types occur, under projective exponential separation. Under these respective separation
assumptions, the combined results give $\dimh(X)=\min\{2,\dimaff(X)\}$ for every
proximal irreducible planar self-affine set $X$.
In particular, \cref{thm:main2} extends beyond strong irreducibility both the dimension formula
of B\'ar\'any--Hochman--Rapaport under the strong open set condition, through condition
\cref{it:main2_sosc}, and that of Hochman--Rapaport under exponential separation, through
condition \cref{it:main2_pesc}.
Conditions \cref{it:main2_threshold} and \cref{it:main2_sosc} impose on the projective system
only the absence of exact overlaps, as does the result of Rapaport--Ren mentioned above for
homogeneous diagonal systems with rational translations.

\subsection{Relation to Fraser's formula} \label{sec:intro_fraser}

At the level of dimension, \cref{thm:main2} is the Hausdorff-dimension counterpart of Fraser's
formula in \cref{thm:fraser}, which, under the rectangular open set condition, computes the
packing and Minkowski dimensions through his \emph{modified} singular value function. For sets
covered by \cref{thm:main2} with $\dimaff(X) \le 2$, we have $\dimh(X) = \dimaff(X)$, and the
first part of \cref{thm:fraser} together with $\dimmaff(X) \le \dimaff(X)$ yields
\begin{equation*}
  \dimaff(X) = \dimh(X) \le \ldimm(X) \le \udimm(X) = \dimp(X) \le \dimmaff(X) \le \dimaff(X),
\end{equation*}
where $\ldimm$ denotes the lower Minkowski dimension. Hence the Minkowski dimension of $X$ exists
and $\dimh(X) = \dimp(X) = \dimm(X) = \dimmaff(X) = \dimaff(X)$. For these sets, this extends the
equality $\dimp(X) = \dimm(X) = \dimmaff(X)$ of \cref{thm:fraser} beyond the rectangular open set
condition. Up to the threshold of condition \cref{it:main2_threshold}, the only separation
hypothesis is the absence of projective exact overlaps, and above it, the strong open set condition
or projective exponential separation is added in place of the rectangular one. \Cref{thm:main2}
does not generalise \cref{thm:fraser}, however, since the rectangular open set condition allows
projective exact overlaps, as happens when two diagonal generators act identically on at least one
coordinate axis.

For $\dimaff(X) < 2$, the equality $\dimmaff(X) = \dimaff(X)$ needs neither the threshold nor
the strong open set condition: the diagonal subsystems in the proof of \cref{thm:main2} show
that, without projective exact overlaps, both coordinate projections of $X$ have Hausdorff and
Minkowski dimension $\min\{1,\dimaff(X)\}$, so the modified and the ordinary singular value
functions agree at $s = \dimaff(X)$. For $\dimaff(X) \ge 1$, this is the projection hypothesis
in the last part of \cref{thm:fraser}.

Conversely, projective exact overlaps can push the dimension of the coordinate projections below
$\min\{1,\dimaff(X)\}$. For proximal $X$ with $\dimaff(X) < 2$, if the coordinate projections
have Minkowski dimension below this value, then $\dimmaff(X) < \dimaff(X)$, and hence
$\dimh(X) \le \udimm(X) < \dimaff(X)$ by \cref{thm:fraser}, without any separation condition.
We prove this in \cref{sec:dimension-drop}. Consequently, the equality
$\dimh(X) = \dimaff(X)$ requires both coordinate projections to have Minkowski dimension
$\min\{1,\dimaff(X)\}$, which by the previous paragraph holds in the absence of projective
exact overlaps.
For instance, the proximal box-like self-affine set $X$ with rotations constructed by
Fraser~\cite[\S3.1]{Fraser2012} satisfies the strong separation condition, but its coordinate
projections have Minkowski dimension approximately $0.891$, whereas $\dimm(X) > 1$. It follows
that $\dimh(X) \le \dimm(X) < \dimaff(X)$, so the absence of projective exact overlaps cannot
be dropped from \cref{thm:main2}, even under the strong separation condition. The example plays
the role that the Bedford--McMullen carpets play for \cref{thm:rapaport}.

\subsection{Strategy of the proofs} \label{sec:intro_methods}

To prove \cref{thm:main1}, we express each coordinate projection of the equilibrium state as a
sum of two exact-dimensional measures and compare their dimensions. The difficulty is that, for
$0<s<2$, the equilibrium state is not quasi-Bernoulli, so the results of
B\'ar\'any--K\"aenm\"aki~\cite{BaranyKaenmaki2017} do not apply directly. We use the decomposition
of Fraser--Jordan--Jurga~\cite{FraserJordanJurga2020} into two Gibbs measures on a sub-shift
recording orientation, as described in \cref{sec:subshift,sec:gibbs_measures}.
Feng--Hu's theorem~\cite[Theorem~2.8]{FengHu2009} gives exact-dimensionality of their projections
without any separation assumption. Under projective exponential separation,
Jordan--Rapaport's theorem~\cite{JordanRapaport2021} then gives their dimensions, which differ
precisely below the threshold in \cref{thm:main1}.

To prove the dimension formula in \cref{thm:main2}, we construct Bernoulli measures on finite
diagonal subsystems of $X$. In cases \cref{it:main2_threshold,it:main2_sosc}, Kittle--Kogler's
solution of the exact overlaps conjecture~\cite[Theorem~1.1]{KittleKogler2026-preprint}
determines the dimensions of their coordinate projections. We then apply D.-J.~Feng's
Ledrappier--Young formula~\cite{Feng2023} to compute the dimensions of these measures, which tend
to the affinity dimension and give the required lower bound for $\dimh(X)$. Under projective
exponential separation, we apply \cref{thm:rapaport} directly to the subsystems to obtain the
lower bound $\min\{2,\dimaff(X)\}$ without a threshold restriction.

For the last claim of \cref{thm:main2}, we apply D.-J.~Feng's Ledrappier--Young
formula~\cite{Feng2023} to the equilibrium state itself. This takes us beyond the strong
separation required by Fraser--Jordan--Jurga~\cite{FraserJordanJurga2020}, since Feng's formula
holds for \emph{every} ergodic invariant measure on a planar self-affine set, without any
separation condition. We use the dimensions of the projected Gibbs measures computed above to
evaluate this formula, with the strong open set condition used beyond the threshold. 

\subsection{Open questions} \label{sec:intro_questions}

It is natural to ask how far the separation hypotheses in \cref{thm:main1} can be relaxed.
Without projective exponential separation, we cannot rule out that the two projected Gibbs
components have the same dimension even below the threshold, in which case the coordinate
projections would be exact-dimensional.

\begin{question} \label{q:coincident-projection-dimensions}
  Can the coordinate projections of $\pi_*\mu$, for the $s$-equilibrium state $\mu$ of a proximal
  box-like self-affine set $X$ with rotations, be exact-dimensional, with their two Gibbs components
  of equal dimension, when $s = \dimaff(X) < 2 - \lambda_1(\A,\mu)/\lambda_2(\A,\mu)$ and the
  projective system has no exact overlaps?
\end{question}

Under the hypotheses of \cref{q:coincident-projection-dimensions}, \cref{thm:main2} gives
$\dimh(X)=\dimaff(X)$. If the two component dimensions coincide below the threshold, their
common value is strictly less than $\min\{1,\dimaff(X)\}$. Feng's Ledrappier--Young formula
in \cref{thm:feng1} would therefore give
\begin{equation*}
  \dim(\pi_*\mu)<\dimh(X)=\dimaff(X).
\end{equation*}
Thus a positive answer would give an equilibrium state whose projection is not a measure of
maximal dimension.
The two components are projections of Markov measures on the sub-shift of finite type
introduced in \cref{sec:subshift}. Kittle--Kogler~\cite[Theorem~1.1]{KittleKogler2026-preprint}
only treat self-similar measures of finitely many similarities and a version for Markov measures
would extend their result as Jordan--Rapaport~\cite{JordanRapaport2021} extended Hochman's
theorem~\cite{Hochman2014} to ergodic measures. It would give both components the dimensions
they have under projective exponential separation, assuming only the absence
of projective exact overlaps. It would therefore allow the projective exponential separation
condition in \cref{thm:main1} to be replaced by the absence of projective exact overlaps, and
it would answer \cref{q:coincident-projection-dimensions} negatively.

For box-like sets with rotations, the Hausdorff dimension question remains open when
$\dimaff(X) > 2 - \lambda_1(\A,\mu)/\lambda_2(\A,\mu)$, the projective system has no exact
overlaps, and $X$ satisfies neither the strong open set condition nor projective exponential
separation. If exponential separation in \cref{thm:rapaport} could be replaced by the absence
of exact overlaps in the induced systems, then \cref{rem:threshold_bound} would give the
dimension formula of \cref{thm:main2} without any of the conditions
\cref{it:main2_threshold}--\cref{it:main2_pesc}. The proof under condition
\cref{it:main2_pesc} applies \cref{thm:rapaport} only to homogeneous systems, so this extension
is needed only in that case. Rapaport--Ren~\cite[Theorem~1.2]{RapaportRen2024-preprint} settled
it for positive diagonal entries and rational translations.

Under the strong open set condition, we can also ask whether the absence of projective exact
overlaps in \cref{thm:main2} can be replaced by the requirement that both coordinate projections
have Minkowski dimension $\min\{1,\dimaff(X)\}$. In the reducible case, it remains to determine
the Hausdorff dimension for upper-triangular tuples that are not simultaneously diagonalisable,
beyond the cases covered by
B\'ar\'any--Hochman--Rapaport~\cite[Proposition~6.6]{BaranyHochmanRapaport2019} and
Hochman--Rapaport~\cite[Theorem~1.7]{HochmanRapaport2022}.

\subsection{Organisation of the paper and back matter}

In \cref{sec:preliminaries}, we define box-like self-affine sets, classify the associated matrix
tuples, and introduce the symbolic model recording orientation and the projective system encoding
the coordinate projections. We also discuss the separation conditions, including an example
showing that strong separation does not imply projective exponential separation. In
\cref{sec:thermodynamic_formalism}, we collect the thermodynamic formalism for the equilibrium
state and its Gibbs decomposition, and prove the dimension-drop assertion from
\cref{sec:intro_fraser}. In \cref{sec:hausdorff_dimension}, we recall Feng's Ledrappier--Young
formula, identify the entropies of the equilibrium state, and prove \cref{thm:main1} in
\cref{sec:dimension_of_projections} and \cref{thm:main2} in \cref{sec:dimension_result_proof}.

\begin{acknowledgements}
  We thank Demi Allen, who took part in the early stages of this project before
  deciding to step back as an author.
  The original results of this paper had been obtained by late 2025, when
  Feng--Xie~\cite{FengXie2025-preprint} and Morris--Sert~\cite{MorrisSert2025-preprint} cited the
  work as being in preparation, with Demi Allen among the authors. In 2026, we formalised the
  results, finalised and shortened the arguments, and extended \cref{thm:main2} using Kittle--Kogler's resolution
  of the exact overlaps conjecture on the line~\cite{KittleKogler2026-preprint}.
\end{acknowledgements}

\begin{declaration}
  In preparing the final version of this paper, the authors used large language models for language editing throughout and, in places, to draft English prose from their own descriptions of the intended mathematical content. The models were also used to help locate and check references. The models contributed to mathematical arguments and calculations only in the arguments for $s \ge 2$ in \cref{thm:bernoulli_equilibrium,thm:eq_state_proj_dim,thm:main_extended}, the construction in \cref{prop:ssc_no_overlap_not_pesc}, the bound in \cref{rem:threshold_bound}, and the computations in \cref{ex:projection_double_dimensional}. In each of these, the mathematical ideas and proof strategies were the authors' own; the models drafted or edited the written arguments, filling in routine details, and assisted with checking the arguments and calculations. Separately, the models were used extensively to formalise the results and their proofs in the accompanying Lean code, and the authors have checked that the formal statements, together with the definitions on which they depend, faithfully express the corresponding results in the paper. The authors have checked every argument themselves and every reference against the source, and take full responsibility for the content of this paper.
\end{declaration}

\begin{leandevelopment}
  The paper is accompanied by a formalisation in the Lean~4 theorem prover over Mathlib.
  Every numbered statement is recorded in the development, and every numbered result proved in
  this paper has been formalised. The proofs use only Lean's standard axioms and no unproved
  literature hypotheses. For Kittle--Kogler's theorem, we use their formalisation.
  The Lean files are available at
  \url{https://github.com/sascha2718/lean-diagonal-antidiagonal} and are registered with
  Palomar~\cite{KaenmakiProkajSimonTroscheit2026-lean}. Some formal proofs use slightly different arguments from those presented here.
  We refer the reader to \texttt{paper-correspondence.yaml} for the correspondence with the paper
  and to \texttt{docs/correspondence.md} for the differences in the proofs.
\end{leandevelopment}

\begin{financial}
  KS was supported by NKFIH grants KKP~144059 and Advanced Grant~152858.
  ST was initially funded by the European Union's Horizon Europe Marie Sk\l{}odowska-Curie European
  Fellowship (MSCA-EF \#101064701) and later by Lisa \& Carl-Gustav
  Esseens fund for mathematics (2025--26) and Kungl.~Vetenskapsakademien (The Royal Swedish Academy
  of Sciences) research stipend MA2025-0070.
  This material is partly based upon work supported by the Swedish Research Council under grant
  no.~2021-06594 while some of the authors were attending the Institut Mittag-Leffler in Djursholm,
  Sweden during the semester event \emph{Interactions between fractal geometry, harmonic
  analysis, and dynamics in 2026}
\end{financial}

\section{Preliminaries} \label{sec:preliminaries}

\subsection{Box-like self-affine sets} \label{sec:box_like_self_affine}
We denote by $\GL$ the general linear group of degree $2$ over real numbers.
By an affine iterated
function system on $\R^2$ we mean a finite tuple $\Phi = (\fii_1,\ldots,\fii_N)$ of contracting
affine maps $\fii_i \colon \R^2 \to \R^2$,
\begin{equation*}
  \fii_i(x) = A_ix + v_i,
\end{equation*}
where $A_i \in \GL$ has operator norm strictly less than one, $\|A_i\|<1$, and $v_i \in
\R^2$. It is a classical result of Hutchinson \cite{Hutchinson1981} that for each affine iterated
function system there exists a unique non-empty compact set $X \subset \R^2$, called the
\emph{self-affine set}, satisfying
\begin{equation} \label{eq:self_affine_set_def}
  X = \bigcup_{i=1}^N \fii_i(X).
\end{equation}
We use the convention that whenever we speak about a self-affine set $X$, then it is automatically
accompanied with a tuple $\Phi = (\fii_1,\ldots,\fii_N)$ of affine maps which defines it. The
associated tuple of contractive matrices is denoted by $\A = (A_1,\ldots,A_N) \in \GL^N$. If the
matrices $A_i$ are of the form $D_iR_i$, where $D_i \in \GL$ is a diagonal matrix with $\|D_i\|<1$
and $R_i$ is a linear isometry satisfying $R_i([-1,1]^2)=[-1,1]^2$, then we say that the tuple $\A$
is \emph{box-like}. We call such an isometry $R_i$ \emph{rotating} if it exchanges the two
coordinate axes, and \emph{non-rotating} if it preserves each axis. The axis-exchanging
isometries include both rotations and reflections in $\R^2$. For convenience, we refer to both
as rotations. If both rotating and non-rotating linear isometries are present, then we say that
the box-like tuple $\A$ is \emph{with rotations}. By our convention, writing e.g.\ ``$X$ is a
box-like self-affine set with rotations'' means that ``the associated tuple $\A$ of matrices is
box-like with rotations''.

We also recall the standard separation conditions. The affine iterated function system $\Phi$
satisfies the \emph{strong separation condition} if the sets $\fii_1(X),\ldots,\fii_N(X)$ are
pairwise disjoint, and the \emph{strong open set condition} if there exists an open set $U \subseteq
\R^2$ such that $U \cap X \ne \emptyset$, $\bigcup_{i=1}^N \fii_i(U) \subseteq U$, and $\fii_i(U)
\cap \fii_j(U) = \emptyset$ whenever $i \ne j$. Following Feng--Wang~\cite{FengWang2005}, the system
$\Phi$ satisfies the \emph{rectangular open set condition} if there exists a non-empty open
rectangle $R \subseteq \R^2$ with sides parallel to the coordinate axes such that the sets
$\fii_1(R),\ldots,\fii_N(R)$ are pairwise disjoint subsets of $R$. The strong separation condition
implies the strong open set condition, as is seen by choosing $U$ to be a sufficiently small open
neighbourhood of $X$. As before, we say that $X$ satisfies these conditions if the defining system
$\Phi$ does.

\begin{lemma} \label{thm:box_like_sa}
  A planar self-affine set $X$ is box-like with rotations if and only if the matrices in $\A$ are
  simultaneously diagonal or anti-diagonal so that $\A$ contains both diagonal and anti-diagonal
  matrices.
\end{lemma}

\begin{proof}
  Suppose first that each $A_i$ is of the form $D_iR_i$, where $D_i \in \GL$ is a diagonal matrix
  with $\|D_i\|<1$ and $R_i$ is a linear isometry satisfying $R_i([-1,1]^2)=[-1,1]^2$. Recall that
  $R_i$ is a linear isometry if and only if $R_i$ preserves the inner product. Hence $R_ie_1 \cdot
  R_ie_2 = e_1 \cdot e_2 = 0$, where $e_1,e_2$ are the standard base vectors $e_1=(1,0)$ and
  $e_2=(0,1)$. Since $R_i([-1,1]^2) = [-1,1]^2$ and $R_i$ is linear, we have
  $R_ie_1, R_ie_2 \in \{\pm e_1,\pm e_2\}$ and therefore,
  \begin{equation*}
    R_i=
    \begin{pmatrix}
      a & 0 \\ 
      0 & b
    \end{pmatrix}
    \quad\text{or}\quad
    R_i=
    \begin{pmatrix}
      0 & a \\ 
      b & 0
    \end{pmatrix},
  \end{equation*}
  where $a,b \in \{1,-1\}$. It follows that $A_i = D_iR_i$ is either diagonal or anti-diagonal
  matrix. Since $X$ is a box-like self-affine set with rotations, $\A$ contains both diagonal and
  anti-diagonal matrices.

  Suppose now  that the matrices in $\A$ are simultaneously diagonal or anti-diagonal. Since
  \begin{equation*}
    \begin{pmatrix}
      0 & a \\ 
      b & 0
    \end{pmatrix}
    =
    \begin{pmatrix}
      a & 0 \\ 
      0 & b
    \end{pmatrix}
    \begin{pmatrix}
      0 & 1 \\ 
      1 & 0
    \end{pmatrix}
  \end{equation*}
  for all $a,b \in \R$, we see that matrices in $\A$ can always be written as a product $D_iR_i$,
  where $D_i$ is diagonal and $R_i$ is a linear isometry satisfying $R_i([-1,1]^2) = [-1,1]^2$. If
  $\A$ contains both diagonal and anti-diagonal matrices, then both rotating and non-rotating linear
  isometries are present, and $X$ is a box-like self-affine set with rotations.
\end{proof}

\subsection{Irreducibility} \label{sec:irreducibility}

Let $\RP$ be the real projective line, that is, the set of all lines through the origin in $\R^2$.
We say that $\A = (A_1,\ldots,A_N) \in \GL^N$ is \emph{irreducible} if there does not exist $V \in
\RP$ such that $A_iV=V$ for all $i \in \{ 1,\ldots,N \}$. Otherwise $\A$ is \emph{reducible}.  The
tuple $\A$ is \emph{strongly irreducible} if there does not exist a nonempty finite set $\cV \subset \RP$
such that $A_i\cV=\cV$ for all $i\in\{1,\ldots,N\}$. A matrix $A \in \GL$ is called \emph{proximal}
if it has two real eigenvalues with different absolute values. We say that $\mathsf{A}$ is
\emph{proximal} if the semigroup generated by $\A$ contains a proximal matrix.  The following lemma,
which classifies the reducibility conditions, follows immediately from
B\'ar\'any--K\"aenm\"aki--Rossi~\cite[proof of Lemma~2.2]{BaranyKaenmakiRossi2021} and is stated in
B\'ar\'any--K\"aenm\"aki--Yu~\cite[Lemma~2.7]{BaranyKaenmakiYu2021-preprint}.

\begin{lemma} \label{thm:A_classification}
  If $\A \in \GL^N$ is proximal, then precisely one of the following conditions holds:
  \begin{enumerate}
    \item\label{it:A_class1} $\A$ is strongly irreducible,
    \item\label{it:A_class2} $\A$ is irreducible but not strongly irreducible, i.e., the matrices in
      $\A$ are simultaneously diagonal or anti-diagonal in some basis, with at least one
      anti-diagonal matrix,
    \item\label{it:A_class3} $\A$ is reducible, i.e., the matrices in $\A$ are simultaneously upper
      triangular in some basis.
  \end{enumerate}
\end{lemma}

Note that if $\A$ consists only of diagonal matrices, then $\A$ is reducible. Let $\A^2 =
(A_iA_j)_{i,j \in \{1,\ldots,N\}}$ be the tuple of all pairwise products of matrices in $\A$. If
$\A$ consists only of anti-diagonal matrices, then all the matrices in $\A^2$ are diagonal. 

\begin{lemma} \label{thm:box_like_proximal}
  If $\A \in \GL^N$ is proximal and box-like, then $\A$ or $\A^2$ contains a diagonal
  matrix whose diagonal entries have distinct absolute values.
\end{lemma}

\begin{proof}
  Since $\A$ is proximal there exists a matrix $B$ in the semigroup generated by $\A$ such that $B$
  has two real eigenvalues with different absolute values. Since $\A$ is box-like, each $A_i$ is
  diagonal or anti-diagonal (recall the proof of \cref{thm:box_like_sa}) and hence so is every
  matrix in the semigroup, so such a matrix $B$ has to be diagonal.
  Therefore, since $B = A_{i_1} \cdots A_{i_n}$ for some $i_1,\ldots,i_n \in \{1,\ldots,N\}$ and $n
  \in \N$, we infer that the number of anti-diagonal matrices in the product $A_{i_1} \cdots
  A_{i_n}$ is even. We may assume all diagonal matrices in the product have diagonal
  entries of equal absolute value, as otherwise the proof is complete. For a diagonal or
  anti-diagonal matrix $C$, write $|C|$ for its entrywise absolute value. Products of such
  matrices satisfy $|CD|=|C||D|$, and the absolute value of each diagonal factor is a positive
  scalar matrix. In particular, the product contains at least one anti-diagonal matrix, since
  otherwise $|B|$ would be a positive scalar matrix. If the
  anti-diagonal matrices in the product are $B_1, \ldots, B_{2k}$, where $k \in \N^+$ such that $2k
  \le n$, then there is a constant $c>0$ such that
  \begin{equation*}
    |B| = c|B_1B_2||B_3B_4| \cdots |B_{2k-1}B_{2k}|.
  \end{equation*}
  Notice that all the products $B_{2i-1}B_{2i}$, $i \in \{1,\ldots,k\}$, are diagonal. Since the
  absolute values of the diagonal entries of $B$ are distinct, there exists $i \in
  \{1,\ldots,k\}$ such that the absolute values of the diagonal entries of
  $B_{2i-1}B_{2i}$ are distinct.
\end{proof}

By \cref{thm:box_like_sa,thm:A_classification}, apart from the all-anti-diagonal case discussed
above, a proximal tuple $\A \in \GL^N$ is irreducible but not strongly irreducible if and only if
$\A$ is box-like with rotations in a basis of unit vectors. In this case, by
\cref{thm:box_like_proximal}, $\A$ or $\A^2$ contains a diagonal matrix whose diagonal entries
have distinct absolute values. Since neither changing the basis nor passing to second iterates
affects the dimension theory of self-affine sets, these settings are interchangeable. We
therefore state our results for proximal box-like self-affine sets with rotations and, in the
following sections, adopt the most suitable normalisation without explicit mention
(\cref{thm:two-block-normalisation}): the matrices in $\A$ are diagonal or anti-diagonal with
both types present and, whenever needed, $\A$ itself contains a diagonal matrix whose diagonal
entries have distinct absolute values. The strong separation and strong open set conditions
are preserved under both reductions, and the projective exponential separation condition,
introduced in \cref{sec:exponential_separation}, always refers to the system in this normalisation.

\subsection{Two-sided shift space} \label{sec:two_sided_shift_space}
Let $(\Sigma,\sigma)$ be the two-sided full shift over the alphabet $\{1,\ldots,N\}$, that is,
$\Sigma = \{1,\ldots,N\}^\Z$ and $\sigma \colon \Sigma \to \Sigma$, where $\sigma(\iii) =
(i_{n+1})_{n \in \Z}$ for all $\iii = \cdots i_{-1}i_0i_1i_2 \cdots = (i_n)_{n \in \Z} \in \Sigma$,
is the \emph{left shift}. The corresponding one-sided full shift is $(\Sigma^+,\sigma^+)$, where
$\Sigma^+ = \{1,\ldots,N\}^\N$ and $\sigma^+ \colon \Sigma^+ \to \Sigma^+$, $\sigma^+(i_0i_1\cdots)
= i_1i_2\cdots$ is the \emph{left shift}. We endow both $\Sigma$ and $\Sigma^+$ with the product
topology and observe that the basis for the topology is the family of cylinder sets of the form
$[i_k\cdots i_\ell] = \{(j_n)_{n \in \Z} \in \Sigma : j_k=i_k,\ldots,j_\ell=i_\ell\}$ and similarly
for $\Sigma^+$. We can embed $\Sigma^+$ into $\Sigma$ by the map $e\colon\Sigma^+\to\Sigma$,
$e(i_0i_1\cdots) = \cdots 11i_0i_1\cdots$. This embedding is not a
conjugacy since $e \circ \sigma^+ \ne \sigma \circ e$, but it is injective and continuous (i.e.\ a
topological embedding). We denote $\sigma^+$ also by $\sigma$ for simplicity.

If $\iii = (i_n)_{n \in \Z} \in \Sigma$, then we write $\iii^+ = (i_n)_{n \in \N}$
and $\iii^- = (i_n)_{n \in \Z \setminus \N}$, where $\N=\{0,1,2,\dots\}$ is
the set of non-negative integers. Let $\eta \colon \Sigma \to \Sigma^+$ be the natural projection
defined by $\eta(\iii) = \iii^+$ for all $\iii \in \Sigma$. Notice that $\eta \circ \sigma = \sigma
\circ \eta$. Let $\Phi = (\fii_1,\ldots,\fii_N)$ be an affine iterated function system consisting of
maps $\fii_i \colon \R^2 \to \R^2$, $\fii_i(x) = A_ix+v_i$. The one-sided
\emph{canonical projection} $\pi^+\colon\Sigma^+\to\R^2$ associated to $\Phi$ is defined by
\begin{equation}\label{eq:canonical_projection}
  \pi^+(\iii) = v_{i_0} + \sum_{n=0}^\infty A_{i_0} \cdots A_{i_n}v_{i_{n+1}} = \lim_{n \to \infty}
  \fii_{i_0} \circ \cdots \circ \fii_{i_n}(0)
\end{equation}
for all $\iii=(i_n)_{n\in\N}\in\Sigma^+$. We define the two-sided canonical projection by
$\pi=\pi^+\circ\eta\colon\Sigma\to\R^2$. Thus $\pi(\iii)=\pi^+(\iii^+)$ depends only on
the future coordinates, and $\pi(\Sigma)=\pi^+(\Sigma^+)=X$, where $X$ is the associated
self-affine set.

Given a measure $\mu$ and a measurable function $f$, we denote the push-forward measure by $f_*\mu =
\mu \circ f^{-1}$. Let $\MM_\sigma(\Sigma)$ and $\MM_\sigma(\Sigma^+)$ be the spaces of all
$\sigma$-invariant Borel probability measures on $\Sigma$ and $\Sigma^+$, respectively. For
$\mu\in\MM_\sigma(\Sigma)$, we call $\pi_*\mu$ \emph{invariant}, and \emph{ergodic} if
$\mu$ is ergodic. We use the same terminology for $(\pi^+)_*\mu^+$ when
$\mu^+\in\MM_\sigma(\Sigma^+)$.
If $\mu \in
\MM_\sigma(\Sigma)$, then
\begin{align*}
  \eta_*\mu([i_0\cdots i_n]) &= \mu(\sigma^{-1}\eta^{-1}([i_0\cdots i_n])) \\ 
			     &= \mu(\eta^{-1}\sigma^{-1}([i_0\cdots i_n])) = \eta_*\mu(\sigma^{-1}([i_0\cdots i_n]))
\end{align*}
and $\eta_*\mu \in \MM_\sigma(\Sigma^+)$. On the other hand, if $\mu^+ \in \MM_\sigma(\Sigma^+)$,
then by setting $\mu([i_k\cdots i_\ell]) = \mu^+(\{(j_n)_{n \in \N} \in \Sigma^+ :
j_0=i_k,\ldots,j_{\ell-k}=i_\ell\})$ and extending $\mu$ to a Borel probability measure on $\Sigma$
by the Hahn--Kolmogorov extension theorem, we see that
\begin{align*}
  \mu(\sigma^{-1}([i_k\cdots i_\ell])) 
  &= \mu^+(\sigma^{-1}(\{(j_n)_{n \in \N} \in \Sigma^+ :
  j_0=i_k,\ldots,j_{\ell-k}=i_\ell\})) \\ 
  &= \mu([i_k\cdots i_\ell])
\end{align*}
and $\mu \in \MM_\sigma(\Sigma)$. Therefore the push-forward map $\eta_* \colon \MM_\sigma(\Sigma)
\to \MM_\sigma(\Sigma^+)$ defined by $\eta_*(\mu) = \eta_*\mu$ is bijective: invariance determines
the measure of every two-sided cylinder from its one-sided counterpart. This correspondence
also preserves ergodicity in both directions.
Furthermore, if $\mu^+ \in \MM_\sigma(\Sigma^+)$, then $\mu = (\eta_*)^{-1}\mu^+ \in
\MM_\sigma(\Sigma)$ and $(\pi^+)_*\mu^+ = \pi_*\mu$. As far as we are interested in the self-affine set,
we can equivalently study invariant measures on the one-sided full shift. The role of negative
sequences will be clarified in \cref{sec:lyapunov_exponents}.

\subsection{Sub-shift of finite type} \label{sec:subshift}
Let $\A = (A_1,\ldots,A_N) \in \GL^N$ be proximal and box-like with rotations as defined in
\cref{sec:box_like_self_affine,sec:irreducibility}. Recalling \cref{thm:box_like_sa}, after
relabelling the maps if necessary, there are $\ell \in \{1,\ldots,N-1\}$ and $a_i,b_i \in
(-1,1)\setminus \{0\}$ such that
\begin{equation*}
  A_i =
  \begin{pmatrix}
    a_i & 0 \\ 
    0 & b_i
  \end{pmatrix}
\end{equation*}
for all $i \in \{1,\ldots,\ell\}$ and
\begin{equation*}
  A_i =
  \begin{pmatrix}
    0 & a_i \\ 
    b_i & 0
  \end{pmatrix}
\end{equation*}
for all $i \in \{\ell+1,\ldots,N\}$. Let $Z =
(Z(i,j))_{i,j=1}^{2N}$ be the $2N \times 2N$-matrix defined by
\begin{equation*}
  Z(i,j) =
  \begin{cases}
    1, &\text{if } i \in \{1,\ldots,\ell\} \cup \{N+\ell+1,\ldots,2N\} \text{ and } j \in
    \{1,\ldots,N\}, \\ 
    1, &\text{if } i \in \{\ell+1,\ldots,N+\ell\} \text{ and } j \in \{N+1,\ldots,2N\}, \\ 
    0, &\text{otherwise},
  \end{cases}
\end{equation*}
and let
\begin{equation*}
  \Sigma_Z = \{\iii \in \{1,\ldots,2N\}^\Z : Z(i_n,i_{n+1})=1 \text{ for all }n \in \Z\}
\end{equation*}
be the two-sided sub-shift of finite type associated to $Z$. We endow $\Sigma_Z$ with the
restriction of the product topology and observe that the basis for the topology is the family of
cylinder sets. We denote the left shifts acting on $\Sigma_Z$ by $\sigma$. Since $Z^2$ is a positive
matrix, $\Sigma_Z$ is mixing, see, e.g., Walters \cite[Theorem 1.31]{Walters1982}.
The purpose of the sub-shift is to precisely
record at which times the orientation under the action of $\A$ is preserved. 

Let $\Sigma$ be as in \cref{sec:two_sided_shift_space}. We define $\tau, \omega \colon \Sigma \to
\Sigma_Z$ by setting
\begin{equation*}
  \tau(\iii) = (\tau_n(\iii))_{n \in \Z} \quad\text{and}\quad \omega(\iii) = (\omega_n(\iii))_{n \in \Z}
\end{equation*}
for all $\iii = (i_n)_{n \in \Z} \in \Sigma$, where $\tau_0(\iii)=i_0$, $\omega_0(\iii)=i_0+N$, and
\begin{equation} \label{eq:def_of_tau_omega}
  \begin{split}
    \tau_n(\iii) &= i_n + N\eps_n(\iii), \\ 
    \tau_{-n}(\iii) &= i_{-n} + N\eps_{-n}(\iii), \\ 
    \omega_n(\iii) &= i_n + N(1-\eps_n(\iii)), \\ 
    \omega_{-n}(\iii) &= i_{-n} + N(1-\eps_{-n}(\iii)),
  \end{split}
\end{equation}
where
\begin{equation} \label{eq:def_of_eps_mod2}
  \begin{split}
    \eps_n(\iii) &= \#\{j \in \{0,\ldots,n-1\} : i_j \ge \ell+1\} \mod 2, \\ 
    \eps_{-n}(\iii) &= \#\{j \in \{-1,\ldots,-n\} : i_j \ge \ell+1\} \mod 2,
  \end{split}
\end{equation}
for all $n \in \N^{+}$. By setting $\eps_0(\iii)=0$ we see that the formulae in
\cref{eq:def_of_tau_omega} are valid for all $n \in \Z$. Observe that $\eps_n((i_k)_{k \in \Z}) = 1$
if and only if the product $A_{i_0} \cdots A_{i_{n-1}}$ contains odd number of anti-diagonal
matrices for all $n \in \N^{+}$. Notice also that $i_k \le \ell$ if and only if either
both $\tau_k(\iii)$ and $\tau_{k+1}(\iii)$ are less than or equal to $N$ or both $\tau_k(\iii)$ and
$\tau_{k+1}(\iii)$ are greater than or equal to $N+1$. Analogously, $i_k \ge \ell+1$ if and only if
exactly one of $\tau_k(\iii)$ and $\tau_{k+1}(\iii)$ is less than or equal to $N$. The same
assertion holds if we replace $\tau$ by $\omega$. Hence, given that $\tau_0(\iii)=i_0$ and
$\omega_0(\iii)=i_0+N$, we obtain that $\tau$ and $\omega$ are injective with
\begin{equation} \label{eq:tau_omega_partition}
  \Sigma_Z = \tau(\Sigma) \cup \omega(\Sigma) \quad\text{and}\quad \tau(\Sigma) \cap \omega(\Sigma)
  = \emptyset.
\end{equation}
Let $\overline{\tau} = \tau^{-1}$ be the inverse of $\tau$ defined on $\tau(\Sigma)$ and, similarly,
let $\overline{\omega} = \omega^{-1}$ be the inverse defined on $\omega(\Sigma)$.
Both maps are homeomorphisms onto their images. For a finite Borel measure $\nu$ on $\Sigma_Z$, we
write $\overline{\tau}_*\nu$ and $\overline{\omega}_*\nu$ for the push-forwards of the
unnormalised restrictions $\nu|_{\tau(\Sigma)}$ and $\nu|_{\omega(\Sigma)}$, respectively. Thus
\begin{equation} \label{eq:inverse-coding-restrictions}
  \tau_*(\overline{\tau}_*\nu) = \nu|_{\tau(\Sigma)}
  \quad\text{and}\quad
  \omega_*(\overline{\omega}_*\nu) = \nu|_{\omega(\Sigma)}.
\end{equation}
The following lemma describes how $\tau$ and $\omega$ interact with the shift map, commuting with it when
$\eps_n(\iii)=0$ and interchanging when $\eps_n(\iii)=1$.

\begin{lemma} \label{thm:tau_omega_sigma_commutative}
  Let $\iii \in \Sigma$ and $n \in \Z$. If $\eps_n(\iii)=0$, then
  \begin{equation*}
    \tau(\sigma^n(\iii)) = \sigma^n(\tau(\iii)) \quad\text{and}\quad
    \omega(\sigma^n(\iii)) = \sigma^n(\omega(\iii)),
  \end{equation*}
  whereas if $\eps_n(\iii)=1$, then
  \begin{equation*}
    \tau(\sigma^n(\iii)) = \sigma^n(\omega(\iii)) \quad\text{and}\quad
    \omega(\sigma^n(\iii)) = \sigma^n(\tau(\iii)).
  \end{equation*}
\end{lemma}

\begin{proof}
  Let $S \colon \Z \to \Z$ be defined by $S(0)=0$ and $S(m+1)-S(m) = \mathds{1}_{\{i_m \ge \ell+1\}}$
  for all $m \in \Z$. Comparing with \cref{eq:def_of_eps_mod2} gives $\eps_m(\iii) = S(m) \bmod 2$
  for every $m \in \Z$. If $\jjj = \sigma^n(\iii)$, then $j_m = i_{n+m}$, so the map associated to
  $\jjj$ is $m \mapsto S(n+m)-S(n)$, and hence
  \begin{equation} \label{eq:tau_omega_sigma_commutative1}
    \eps_k(\sigma^n(\iii)) = S(n+k)-S(n) = \eps_n(\iii) + \eps_{n+k}(\iii) \mod 2
  \end{equation}
  for all $n,k \in \Z$. Fix $k \in \Z$. By \cref{eq:def_of_tau_omega} and $(\sigma^n(\iii))_k = i_{n+k}$,
  \begin{equation*}
    \tau_k(\sigma^n(\iii)) = i_{n+k} + N\eps_k(\sigma^n(\iii)) \quad\text{and}\quad
    (\sigma^n(\tau(\iii)))_k = \tau_{n+k}(\iii) = i_{n+k} + N\eps_{n+k}(\iii).
  \end{equation*}
  If $\eps_n(\iii)=0$, then \cref{eq:tau_omega_sigma_commutative1} yields $\eps_k(\sigma^n(\iii)) =
  \eps_{n+k}(\iii)$ and thus $\tau_k(\sigma^n(\iii)) = (\sigma^n(\tau(\iii)))_k$. If $\eps_n(\iii)=1$,
  then $\eps_k(\sigma^n(\iii)) = 1-\eps_{n+k}(\iii)$, so $\tau_k(\sigma^n(\iii)) = i_{n+k} +
  N(1-\eps_{n+k}(\iii)) = \omega_{n+k}(\iii) = (\sigma^n(\omega(\iii)))_k$. The identities for
  $\omega$ follow in the same way, and $k \in \Z$ was arbitrary.
\end{proof}

\subsection{Projective system} \label{sec:projective_system}
Let $\proj_x \colon \R^2 \to \R$, $\proj_x(x,y)=x$, and $\proj_y \colon \R^2 \to \R$,
$\proj_y(x,y)=y$. Interpreting the codomains as the $x$-axis and $y$-axis, we see that $\proj_x$ and
$\proj_y$ correspond to the orthogonal projections onto the $x$-axis and $y$-axis, respectively.
We use the coordinate embeddings $\iota_x,\iota_y\colon\R\to\R^2$ defined by
\begin{equation*}
  \iota_x(t)=(t,0)\quad\text{and}\quad\iota_y(t)=(0,t).
\end{equation*}
Let $\Phi = (\fii_1,\ldots,\fii_N)$ be an affine iterated function system consisting of maps $\fii_i
\colon \R^2 \to \R^2$, $\fii_i(x) = A_ix+(v_i^{(1)},v_i^{(2)})$, such that there are $a_i,b_i \in
(-1,1)\setminus\{0\}$ and $\ell \in \{1,\ldots,N-1\}$ for which
\begin{equation*}
  A_i =
  \begin{pmatrix}
    a_i & 0 \\ 
    0 & b_i
  \end{pmatrix}
\end{equation*}
for all $i \in \{1,\ldots,\ell\}$ and
\begin{equation*}
  A_i =
  \begin{pmatrix}
    0 & a_i \\ 
    b_i & 0
  \end{pmatrix}
\end{equation*}
for all $i \in \{\ell+1,\ldots,N\}$. Following Fraser--Jordan--Jurga \cite{FraserJordanJurga2020},
we define an iterated function system $\Psi = (g_1,\ldots,g_{2N})$ by setting $g_j \colon \R \to
\R$,
\begin{equation} \label{eq:def_of_g_system}
  g_j(x) =
  \begin{cases}
    a_jx + v_j^{(1)}, &\text{if } j \in \{1,\ldots,N\}, \\ 
    b_{j-N}x + v_{j-N}^{(2)}, &\text{if } j \in \{N+1,\ldots,2N\}.
  \end{cases}
\end{equation}
Observe that each $g_j$ corresponds to the action of some $\fii_i$ on one of the coordinate axes.
Namely, 
\begin{equation}\label{eq:phi_def_by_g}
  \fii_i(x,y) =
  \begin{cases}
    (g_i(x),g_{i+N}(y)), &\text{if } i \in \{1,\ldots,\ell\}, \\ 
    (g_i(y),g_{i+N}(x)), &\text{if } i \in \{\ell+1,\ldots,N\}.
  \end{cases}
\end{equation}
It follows that the functions in $\Psi$ arise as projections of the functions of $\Phi$ as
\begin{equation*}
  g_j(x) =
  \begin{cases}
    (\proj_x\circ\fii_j\circ\iota_x)(x), &\text{if } j \in \{1,\ldots,\ell\}, \\
    (\proj_x\circ\fii_j\circ\iota_y)(x), &\text{if } j \in \{\ell+1,\ldots,N\}, \\
    (\proj_y\circ\fii_{j-N}\circ\iota_y)(x), &\text{if } j \in
    \{N+1,\ldots,N+\ell\}, \\ 
    (\proj_y\circ\fii_{j-N}\circ\iota_x)(x), &\text{if } j \in \{N+\ell+1,\ldots,2N\}.
  \end{cases}
\end{equation*}

\begin{figure}
  \begin{tikzpicture}
    \begin{scope}[scale=0.55]
      %
      \draw (0,0) -- (10,0) -- (10,10) -- (0,10) -- cycle;
      %
      \filldraw[blue, opacity=0.4, very thick] (0,0) -- (3,0) -- (3,2) -- (0,2) -- cycle;
      \filldraw[blue, opacity=0.4, very thick] (4.5,1.2) -- (7.5,1.2) -- (7.5,2.5) -- (4.5,2.5) -- cycle;
      \filldraw[blue, opacity=0.4, very thick] (2,7.5) -- (4,7.5) -- (4,9.5) -- (2,9.5) -- cycle;
      \filldraw[blue, opacity=0.4, very thick] (7,6.3) -- (8,6.3) -- (8,9) -- (7,9) -- cycle;
      \filldraw[red, opacity=0.15, very thick] (1,1.5) -- (2.5,1.5) -- (2.5,4) -- (1,4) -- cycle;
      \filldraw[red, opacity=0.15, very thick] (3,5.7) -- (5,5.7) -- (5,6.7) -- (3,6.7) -- cycle;
      \filldraw[red, opacity=0.15, very thick] (8,7) -- (10,7) -- (10,10) -- (8,10) -- cycle;
      \node[scale=1] at (1.5,1) {$f_1$};
      \node[scale=1] at (6,1.8) {$f_2$};
      \node[scale=1] at (3,8.5) {$f_3$};
      \node[scale=1] at (7.5,7.75) {$f_4$};
      %
      \draw[blue, very thick] (0,-0.4) -- (3,-0.4);
      \node[blue, scale=1] at (1,-0.8) {$g_1$};
      \draw[dashed, opacity=0.4] (4.5,1.2) -- (4.5,0)  (7.5,1.2) -- (7.5,0);
      \draw[blue, very thick] (4.5,-0.4) -- (7.5,-0.4);
      \node[blue, scale=1] at (6,-0.8) {$g_2$};
      \draw[dashed, opacity=0.4] (2,7.5) -- (2,0)  (4,7.5) -- (4,0);
      \draw[blue, very thick] (2,-0.7) -- (4,-0.7);
      \node[blue, scale=1] at (3,-1.1) {$g_3$};
      \draw[dashed, opacity=0.4] (7,6.3) -- (7,0)  (8,6.3) -- (8,0);
      \draw[blue, very thick] (7,-0.7) -- (8,-0.7);
      \node[blue, scale=1] at (7.5,-1.1) {$g_4$};
      %
      \draw[blue, very thick] (-0.4,0) -- (-0.4,2);
      \node[blue, scale=1] at (-0.9,0.6) {$g_8$};
      \draw[dashed, opacity=0.4] (0,1.2) -- (4.5,1.2)  (0,2.5) -- (4.5,2.5);
      \draw[blue, very thick] (-0.7,1.2) -- (-0.7,2.5);
      \node[blue, scale=1] at (-1.2,2) {$g_9$};
      \draw[dashed, opacity=0.4] (0,7.5) -- (2,7.5)  (0,9.5) -- (2,9.5);
      \draw[blue, very thick] (-0.7,7.5) -- (-0.7,9.5);
      \node[blue, scale=1] at (-1.3,9) {$g_{10}$};
      \draw[dashed, opacity=0.4] (0,6.3) -- (7,6.3)  (0,9) -- (7,9);
      \draw[blue, very thick] (-0.4,6.3) -- (-0.4,9);
      \node[blue, scale=1] at (-1,6.8) {$g_{11}$};
      %
      \begin{scope}[shift={(12.5,0)}]
	\draw (0,0) -- (10,0) -- (10,10) -- (0,10) -- cycle;
	%
	\filldraw[blue, opacity=0.15, very thick] (0,0) -- (3,0) -- (3,2) -- (0,2) -- cycle;
	\filldraw[blue, opacity=0.15, very thick] (4.5,1.2) -- (7.5,1.2) -- (7.5,2.5) -- (4.5,2.5) -- cycle;
	\filldraw[blue, opacity=0.15, very thick] (2,7.5) -- (4,7.5) -- (4,9.5) -- (2,9.5) -- cycle;
	\filldraw[blue, opacity=0.15, very thick] (7,6.3) -- (8,6.3) -- (8,9) -- (7,9) -- cycle;
	\filldraw[red, opacity=0.4, very thick] (1,1.5) -- (2.5,1.5) -- (2.5,4) -- (1,4) -- cycle;
	\filldraw[red, opacity=0.4, very thick] (3,5.7) -- (5,5.7) -- (5,6.7) -- (3,6.7) -- cycle;
	\filldraw[red, opacity=0.4, very thick] (8,7) -- (10,7) -- (10,10) -- (8,10) -- cycle;
	\node[scale=1] at (1.75,2.75) {$f_5$};
	\node[scale=1] at (4,6.15) {$f_6$};
	\node[scale=1] at (9,8.5) {$f_7$};
	%
	\draw[dashed, opacity=0.4] (1,1.5) -- (1,0)  (2.5,1.5) -- (2.5,0);
	\draw[red, very thick] (1,-0.4) -- (2.5,-0.4);
	\node[red, scale=1] at (1.75,-0.8) {$g_5$};
	\draw[dashed, opacity=0.4] (3,5.7) -- (3,0)  (5,5.7) -- (5,0);
	\draw[red, very thick] (3,-0.4) -- (5,-0.4);
	\node[red, scale=1] at (4,-0.8) {$g_6$};
	\draw[dashed, opacity=0.4] (8,7) -- (8,0);
	\draw[red, very thick] (8,-0.4) -- (10,-0.4);
	\node[red, scale=1] at (9,-0.8) {$g_7$};
	%
	\draw[dashed, opacity=0.4] (0,1.5) -- (1,1.5)  (0,4) -- (1,4);
	\draw[red, very thick] (-0.4,1.5) -- (-0.4,4);
	\node[red, scale=1] at (-1,2.75) {$g_{12}$};
	\draw[dashed, opacity=0.4] (0,5.7) -- (3,5.7)  (0,6.7) -- (3,6.7);
	\draw[red, very thick] (-0.4,5.7) -- (-0.4,6.7);
	\node[red, scale=1] at (-1,6.2) {$g_{13}$};
	\draw[dashed, opacity=0.4] (0,7) -- (8,7);
	\draw[red, very thick] (-0.4,7) -- (-0.4,10);
	\node[red, scale=1] at (-1,8.5) {$g_{14}$};
      \end{scope}
    \end{scope}
  \end{tikzpicture}
  \caption{Notations for the self-affine IFS, $N=7$, $\ell=4$. On the left: functions with diagonal
  linear part. On the right: functions with anti-diagonal linear part.}\label{fig:boxes_projections}
\end{figure}
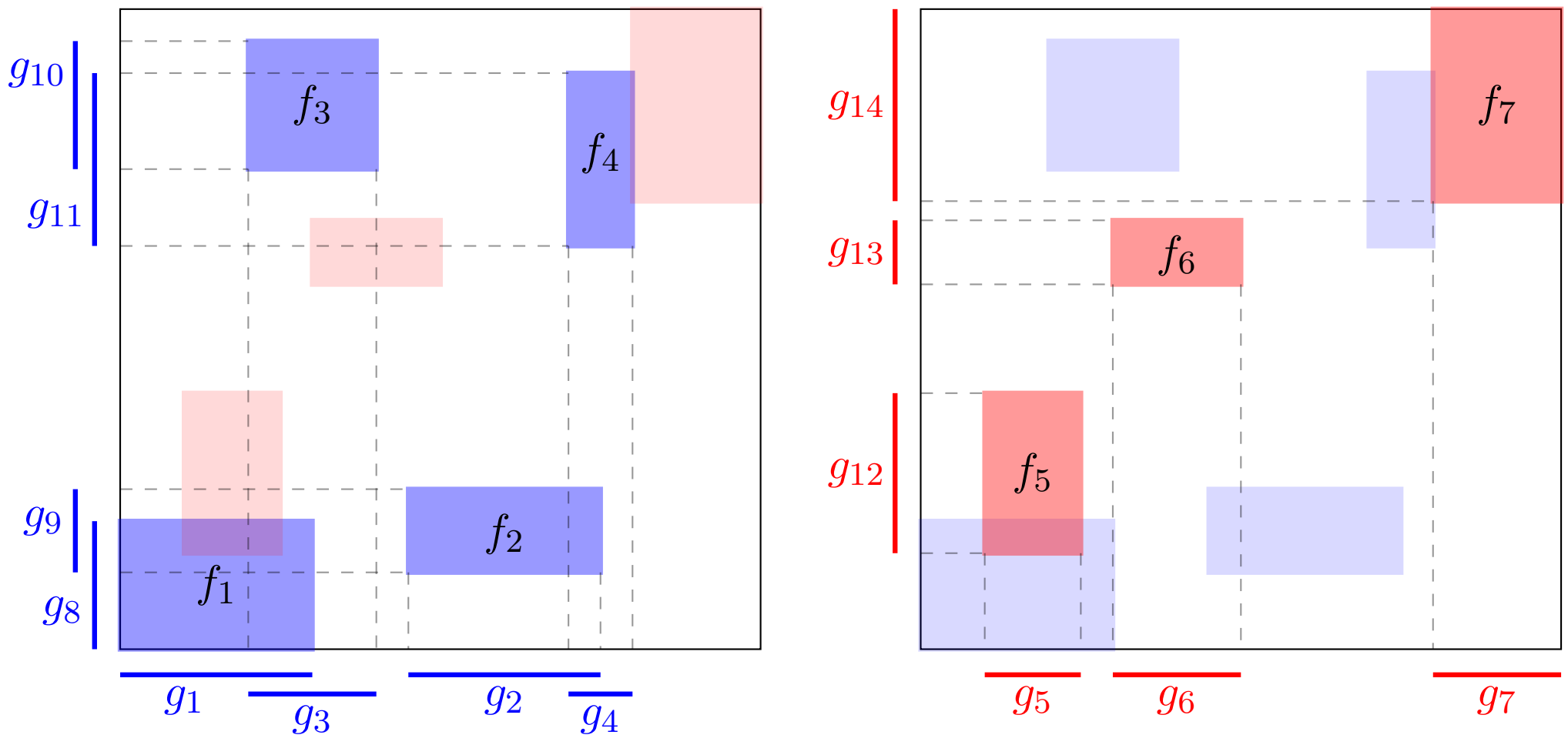

The following lemma will be useful in \cref{sec:lyapunov_exponents} when we calculate the
Lyapunov exponents of ergodic measures.

\begin{lemma} \label{thm:matrix_products_g_derivatives}
  If $X$ is a proximal box-like self-affine set with rotations, then
  \begin{align*}
    | A_{i_{-n}} \cdots A_{i_{-1}}(1,0)| 
      &= |g'_{\tau_{-n}(\iii) \cdots
      \tau_{-1}(\iii)}(0)|, \\ 
      | A_{i_{-n}} \cdots A_{i_{-1}}(0,1)| 
      &= |g'_{\omega_{-n}(\iii) \cdots
      \omega_{-1}(\iii)}(0)|,
  \end{align*}
  for all $\iii = (i_n)_{n \in \Z} \in \Sigma$ and $n \in \N$.
\end{lemma}

\begin{proof}
  Fix $\iii = (i_n)_{n \in \Z} \in \Sigma$ and for each $n \in \N$ let
  \begin{align*}
    e_n = | A_{i_{-n}} \cdots A_{i_{-1}}(1,0)|, \\ 
    f_n = | A_{i_{-n}} \cdots A_{i_{-1}}(0,1)|,
  \end{align*}
  be the absolute values of the non-zero elements of the first and second columns of $A_{i_{-n}}
  \cdots A_{i_{-1}}$, respectively. Notice that, by \cref{eq:def_of_tau_omega} and
  \cref{eq:def_of_g_system}, we have $e_1 = |g'_{\tau_{-1}(\iii)}(0)|$ and $f_1 =
  |g'_{\omega_{-1}(\iii)}(0)|$. Suppose that for a fixed $n \ge 2$ it holds that
  \begin{align*}
    e_{n-1} &= |g'_{\tau_{-(n-1)}(\iii) \cdots \tau_{-1}(\iii)}(0)|, \\ 
    f_{n-1} &= |g'_{\omega_{-(n-1)}(\iii) \cdots \omega_{-1}(\iii)}(0)|.
  \end{align*}
  By \cref{eq:def_of_g_system}, we see that the derivative of $g_j$ is constant and we have
  \begin{equation} \label{eq:derivative_g_equals1}
    \begin{split}
      |g'_{\tau_{-n}(\iii) \tau_{-(n-1)}(\iii) \cdots \tau_{-1}(\iii)}(0)| &=
      |g'_{\tau_{-n}(\iii)}(0)| e_{n-1}, \\ 
      |g'_{\omega_{-n}(\iii) \omega_{-(n-1)}(\iii) \cdots \omega_{-1}(\iii)}(0)| &=
      |g'_{\omega_{-n}(\iii)}(0)| f_{n-1}.
    \end{split}
  \end{equation}
  Notice that \cref{eq:def_of_tau_omega} and \cref{eq:def_of_g_system} give
  \begin{equation} \label{eq:derivative_g_equals2}
    \begin{split}
      g'_{\tau_{-n}(\iii)}(0) = g'_{i_{-n} + N\eps_{-n}(\iii)}(0) &= 
      \begin{cases}
	g'_{i_{-n}}(0) = a_{i_{-n}}, &\text{if } \eps_{-n}(\iii)=0, \\
	g'_{i_{-n} + N}(0) = b_{i_{-n}}, &\text{if } \eps_{-n}(\iii)=1,
      \end{cases} \\ 
      g'_{\omega_{-n}(\iii)}(0) = g'_{i_{-n} + N(1-\eps_{-n}(\iii))}(0) &= 
      \begin{cases}
	g'_{i_{-n} + N}(0) = b_{i_{-n}}, &\text{if } \eps_{-n}(\iii)=0, \\ 
	g'_{i_{-n}}(0) = a_{i_{-n}}, &\text{if } \eps_{-n}(\iii)=1.
      \end{cases}
    \end{split}
  \end{equation}
  Recall from \cref{eq:def_of_eps_mod2} that $\eps_{-n}(\iii)=0$ is equivalent to the matrix
  $A_{i_{-n}} \cdots A_{i_{-1}}$ being diagonal and $\eps_{-n}(\iii)=1$ is equivalent to the fact it
  is anti-diagonal. Recall also that $\pm e_{n-1}$ and $\pm f_{n-1}$ are the non-zero elements of
  the first and second columns of $A_{i_{-(n-1)}} \cdots A_{i_{-1}}$, respectively. Here $\pm$
  indicates that the claim holds with one of the signs. Therefore, if $\eps_{-n}(\iii)=0$, the
  product $A_{i_{-n}} \cdot (A_{i_{-(n-1)}} \cdots A_{i_{-1}})$, which gives the matrix $A_{i_{-n}} \cdots
  A_{i_{-1}}$, can have two forms:
  \begin{align*}
    \begin{pmatrix}
      a_{i_{-n}} & 0 \\ 
      0 & b_{i_{-n}}
    \end{pmatrix}
    \begin{pmatrix}
      \pm e_{n-1} & 0 \\ 
      0 & \pm f_{n-1}
    \end{pmatrix} &=
    \begin{pmatrix}
      \pm a_{i_{-n}} e_{n-1} & 0 \\ 
      0 & \pm b_{i_{-n}} f_{n-1}
    \end{pmatrix}, \\
    &=
    \begin{pmatrix}
      0 & a_{i_{-n}} \\ 
      b_{i_{-n}} & 0
    \end{pmatrix}
    \begin{pmatrix}
      0 & \pm f_{n-1} \\ 
      \pm e_{n-1} & 0
    \end{pmatrix}.
  \end{align*}
  Similarly, if $\eps_{-n}(\iii)=1$, then the product can have the following two forms:
  \begin{align*}
    \begin{pmatrix}
      a_{i_{-n}} & 0 \\ 
      0 & b_{i_{-n}}
    \end{pmatrix}
    \begin{pmatrix}
      0 & \pm f_{n-1} \\ 
      \pm e_{n-1} & 0
    \end{pmatrix} &=
    \begin{pmatrix}
      0 & \pm a_{i_{-n}} f_{n-1} \\ 
      \pm b_{i_{-n}} e_{n-1} & 0
    \end{pmatrix}, \\
    &=
    \begin{pmatrix}
      0 & a_{i_{-n}} \\ 
      b_{i_{-n}} & 0
    \end{pmatrix}
    \begin{pmatrix}
      \pm e_{n-1} & 0 \\ 
      0 & \pm f_{n-1}
    \end{pmatrix}.
  \end{align*}
  From this we see that
  \begin{equation} \label{eq:derivative_g_equals3}
    \begin{split}
      e_n &=
      \begin{cases}
	|a_{i_{-n}}| e_{n-1}, &\text{if } \eps_{-n}(\iii)=0, \\ 
	|b_{i_{-n}}| e_{n-1}, &\text{if } \eps_{-n}(\iii)=1,
      \end{cases} \\ 
      f_n &=
      \begin{cases}
	|b_{i_{-n}}| f_{n-1}, &\text{if } \eps_{-n}(\iii)=0, \\ 
	|a_{i_{-n}}| f_{n-1}, &\text{if } \eps_{-n}(\iii)=1.
      \end{cases}
    \end{split}
  \end{equation}
  Combining \cref{eq:derivative_g_equals1}, \cref{eq:derivative_g_equals2}, and
  \cref{eq:derivative_g_equals3} gives
  \begin{align*}
    e_n &= |g'_{\tau_{-n}(\iii)}(0)| e_{n-1} = |g'_{\tau_{-n}(\iii) \cdots \tau_{-1}(\iii)}(0)|, \\ 
    f_n &= |g'_{\omega_{-n}(\iii)}(0)| f_{n-1} = |g'_{\omega_{-n}(\iii) \cdots \omega_{-1}(\iii)}(0)|.
  \end{align*}
  This finishes the proof.
\end{proof}

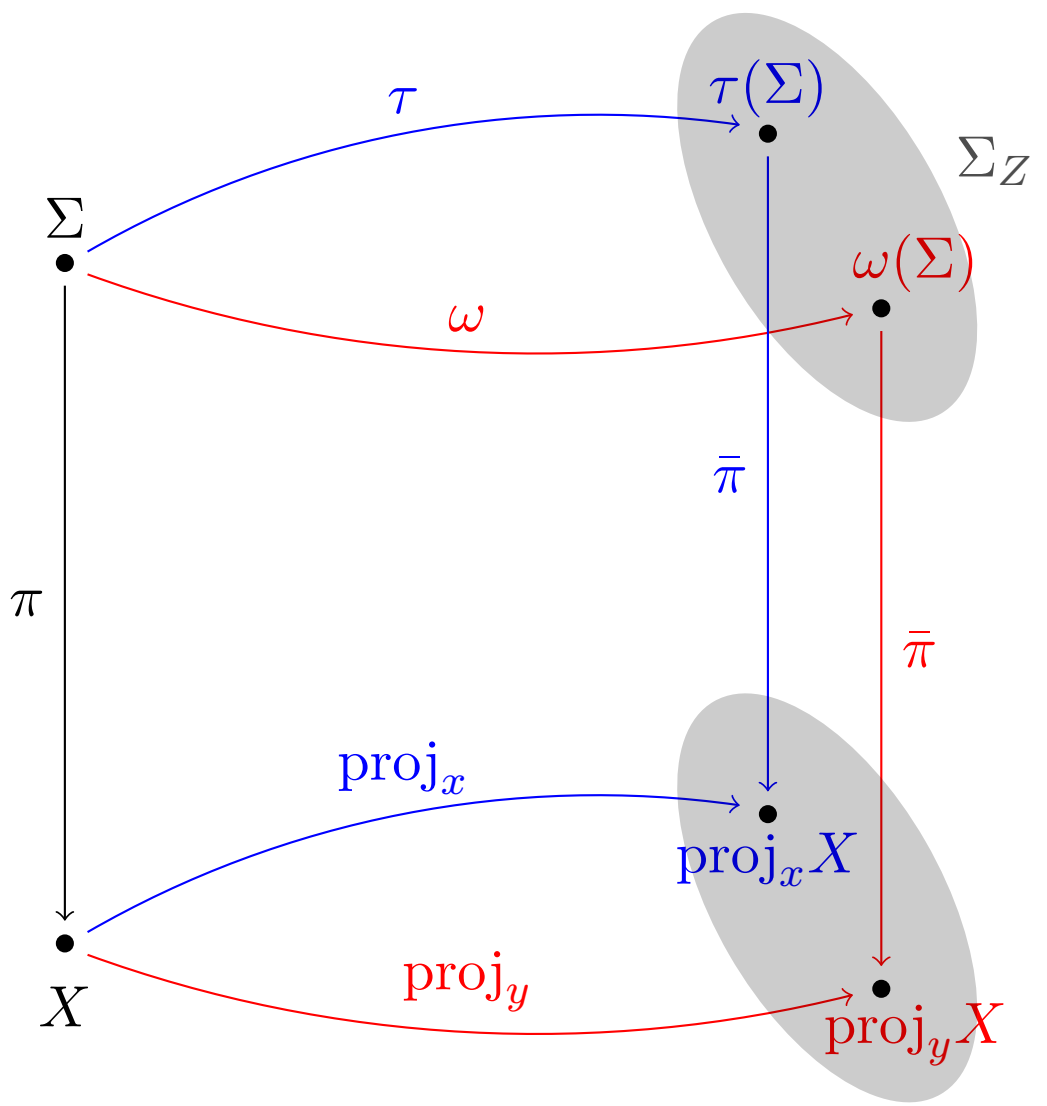
\begin{figure}
  \begin{tikzpicture}
    \begin{scope}[scale=1.5]
      \filldraw (-0.1,-0.05) circle (1pt);
      \node[below] at (-0.1,-0.15) {$X$};
      \filldraw (3,0.52) circle (1pt);
      \node[blue] at (3,0.32) {$\mbox{proj}_x X$};
      \filldraw (3.5,-0.25) circle (1pt);
      \node[red] at (3.65,-0.45) {$\mbox{proj}_y X$};
      \begin{scope}[shift={(0,3)}]
	\filldraw (-0.1,-0.05) circle (1pt);
	\node at (-0.1,0.15) {$\Sigma$};
	\filldraw (3,0.52) circle (1pt);
	\node[blue] at (3,0.72) {$\tau(\Sigma)$};
	\filldraw (3.5,-0.25) circle (1pt);
	\node[red] at (3.65,-0.05) {$\omega(\Sigma)$};
      \end{scope}
      \draw[->, blue] (0,0) arc (120:82:4.5) 
	node[midway, above] {$\mbox{proj}_x$};
      \draw[->, red] (0,-0.1) arc (250:284:5.78) 
	node[midway, above] {$\mbox{proj}_y$};
      \begin{scope}[shift={(0,3)}]
	\draw[->, blue] (0,0) arc (120:82:4.5) 
	  node[midway, above] {$\tau$};
	\draw[->, red] (0,-0.1) arc (250:284:5.78)
	  node[midway, above] {$\omega$};
      \end{scope}
      \draw[->] (-0.1,2.85) -- (-0.1,0.05) 
	node[midway, left] {$\pi$};
      \draw[->, blue] (3,3.42) -- (3,0.62) 
	node[midway, left] {$\bar{\pi}$};
      \draw[->, red] (3.5,2.65) -- (3.5,-0.15)
	node[midway, right] {$\bar{\pi}$};
      \fill[rotate=30, opacity=0.2] (2.9,-1.5) ellipse (0.5cm and 1cm);
      \begin{scope}[shift={(0,3)}]
	\fill[rotate=30, opacity=0.2] (2.9,-1.5) ellipse (0.5cm and 1cm);
	\node[opacity=0.7] at (4,0.4) {$\Sigma_Z$};
      \end{scope}
    \end{scope}
  \end{tikzpicture}
  \caption{Illustration for various spaces and maps between them.}\label{fig:symbolic_spaces}
\end{figure}

Let $\Sigma_Z$, $\tau$, and $\omega$ be as in
\cref{sec:subshift}. The \emph{canonical projection} $\overline{\pi} \colon \Sigma_Z \to \R$
associated to $\Psi$ is defined by setting
\begin{equation*}
  \overline{\pi}(\iii) = \lim_{n \to \infty} g_{i_0} \circ \cdots \circ g_{i_n}(0)
\end{equation*}
for all $\iii = (i_n)_{n \in \Z} \in \Sigma_Z$. By \cref{eq:phi_def_by_g} and
\cref{eq:canonical_projection} it is straightforward to see that
\begin{equation} \label{eq:commutative_diagram}
  \begin{split}
    \proj_x(\pi(\iii)) &= \overline{\pi}(\tau(\iii)), \\ 
    \proj_y(\pi(\iii)) &= \overline{\pi}(\omega(\iii))
  \end{split}
\end{equation}
for all $\iii \in \Sigma$. Therefore, the projective system $\Psi$ together with the sub-shift
$\Sigma_Z$ introduces us a way to study the coordinate projections of the box-like self-affine set.
The introduced symbolic spaces and projections are presented on \cref{fig:symbolic_spaces}.

\subsection{Exponential separation} \label{sec:exponential_separation}

The projective system $\Psi$ introduced in \cref{sec:projective_system} encodes the coordinate
projections of $X$ through the canonical projection $\overline{\pi}$. To control the dimension of
these projections, we now impose a quantitative separation condition on $\Psi$, requiring that
distinct compositions of the maps $g_i$ with a common derivative stay exponentially far apart, and
we compare it with the standard separation condition for the planar system $\Phi$. Write
\begin{equation*}
  d(g_i,g_j) =
  \begin{cases}
    |g_i(0)-g_j(0)|, &\text{if } g_i' = g_j', \\
    \infty, &\text{if } g_i' \ne g_j'.
  \end{cases}
\end{equation*}
For a finite word $\iii = i_1 \cdots i_n \in \{1,\ldots,2N\}^n$, write $g_\iii = g_{i_1} \circ
\cdots \circ g_{i_n}$ and
\begin{equation*}
  \Delta_n^\Psi = \min\{d(g_\iii,g_\jjj) : \iii \ne \jjj \in \{1,\ldots,2N\}^n \text{ and } g_\iii'
  = g_\jjj'\}.
\end{equation*}
Here the minimum of an empty set is understood to be $\infty$.
The projective system has \emph{no exact overlaps} if $g_\iii\ne g_\jjj$ for any two distinct
finite words $\iii,\jjj$, or equivalently, if its maps generate a free semigroup.
We say that $X$ satisfies the \emph{projective exponential separation condition} if there exist a
constant $c>0$ and a strictly increasing sequence $(n_k)_{k \in \N}$ of natural numbers such that
$\Delta_{n_k}^\Psi \ge c^{n_k}$ for every $k \in \N$. This condition is satisfied for instance if
$a_i,b_i,v_i^{(1)},v_i^{(2)}$ are all algebraic numbers and the maps in $\Psi$ generate a free
semigroup \cite[proof of Theorem~1.5]{Hochman2014}.

The projective exponential separation condition implies that the projective system has no exact
overlaps. To see this, let $\iii$ and $\jjj$ be distinct words of length $n$ with $g_\iii'=g_\jjj'$,
and let $q\in\{1,\ldots,2N\}$. The words $\iii q$ and $\jjj q$ are distinct and have equal
derivatives, and since $g_\iii-g_\jjj$ is constant,
\begin{equation*}
  g_{\iii q}(0)-g_{\jjj q}(0) = g_\iii(g_q(0))-g_\jjj(g_q(0)) = g_\iii(0)-g_\jjj(0).
\end{equation*}
Hence $\Delta_{n+1}^\Psi\le\Delta_n^\Psi$ for every $n$. Suppose now that $g_\iii=g_\jjj$ for
distinct finite words $\iii$ and $\jjj$, not necessarily of the same length. Neither word is a
prefix of the other, since $\jjj=\iii\kkk$ for a nonempty word $\kkk$ would give $g_\iii\circ
g_\kkk=g_\iii$, and hence $g_\kkk=\mathrm{id}$ by the injectivity of $g_\iii$, which is impossible
because $g_\kkk$ is a contraction. The same argument applies with $\iii$ and $\jjj$ interchanged.
Thus $\iii$ and $\jjj$ differ at some position $r\le\min\{|\iii|,|\jjj|\}$, and the words $\iii\jjj$
and $\jjj\iii$, of common length $m=|\iii|+|\jjj|$, also differ at position $r$. As
$g_{\iii\jjj}=g_\iii\circ g_\jjj=g_\jjj\circ g_\iii=g_{\jjj\iii}$, we obtain $\Delta_m^\Psi=0$, and
the monotonicity above gives $\Delta_n^\Psi=0$ for all $n\ge m$. This is incompatible with the
projective exponential separation condition.

We finish the section with an example showing that the strong separation condition does not force
the projective exponential separation condition, even when the projective maps generate a free
semigroup. To violate the projective exponential separation condition, it suffices to make the
projective gap $\Delta_n^\Psi$ decay faster than every exponential. The next lemma reduces this to a
one-dimensional problem: when a sub-collection of the maps shares its first-row entry, the
equal-derivative words it produces bound $\Delta_n^\Psi$ by the minimal
first-coordinate gap of that sub-collection.

\begin{lemma} \label{lem:proj_gap_dominated}
  Suppose that $m \in \{2,\ldots,N\}$ and that $\fii_1,\ldots,\fii_m$ have a common first-row entry
  $a$, that is, $a_1 = \cdots = a_m = a$, so that
  \begin{equation*}
    g_i(x) = ax + v_i^{(1)}
  \end{equation*}
  for every $i \in \{1,\ldots,m\}$. If
  \begin{equation*}
    \Delta_n^{\mathrm{diag}} = \min_{\iii \ne \jjj \in \{1,\ldots,m\}^n} |g_\iii(0)-g_\jjj(0)|,
  \end{equation*}
  then $\Delta_n^\Psi \le \Delta_n^{\mathrm{diag}}$ for every $n$.
\end{lemma}

\begin{proof}
  Let $\iii = i_1 \cdots i_n \in \{1,\ldots,m\}^n$. Since $m \le N$, the word $\iii$ gives a
  length-$n$ word in $\{1,\ldots,m\} \subseteq \{1,\ldots,2N\}$ of the projective system, and, as
  each of $g_1,\ldots,g_m$ has derivative $a$, its derivative is $a^n$. Two
  distinct $\iii,\jjj \in \{1,\ldots,m\}^n$ therefore give distinct words in $\{1,\ldots,2N\}^n$ of
  common derivative $a^n$, so $d(g_\iii,g_\jjj) = |g_\iii(0)-g_\jjj(0)|$. These pairs form a
  sub-collection of all distinct equal-derivative pairs in $\{1,\ldots,2N\}^n$, so the minimum
  defining $\Delta_n^\Psi$ is at most $\Delta_n^{\mathrm{diag}}$.
\end{proof}

We are now ready to construct the example announced above. A diagonal triple with a
super-exponentially small first-coordinate gap lets \cref{lem:proj_gap_dominated} push the
projective gap below every exponential, while pairwise disjoint first-level images yield the strong
separation condition. The projective system nonetheless generates a free semigroup.

\begin{proposition} \label{prop:ssc_no_overlap_not_pesc}
  There exists a system $\Phi = (\fii_1,\fii_2,\fii_3,\fii_4)$ with $N = 4$ and $\ell = 3$, for
  which $(A_1,A_2,A_3,A_4)$ is proximal and irreducible but not strongly irreducible and $|a_{i_0}|
  \ne |b_{i_0}|$ for some $i_0 \le \ell$, such that $X$ satisfies the strong separation condition
  and the maps in the projective system $\Psi$ generate a free semigroup, yet $X$ does not satisfy
  the projective exponential separation condition.
\end{proposition}

\begin{proof}
  Let $\eta_n = \exp(-n^2)$. By \cite[Theorem~2.1]{BaranyKaenmaki2021}, there is a pair
  $(\lambda,t)$ such that $0 < \lambda < \tfrac13$, $0 < t <
  \lambda/(1-\lambda)$, the homogeneous system $(x \mapsto \lambda x,\ x \mapsto \lambda x + t,\ x
  \mapsto \lambda x + 1)$ generates a free semigroup, and its minimal length-$n$ gap
  $\Delta_n^{\lambda,t}$ is bounded above by $C_0\eta_n$ for all $n \ge 1$, where $C_0 > 0$ is
  independent of $n$. Fix such a pair and write $\Delta_n = \Delta_n^{\lambda,t}$. Choose ratios
  \begin{align*}
    r_1, r_2, r_3 &\in \biggl(\frac1{20},\frac1{10}\biggr), & r_4 &\in
    \biggl(\frac1{10},\frac15\biggr), & r_5 &\in \biggl(\frac1{25},\frac1{20}\biggr)
  \end{align*}
  and translations
  \begin{align*}
    q &\in \biggl(1,\frac43\biggr),
    & u &\in \biggl(\frac32,\frac85\biggr),
    & v &\in \biggl(\frac34,\frac{19}{25}\biggr), \\
    s_1 &\in \biggl(\frac1{25},\frac1{20}\biggr),
    & s_2 &\in \biggl(\frac14,\frac3{10}\biggr),
    & s_3 &\in \biggl(\frac12,\frac{11}{20}\biggr),
  \end{align*}
  all algebraically independent over $\Q(\lambda,t)$, the subfield of $\R$ generated by $\lambda$
  and $t$ over $\Q$. Such a choice exists because the forbidden tuples form a
  countable union of proper algebraic hypersurfaces in the displayed open box. This means that no
  nonzero polynomial with coefficients in $\Q(\lambda,t)$ vanishes at the tuple
  $(r_1,\ldots,r_5,q,u,v,s_1,s_2,s_3)$. Let $\fii_1,\fii_2,\fii_3$ be diagonal and $\fii_4$
  anti-diagonal, given on $R = [0,4] \times [0,1]$ by
  \begin{align*}
    \fii_1(x,y) &= (\lambda x + q,\ r_1 y + s_1),     & \fii_2(x,y) & = (\lambda x + q + t,\ r_2 y + s_2), \\
    \fii_3(x,y) &= (\lambda x + q + 1,\ r_3 y + s_3), & \fii_4(x,y) & = (r_4 y + u,\ r_5 x + v),
  \end{align*}
  so that
  \begin{equation*}
    A_1 = \operatorname{diag}(\lambda,r_1), \quad A_2 = \operatorname{diag}(\lambda,r_2), \quad
    A_3 = \operatorname{diag}(\lambda,r_3), \quad\text{and}\quad A_4 =
    \begin{pmatrix}
      0 & r_4 \\
      r_5 & 0
    \end{pmatrix},
  \end{equation*}
  with all $a_i,b_i \in (0,1)$.

  We first check the standing hypotheses. Since $r_1$ is algebraically independent over
  $\Q(\lambda,t)$, we have $|a_1| = \lambda \ne r_1 = |b_1|$ with $i_0 = 1 \le \ell$, and $A_1$ has
  the two real eigenvalues $\lambda,r_1$ of distinct modulus, so $(A_1,A_2,A_3,A_4)$ is proximal.
  The only lines invariant under $A_1$ are the coordinate axes, which $A_4$ interchanges, so the
  tuple is irreducible. The set $\{\linspan\{e_1\}, \linspan\{e_2\}\}$ is invariant under each
  $A_i$, so it is not strongly irreducible. Thus $\Phi$ lies in the full hypothesis class.

  The induced first-coordinate maps $g_1(x) = \lambda x + q$, $g_2(x) = \lambda x + q + t$, and
  $g_3(x) = \lambda x + q + 1$ are conjugate to the homogeneous system fixed above by a common
  translation. Hence they share the ratio $\lambda$, generate a free semigroup, and their minimal
  gap $\Delta_n^{\mathrm{diag}}$ over $\{1,2,3\}^n$ equals $\Delta_n$, so
  \begin{equation*}
    \lim_{n \to \infty} \frac1n \log \Delta_n^{\mathrm{diag}} = -\infty.
  \end{equation*}
  The strong separation condition holds. The images $\fii_i(R)$ are contained in $R$, and their
  second coordinates lie in the pairwise disjoint intervals $[s_1,s_1+r_1]$, $[s_2,s_2+r_2]$,
  $[s_3,s_3+r_3]$, and $[v,v+4r_5]$ for $i \in \{1,2,3,4\}$, respectively. Thus the four first-level
  images are pairwise disjoint, with the separation carried entirely by the four $y$-bands. In
  particular, the $x$-images $[q,q+4\lambda]$ and $[q+t,q+t+4\lambda]$ of $\fii_1$ and $\fii_2$
  overlap, since $t < \lambda/(1-\lambda) < 4\lambda$. The first-coordinate separation is degenerate
  while the planar separation is intact. As $X \subseteq R$, the first-level pieces $\fii_i(X)$ are
  pairwise disjoint, so $X$ satisfies the strong separation condition.

  The projective maps generate a free semigroup. Since $r_1,\ldots,r_5$ are algebraically
  independent over $\Q(\lambda,t)$ and $\lambda \in (0,1)$, the ratios $\lambda,r_1,\ldots,r_5$ are
  multiplicatively independent. By \cref{eq:def_of_g_system}, the non-$\lambda$ maps are $g_4(x) =
  r_4x+u$, $g_5(x) = r_1x+s_1$, $g_6(x) = r_2x+s_2$, $g_7(x) = r_3x+s_3$, and $g_8(x) = r_5x+v$. Let
  $\iii \ne \jjj \in \{1,\ldots,8\}^*$ satisfy $g_\iii = g_\jjj$. Their derivatives agree, so $\iii$
  and $\jjj$ use each of the symbols $4,5,6,7,8$ equally often and contain the same total number of
  the $\lambda$-symbols $1,2,3$. In particular, $|\iii| = |\jjj|$. If both words lie in
  $\{1,2,3\}^*$, then $g_\iii = g_\jjj$ contradicts freeness of the three-map subsystem conjugate to
  the homogeneous system fixed above. We may therefore suppose that at least one of the words
  contains a symbol in $\{4,5,6,7,8\}$.

  Derivative equality ensures that both words have a first non-$\lambda$ symbol. Write the
  corresponding decompositions as
  \begin{equation*}
    \iii = \kkk j \lll \quad \text{and} \quad \jjj = \kkk' j' \lll',
  \end{equation*}
  where $\kkk,\kkk' \in \{1,2,3\}^*$ and $j,j' \in \{4,5,6,7,8\}$. The identity $g_\iii(0)-g_\jjj(0)
  = 0$ is a polynomial identity in the algebraically independent variables
  $r_1,\ldots,r_5,q,u,v,s_1,s_2,s_3$ over $\Q(\lambda,t)$. In that identity, the terms involving one
  of the translation variables $u,v,s_1,s_2,s_3$ and no factor from $r_1,\ldots,r_5$ come only from
  the first non-$\lambda$ symbol of the word: earlier symbols contribute only translations from
  $g_1,g_2,g_3$, while every later contribution is multiplied by the ratio of this first
  non-$\lambda$ symbol. Thus the no-$r$ coefficient of the relevant translation variable is
  $\lambda^{|\kkk|}$ on the first word and $\lambda^{|\kkk'|}$ on the second word, after matching
  the variable determined by the first non-$\lambda$ symbol. Hence $j = j'$ and $|\kkk| = |\kkk'|$.
  Comparing the remaining terms with no factor from $r_1,\ldots,r_5$ gives $g_\kkk(0) =
  g_{\kkk'}(0)$, and since $g_\kkk' = g_{\kkk'}' = \lambda^{|\kkk|}$, freeness of the three-map
  subsystem gives $\kkk = \kkk'$. Thus the common prefix $\kkk j$ agrees in the two words.
  Cancelling the invertible affine map $g_{\kkk j}$ and repeating this argument inductively gives
  $\iii = \jjj$, a contradiction. Hence the maps in $\Psi$ generate a free semigroup.

  Finally, the projective exponential separation condition fails. By \cref{lem:proj_gap_dominated}
  with $m = 3$ we have $\Delta_n^\Psi \le \Delta_n^{\mathrm{diag}}$ for every $n$, and therefore
  \begin{equation*}
    \limsup_{n \to \infty} \frac1n \log \Delta_n^\Psi = -\infty.
  \end{equation*}
  Thus no $c > 0$ can give the projective exponential separation lower bound along an unbounded
  sequence of lengths.
\end{proof}

The example in \cref{prop:ssc_no_overlap_not_pesc} satisfies the hypotheses of \cref{thm:main2}:
the interior of $R$ is a feasible open set meeting $X$, and
$\sum_{i=1}^4|\det A_i|=\lambda(r_1+r_2+r_3)+r_4r_5<11/100<1$, so $\dimaff(X)<2$.
Thus $\dimh(X)=\dimaff(X)$ despite the failure of projective exponential separation.

\Cref{prop:ssc_no_overlap_not_pesc} also shows that exponential separation of the planar affine
maps does not imply projective exponential separation. Indeed, under strong separation, put
$\delta=\min_{i\ne j}\dist(\fii_i(X),\fii_j(X))>0$ and
$r=\min_i\alpha_2(A_i)>0$. Distinct words $\iii,\jjj$ of length $n$ have a common prefix of
length at most $n-1$, and hence
\begin{equation*}
  \dist(\fii_\iii(X),\fii_\jjj(X))\ge\delta r^{n-1}.
\end{equation*}
Evaluating the maps at a fixed $x\in X$ gives the corresponding exponential lower bound on their
distance in the norm on affine coefficients, and hence exponential separation of the planar maps, as
defined before \cref{thm:hochman_rapaport}.

\section{Thermodynamic formalism} \label{sec:thermodynamic_formalism}
\subsection{Singular value function} \label{sec:Singular_value_function}
Recall that the singular values of $A \in \GL$ are defined as the non-negative square roots of
the eigenvalues of the positive-semidefinite matrix $A^\top A$ and are denoted by $\alpha_1(A)$ and
$\alpha_2(A)$ in non-increasing order. The identities $\alpha_1(A) = \|A\|$ and $\alpha_2(A) =
\|A^{-1}\|^{-1}$ are standard. We also have $|\det(A)| = \alpha_1(A)\alpha_2(A)$ since the singular
values $\alpha_1(A)$ and $\alpha_2(A)$ are the lengths of the major and minor axes of the ellipse
$A(B(0,1))$. For each $A \in \GL$ and $s \ge 0$, we define the \emph{singular value function} by
\begin{equation*}
  \fii^s(A) =
  \begin{cases}
    \|A\|^s = \alpha_1(A)^s, &\text{if } 0 \le s \le 1, \\
    \|A\|^{2-s}|\det(A)|^{s-1} = \alpha_1(A)\alpha_2(A)^{s-1}, &\text{if } 1 < s \le 2, \\
    |\det(A)|^{s/2} = \alpha_1(A)^{s/2}\alpha_2(A)^{s/2}, &\text{if } s > 2.
  \end{cases}
\end{equation*}
The value $\fii^s(A)$ represents an estimate 
of the
$s$-dimensional volume of the ellipse $A(B(0,1))$. Since the determinant is multiplicative and the
operator norm is sub-multiplicative, we have $\varphi^s(AB) \le \varphi^s(A)\varphi^s(B)$ for all
$A,B \in \GL$ and $s \ge 0$. Notice also that the singular value $\alpha_1$ is sub-multiplicative
and the singular value $\alpha_2$ is super-multiplicative as $\alpha_1(AB) \le
\alpha_1(A)\alpha_1(B)$ and $\alpha_2(AB) \ge
\alpha_2(A)\alpha_2(B)$ for all $A,B \in \GL$.

\subsection{Lyapunov exponents} \label{sec:lyapunov_exponents}
Let $\Sigma$ be as in \cref{sec:two_sided_shift_space}. Given a measure $\mu \in \MM_\sigma(\Sigma)$ and
$\mathsf{A} = (A_1,\ldots,A_N) \in \GL^N$, we define the $i^{\textrm{th}}$ \emph{Lyapunov
exponent} of $\mu$ by setting
\begin{equation*}
  \lambda_i(\A,\mu) = \lim_{n \to \infty} \frac{1}{n} \int_{\Sigma} \log\alpha_i(A_{i_0}\cdots
  A_{i_{n-1}}) \dd\mu(\iii)
\end{equation*}
for both $i \in \{ 1,2 \}$. Since the singular value $\alpha_1$ is sub-multiplicative, the singular
value $\alpha_2$ is super-multiplicative, and $\mu$ is invariant, the limits above exist by Fekete's
lemma. Notice also that $-\infty < \lambda_2(\A,\mu) \le \lambda_1(\A,\mu) < \infty$. If
$\|A_i\|<1$ for every $i$, then $\lambda_1(\A,\mu)<0$.
Suppose that $\mu$ is also ergodic and $\lambda_2(\A,\mu) < \lambda_1(\A,\mu)$, then by
\cite[Theorem 4.2]{MorrisShmerkin2019} it follows that $\A$ is proximal.
Furthermore, by Oseledets' multiplicative ergodic theorem \cite{Oseledec1968}, for $\mu$-almost
every $\iii = (i_n)_{n \in \Z} \in \Sigma$ there are measurable linear subspaces
\begin{equation*}
  \{0\} = V_\iii^2 \subsetneq V_\iii^1 \subsetneq V_\iii^0 = \R^2
\end{equation*}
such that
\begin{enumerate}
  \item\label{i:V1} $A_{i_{-1}}V_\iii^i = V_{\sigma^{-1}(\iii)}^i$,
  \item\label{i:V2} $\dim(V_\iii^i) = 2-i$,
  \item\label{i:V3} $\lim_{n \to \infty} \tfrac{1}{n} \log | A_{i_{-n}} \cdots
    A_{i_{-1}}v| = \lambda_{i+1}(\A,\mu)$ for $i\in\{0,1\}$ and all
    $v\in V_\iii^i\setminus V_\iii^{i+1}$.
\end{enumerate}
Notice that there exists a Borel set $\Sigma' \subseteq \Sigma$ with $\sigma(\Sigma')=\Sigma'$ and
$\mu(\Sigma')=1$ such that the above properties \cref{i:V1}--\cref{i:V3} hold for all $\iii \in
\Sigma'$. The subspaces $V_\iii^i$ depend only on $i$ and $\iii^-$ since, by the above properties,
we have
\begin{equation} \label{eq:V_lyapunov}
  V_\iii^i = \{0\} \cup \biggl\{ v \in \R^2 \setminus \{0\} : \lim_{n\to\infty}\frac{1}{n}
  \log| A_{i_{-n}} \cdots A_{i_{-1}} v| \leq \lambda_{i+1}(\A,\mu) \biggr\}
\end{equation}
for $i\in\{0,1\}$ and all $\iii=(i_n)_{n\in\Z}\in\Sigma'$.
For example, see Froyland--Lloyd--Quas~\cite[Theorem~4.1]{FroylandLloydQuas2010}.

We also define $\Pi \colon \Sigma' \to \RP$ by setting
\begin{equation*}
  \Pi(\iii) = V_\iii^1
\end{equation*}
for all $\iii \in \Sigma'$. The push-forward measure $\Pi_*\mu$ is called the \emph{Furstenberg
measure} associated to $\A$ and $\mu$. Notice that if the matrices in $\A$ are diagonal or
anti-diagonal, then by \cref{i:V3}, $V_\iii^1$ is one of the coordinate axes for each
$\iii \in \Sigma'$ and the Furstenberg measure is supported on the union of the axes.
We choose $\Sigma'$ to contain every sequence sharing the past of one of its elements.
This is possible because the defining growth conditions depend only on past coordinates.
In this case, we define
\begin{equation} \label{eq:sigma_tau_omega}
  \begin{split}
    \Sigma^\tau &= \{\iii \in \Sigma' : V_\iii^1 \text{ is the $y$-axis}\}, \\ 
    \Sigma^\omega &= \{\iii \in \Sigma' : V_\iii^1 \text{ is the $x$-axis}\}.
  \end{split}
\end{equation}
Notice that
\begin{equation} \label{eq:tau_omega_partition2}
  \Sigma' = \Sigma^\tau \cup \Sigma^\omega \quad\text{and}\quad \Sigma^\tau \cap \Sigma^\omega =
  \emptyset
\end{equation}
and the sequences $\iii^-$ represent directions where the Lyapunov exponent gets realised. 
Similarly, if $\mu \in \MM_\sigma(\Sigma_Z)$, where $\Sigma_Z$ is the sub-shift of finite type
defined in \cref{sec:subshift}, and $\Psi = (g_1,\ldots,g_{2N})$ is the projective system defined in
\cref{sec:projective_system}, then the \emph{Lyapunov exponent} of $\mu$ is defined as
\begin{equation*}
  \overline{\lambda}(\A,\mu) = \lim_{n \to \infty}\frac{1}{n}\int_{\Sigma_Z}\log |(g_{i_0}
  \circ \cdots \circ g_{i_{n-1}})'(0)|\dd\mu(\iii),
\end{equation*}
where the derivative depends only on $\A$ for all $\iii = (i_k)_{k \in \Z} \in \Sigma_Z$ and $n \in
\N$. To simplify notation, we write $g_{i_\ell \cdots i_m} = g_{i_\ell} \circ \cdots \circ g_{i_m}$
for all $\ell,m \in \Z$ with $\ell<m$ and $(i_k)_{k \in \Z} \in \Sigma_Z$. Observe that if $\mu \in
\MM_\sigma(\Sigma_Z)$ is ergodic, then
\begin{equation} \label{eq:lyapunov_exponent_gsystem}
  \overline{\lambda}(\A,\mu) = \lim_{n \to \infty}\frac{1}{n}\log 
  |g_{i_{-n} \cdots i_{-1}}'(0)|
\end{equation}
for $\mu$-almost all $\iii \in \Sigma_Z$.
As the notation suggests, there is a connection between the sets $\Sigma^\tau,
\Sigma^\omega\subset\Sigma$ and the mappings $\tau,\omega \colon \Sigma \to \Sigma_Z$. The next lemma
explains this connection in terms of Lyapunov exponents.

\begin{lemma} \label{thm:lyapunov_exponents_via_g_system}
  If $X$ is a proximal box-like self-affine set with rotations and $\mu \in \MM_\sigma(\Sigma)$
  is ergodic with $\lambda_2(\A,\mu)<\lambda_1(\A,\mu)$, then
  \begin{align*}
    \lambda_1(\A,\mu) &= 
    \begin{cases}
      \lim_{n \to \infty} \frac{1}{n} \log|g'_{\tau_{-n}(\iii) \cdots
      \tau_{-1}(\iii)}(0)|, &\text{if } \iii \in \Sigma^\tau, \\ 
      \lim_{n \to \infty} \frac{1}{n} \log|g'_{\omega_{-n}(\iii) \cdots
    \omega_{-1}(\iii)}(0)|, &\text{if } \iii \in \Sigma^\omega,
  \end{cases} \\ 
    \lambda_2(\A,\mu) &= 
    \begin{cases}
      \lim_{n \to \infty} \frac{1}{n} \log|g'_{\omega_{-n}(\iii) \cdots
      \omega_{-1}(\iii)}(0)|, &\text{if } \iii \in \Sigma^\tau, \\ 
      \lim_{n \to \infty} \frac{1}{n} \log|g'_{\tau_{-n}(\iii) \cdots
    \tau_{-1}(\iii)}(0)|, &\text{if } \iii \in \Sigma^\omega,
  \end{cases}
  \end{align*}
  for all $\iii \in \Sigma'$.
\end{lemma}

\begin{proof}
  By \cref{thm:matrix_products_g_derivatives}, we have
  \begin{align*}
    | A_{i_{-n}} \cdots A_{i_{-1}}(1,0)| 
      &= |g'_{\tau_{-n}(\iii) \cdots
      \tau_{-1}(\iii)}(0)|, \\ 
      | A_{i_{-n}} \cdots A_{i_{-1}}(0,1)| 
      &= |g'_{\omega_{-n}(\iii) \cdots
      \omega_{-1}(\iii)}(0)|,
  \end{align*}
  for all $\iii = (i_n)_{n \in \Z} \in \Sigma$ and $n \in \N^{+}$. Recall also that, by
  the property \cref{i:V3} above, for each $\iii = (i_n)_{n \in \Z} \in \Sigma'$ the limit
  \begin{equation*}
    \lim_{n \to \infty} \frac{1}{n} \log| A_{i_{-n}} \cdots A_{i_{-1}}v|
  \end{equation*}
  equals $\lambda_{1}(\A,\mu)$ for all $v \in \R^2 \setminus V_\iii^{1}$ and $\lambda_{2}(\A,\mu)$
  whenever $v \in V_\iii^1 \setminus \{0\}$. By the definition of $\Sigma^{\tau}$ and
  $\Sigma^{\omega}$, if $\iii \in \Sigma^\tau$, then $V_\iii^1 = \linspan\{(0,1)\}$, and if $\iii
  \in \Sigma^\omega$, then $V_\iii^1 = \linspan\{(1,0)\}$. Therefore, 
  \begin{align*}
    \lambda_{1}(\A,\mu) 
	  &= \lim_{n \to \infty} \frac{1}{n} \log| A_{i_{-n}} \cdots
	  A_{i_{-1}}(1,0)| = \lim_{n \to \infty} \frac{1}{n} \log |g'_{\tau_{-n}(\iii) \cdots
	  \tau_{-1}(\iii)}(0)| \\ 
	  \lambda_{2}(\A,\mu) 
	  &= \lim_{n \to \infty} \frac{1}{n} \log| A_{i_{-n}} \cdots
	  A_{i_{-1}}(0,1)| = \lim_{n \to \infty} \frac{1}{n} \log |g'_{\omega_{-n}(\iii) \cdots
	  \omega_{-1}(\iii)}(0)|
  \end{align*}
  for all $\iii \in \Sigma^\tau$ and
  \begin{align*}
    \lambda_{1}(\A,\mu) &= \lim_{n \to \infty} \frac{1}{n} \log| A_{i_{-n}} \cdots
    A_{i_{-1}}(0,1)| = \lim_{n \to \infty} \frac{1}{n} \log |g'_{\omega_{-n}(\iii) \cdots
    \omega_{-1}(\iii)}(0)| \\ 
    \lambda_{2}(\A,\mu) &= \lim_{n \to \infty} \frac{1}{n} \log| A_{i_{-n}} \cdots
    A_{i_{-1}}(1,0)| = \lim_{n \to \infty} \frac{1}{n} \log |g'_{\tau_{-n}(\iii) \cdots
    \tau_{-1}(\iii)}(0)|
  \end{align*}
  for all $\iii \in \Sigma^\omega$.
\end{proof}

\subsection{Conditional entropy and measures} \label{sec:conditional}
Let $\Sigma$ be as in \cref{sec:two_sided_shift_space}, $\mu$ be a Borel probability measure on
$\Sigma$, and $\BB$ be the Borel $\sigma$-algebra of $\Sigma$. If $\cA$ is a
sub-$\sigma$-algebra of $\BB$ and $f \in L^1(\Sigma,\BB,\mu)$, then the \emph{conditional
expectation} of $f$ given $\cA$ is the $\cA$-measurable function $\E_\mu(f\,|\,\cA)$,
unique up to $\mu$-null sets, for which
\begin{equation} \label{eq:conditional_expectation_def}
  \int_A f \dd\mu = \int_A \E(f\,|\,\cA) \dd\mu
\end{equation}
for all $A \in \cA$. Recall that $\E(f\,|\,\cN) = \int_\Sigma f \dd\mu$,
where $\cN = \{A \in \BB : \mu(A) \in \{0,1\}\}$ is the trivial $\sigma$-algebra.
If $\cQ$ is a countable $\BB$-measurable partition of $\Sigma$, then
\begin{equation*}
  I_\mu(\cQ\,|\,\cA) = -\sum_{A \in \cQ} \mathds{1}_A
  \log\E_\mu(\mathds{1}_A\,|\,\cA)
\end{equation*}
is the \emph{conditional information} of $\cQ$ given $\cA$. Here $\mathds{1}_A$ is
the characteristic function of $A \subseteq \Sigma$. The \emph{conditional entropy} of $\cQ$
given $\cA$ is
\begin{equation*}
  H_\mu(\cQ\,|\,\cA) = \int_\Sigma I_\mu(\cQ\,|\,\cA) \dd\mu.
\end{equation*}
Notice that $I_\mu(\cQ\,|\,\cN) = -\sum_{A \in \cQ} \mathds{1}_A \log\mu(A)$
and hence,
\begin{equation*}
  H_\mu(\cQ\,|\,\cN) = -\sum_{A \in \cQ} \mu(A) \log\mu(A)
\end{equation*}
is the \emph{Shannon entropy} of $\mu$ with respect to the partition $\cQ$. Since the
Shannon entropy is trivially conditioned, we denote it simply by $H_\mu(\cQ)$.

Let $\mu$ be a Borel probability measure on $\Sigma$ and $\xi$ a measurable partition of $\Sigma$.
Then by Rohlin \cite{Rohlin1952}, for every $\iii$ in a set of full $\mu$-measure, there exists
a Borel probability measure $\mu_\iii^\xi$ supported on $\xi(\iii)$, the element of $\xi$
containing $\iii$. These measures are uniquely characterised up to a null set by the following
two properties: for each Borel set $A\subseteq\Sigma$, the function
$\iii\mapsto\mu_\iii^\xi(A)$ is $\widehat\xi$-measurable and satisfies the decomposition
\begin{equation*}
  \mu(A)=\int_\Sigma\mu_\iii^\xi(A)\dd\mu(\iii),
\end{equation*}
where $\widehat\xi$ is the $\sigma$-algebra generated by $\xi$. We call these the
\emph{conditional measures of $\mu$ associated with $\xi$}. For every $f\in L^1(\Sigma,\BB,\mu)$,
\begin{equation*}
  \int_\Sigma f\dd\mu_\iii^\xi=\E_\mu(f\,|\,\widehat\xi)(\iii)
\end{equation*}
for $\mu$-almost every $\iii$.
These definitions apply equally to $\Sigma_Z$ and to full-measure Borel subsets of either shift.

For a Borel set $\Xi$ with $\mu(\Xi)>0$, we write $\mu_\Xi$ for the normalised restriction
\begin{equation*}
  \mu_\Xi(A)=\frac{\mu(A\cap\Xi)}{\mu(\Xi)}.
\end{equation*}
It agrees with the conditional measure for the partition $\{\Xi,\Xi^c\}$ at $\mu$-almost every point
of $\Xi$. The following restriction identity is also given in D.-J.~Feng~\cite[Lemma~2.9]{Feng2023}.

\begin{lemma} \label{thm:conditional_measures_of_a_conditional_measure}
  Let $\mu$ be a Borel probability measure on $\Sigma_Z$, let $\xi$ be a measurable partition,
  and let $\Xi\subseteq\Sigma_Z$ be Borel with $\mu(\Xi)>0$. Then $\mu_\iii^\xi(\Xi)>0$ and
  \begin{equation*}
    (\mu_\Xi)_\iii^\xi(A)=\frac{\mu_\iii^\xi(A\cap\Xi)}{\mu_\iii^\xi(\Xi)}
  \end{equation*}
  for $\mu_\Xi$-almost every $\iii$, simultaneously for all Borel sets $A\subseteq\Sigma_Z$.
  If $\Xi$ is a union of atoms of $\xi$, then
  \begin{equation*}
    (\mu_\Xi)_\iii^\xi=\mu_\iii^\xi
  \end{equation*}
  for $\mu$-almost every $\iii\in\Xi$. The same statements hold on $\Sigma$.
\end{lemma}

\begin{proof}
  Write $p(\iii)=\mu_\iii^\xi(\Xi)=\E_\mu(\mathds{1}_\Xi\,|\,\widehat\xi)(\iii)$.
  Since
  \begin{equation*}
    \mu_\Xi(\{p=0\})=\frac{1}{\mu(\Xi)}\int_{\{p=0\}}p\dd\mu=0,
  \end{equation*}
  we can define $K_\iii(A)=\mu_\iii^\xi(A\cap\Xi)/p(\iii)$ where $p(\iii)>0$, and assign a
  fixed probability measure where $p(\iii)=0$. This is a $\widehat\xi$-measurable probability
  kernel, supported on $\xi(\iii)$ for $\mu_\Xi$-almost every $\iii$. For each
  $B\in\widehat\xi$, conditional expectation gives
  \begin{align*}
    \int_B K_\iii(A)\dd\mu_\Xi(\iii)
    &=\frac{1}{\mu(\Xi)}\int_B p(\iii)K_\iii(A)\dd\mu(\iii)\\
    &=\frac{1}{\mu(\Xi)}\int_B\mu_\iii^\xi(A\cap\Xi)\dd\mu(\iii)
      =\mu_\Xi(A\cap B).
  \end{align*}
  Uniqueness of disintegration proves the formula. If $\Xi$ is a union of atoms, then
  $p(\iii)=1$ and $\mu_\iii^\xi(A\cap\Xi)=\mu_\iii^\xi(A)$ for $\mu$-almost every
  $\iii\in\Xi$, proving the final assertion.
\end{proof}

The next lemma describes successive conditioning on nested measurable partitions.

\begin{lemma} \label{thm:double_conditional}
  Let $\mu$ be a Borel probability measure on $\Sigma$, and let $\xi,\xi'$ be measurable
  partitions of a Borel set $\Sigma'$ with $\mu(\Sigma')=1$ such that
  $\xi'(\iii) \subseteq \xi(\iii)$ for all $\iii \in \Sigma'$. Then
  \begin{equation*}
    (\mu_\jjj^{\xi})_\iii^{\xi'} = \mu_\iii^{\xi'}
  \end{equation*}
  for $\mu$-almost all $\jjj \in \Sigma$ and for $\mu_\jjj^{\xi}$-almost all $\iii \in \xi(\jjj)$, and 
  \begin{equation*}
    \int_\Sigma \int_\Sigma \mu_\iii^{\xi'}(A) \dd\mu_\jjj^{\xi}(\iii) \dd\mu(\jjj) = \int_\Sigma
    \mu_\iii^{\xi'}(A) \dd\mu(\iii)
  \end{equation*}
  for all Borel sets $A \subseteq \Sigma$.
\end{lemma}
\begin{proof}
  Since $\xi'$ refines $\xi$, we have $\widehat\xi\subseteq\widehat{\xi'}$, so the first
  assertion follows from \cite[Proposition~5.20]{EinsiedlerWard2011}. The second follows from
  \cite[Theorem~5.14]{EinsiedlerWard2011}, applied to the bounded measurable function
  $\iii\mapsto\mu_\iii^{\xi'}(A)$.
\end{proof}

Dirac conditional measures will allow us to recognise vanishing entropies through
\cref{thm:feng_entropy_zero}. We first note that refinements of measurable partitions
inherit this property.

\begin{lemma} \label{thm:finer_partition_dirac_mass}
  Let $\mu$ be a Borel probability measure on $\Sigma$, and let $\xi,\xi'$ be measurable
  partitions of a Borel set $\Sigma'$ with $\mu(\Sigma')=1$ such that
  $\xi'(\iii) \subseteq \xi(\iii)$ for all $\iii \in \Sigma'$. If $\mu(\{\iii \in \Sigma' :
  \mu_\iii^\xi = \delta_\iii\}) = 1$, then
  \begin{equation*}
    \mu(\{\iii \in \Sigma' : \mu_\iii^{\xi'} = \delta_\iii\}) = 1.
  \end{equation*}
\end{lemma}

\begin{proof}
  For $\mu$-almost every $\jjj$, we have $\mu_\jjj^\xi=\delta_\jjj$, and
  \cref{thm:double_conditional} holds for $\mu_\jjj^\xi$-almost every $\iii$. For these
  $\jjj$, the latter identity holds at $\iii=\jjj$. Since the conditional measure of
  $\delta_\jjj$ at $\jjj$ is $\delta_\jjj$ for any measurable partition, we obtain
  $\mu_\jjj^{\xi'}=\delta_\jjj$.
\end{proof}

\subsection{Measurable partitions and entropies} \label{sec:entropies}
Let $\A = (A_1,\ldots,A_N) \in \GL^N$ and let $\mu \in \MM_\sigma(\Sigma)$ be ergodic, with
Lyapunov exponents satisfying $\lambda_2(\A,\mu)<\lambda_1(\A,\mu)$ as defined in
\cref{sec:lyapunov_exponents}. If the canonical projection $\pi$ is as in
\cref{sec:two_sided_shift_space}, then relying on \cref{eq:V_lyapunov}, we may define
measurable partitions $\xi_i$ of the full-measure set $\Sigma'$ chosen for this measure $\mu$
in \cref{sec:lyapunov_exponents}, by setting for each $i\in\{0,1,2\}$
\begin{equation*}
  \xi_i(\iii) = \{ \jjj \in \Sigma' : \jjj^- = \iii^- \text{ and } \pi(\jjj)-\pi(\iii) \in V_\iii^i \}
\end{equation*}
to be the unique element of $\xi_i$ that contains $\iii$. We also define
\begin{equation*}
  \xi^\bot_1(\iii) = \{ \jjj \in \Sigma' : \jjj^- = \iii^- \text{ and } \pi(\jjj)-\pi(\iii) \,\bot\,
  V_\iii^1 \}.
\end{equation*}
The partitions $\xi_0$, $\xi_1$, and $\xi_2$ are measurable partitions of $\Sigma'$ in the sense of
Rokhlin by \cite[Lemma~4.2]{Feng2023}, and the same argument shows that $\xi_1^\bot$ is measurable.
These partitions use the Oseledets subspaces for $\mu$ and are kept fixed when restricting
$\mu$ or transporting it by an orientation map.

Observe that $\xi_2(\iii) \subseteq \xi_1(\iii)
\subseteq \xi_0(\iii)$ and $\xi_2(\iii) \subseteq \xi^\bot_1(\iii) \subseteq \xi_0(\iii)$ for all
$\iii \in \Sigma'$. Let
\begin{equation*}
  \cP = \{ [i] : i \in \{1,\ldots,N\} \}
\end{equation*}
be the canonical partition of $\Sigma$, where $[i]=\{(j_n)_{n\in\Z}\in\Sigma\colon j_0=i\}$.
We use the same notation for its restriction to $\Sigma'$. The $n$-level refinement
$\cP_0^{n-1}=\bigvee_{i=0}^{n-1}\sigma^{-i}(\cP)$ of the
partition $\cP$ is the collection of sets $C_0 \cap \cdots \cap C_{n-1}$ with $C_i \in
\sigma^{-i}(\cP)$. Notice that then $\cP_0^{n-1}(\iii) = [i_0 \cdots i_{n-1}]$ for
all $\iii = (i_n)_{n \in \Z} \in \Sigma$. Writing $\widehat{\xi}_i$ for the $\sigma$-algebra
generated by $\xi_i$, the \emph{pointwise $i$-dimensional entropy} is
\begin{equation*}
  h_i(\mu,\iii) = \E_\mu(I_\mu(\cP\,|\,\widehat{\xi}_i)\,|\,\cI)(\iii)
\end{equation*}
for $\mu$-almost every $\iii\in\Sigma'$, where $I_\mu$ is the conditional information defined in
\cref{sec:conditional} and $\cI = \{A \in \BB : \sigma^{-1}(A)=A\}$, and the
\emph{$i$-dimensional entropy}
\begin{align*}
  h_i(\mu) &= \int_\Sigma h_i(\mu,\iii) \dd\mu(\iii) = \int_\Sigma
  \E_\mu(I_\mu(\cP\,|\,\widehat{\xi}_i)\,|\,\cI) \dd\mu \\ &= \int_\Sigma
  I_\mu(\cP\,|\,\widehat{\xi}_i) \dd\mu = H_\mu(\cP\,|\,\widehat{\xi}_i)
\end{align*}
is the conditional entropy of $\cP$ given $\widehat{\xi}_i$ defined in
\cref{sec:conditional}. It follows from \cite[Equation (4.15)]{Feng2023} that the function $\iii
\mapsto h_i(\mu,\iii)$ is in $L^1(\Sigma,\BB,\mu)$ and from the proof of \cite[Lemma 4.6]{Feng2023}
we see that
\begin{equation} \label{eq:ergodic_entropy_constant}
  h_i(\mu) = h_i(\mu,\iii)
\end{equation}
for $\mu$-almost all $\iii \in \Sigma$ whenever $\mu$ is ergodic. Bear in mind that $h_1(\mu)$
depends also on $\A$ and $h_2(\mu)$ takes into account the overlapping structure of the self-affine
set. 

Let $\A = (A_1,\ldots,A_N) \in \GL^N$ be proximal box-like with rotations as defined in
\cref{sec:box_like_self_affine}. We may assume that the matrices $A_i$ are as in
\cref{sec:subshift}. Let $\Sigma_Z$ be the sub-shift of finite type defined in \cref{sec:subshift}
and $\overline{\pi}$ the canonical projection of the projective system defined in
\cref{sec:projective_system}. If $\nu \in \MM_\sigma(\Sigma_Z)$, then extending it by zero
gives an invariant measure on the full two-sided shift $\{1,\ldots,2N\}^\Z$.
Hence, we may similarly define
measurable partitions $\overline{\xi}_i$ of $\Sigma_Z$ by setting
\begin{align*}
  \overline{\xi}_0(\iii) &= \{ \jjj \in \Sigma_Z : \jjj^- = \iii^- \}, \\
  \overline{\xi}_1(\iii) &= \{ \jjj \in \Sigma_Z : \jjj^- = \iii^- \text{ and } \overline{\pi}(\jjj)
  = \overline{\pi}(\iii) \}.
\end{align*}
Here $\overline\xi_i(\iii)$ is the atom containing $\iii$. Since the projective system acts on $\R$,
it has one Lyapunov exponent of multiplicity one on every ergodic component. Thus the
constant-multiplicity hypothesis in D.-J.~Feng~\cite[Section~4]{Feng2023} holds also when $\nu$ is
not ergodic. Observe that $\overline\xi_1(\iii)\subseteq\overline\xi_0(\iii)$ for all
$\iii\in\Sigma_Z$. Letting $\overline\cP=\{[i]:i\in\{1,\ldots,2N\}\}$ be the canonical partition of
$\Sigma_Z$, we define the \emph{pointwise $i$-dimensional entropy} as
\begin{equation*}
  \overline{h}_i(\nu,\iii) = \E_\nu(I_\nu(\overline{\cP}\,|\,\widehat{\overline{\xi}}_i)
  \,|\, \overline{\cI})(\iii)
\end{equation*}
for $\nu$-almost every $\iii\in\Sigma_Z$, where $\overline\cI$ is the $\sigma$-algebra of
shift-invariant Borel subsets of $\Sigma_Z$, and the \emph{$i$-dimensional entropy}
\begin{equation*}
  \overline{h}_i(\nu) = \int_{\Sigma_Z} \overline{h}_i(\nu,\iii) \dd\nu(\iii) =
  H_\nu(\overline{\cP}\,|\,\widehat{\overline{\xi}}_i)
\end{equation*}
is the conditional entropy of $\overline{\cP}$ given $\widehat{\overline{\xi}}_i$ defined in
\cref{sec:conditional}. We similarly have $\overline{h}_i(\nu) = \overline{h}_i(\nu,\iii)$ for
$\nu$-almost all $\iii \in \Sigma_Z$ whenever $\nu$ is ergodic. 

Let $\mu_\iii^{\xi_i}$ and $\nu_\iii^{\overline{\xi}_i}$ be the conditional measures supported on
$\xi_i(\iii)$ and $\overline{\xi}_i(\iii)$, respectively, as defined in \cref{sec:conditional}. The
following lemma is proven by D.-J.~Feng \cite[Lemma 4.6]{Feng2023}.

\begin{lemma} \label{thm:feng_entropy_as_limit}
  If $X$ is a planar self-affine set and $\mu \in \MM_\sigma(\Sigma)$ is ergodic with
  $\lambda_2(\A,\mu) < \lambda_1(\A,\mu)$, then for each $i \in \{0,1,2\}$ we have
  \begin{equation*}
    h_i(\mu,\iii) = -\lim_{n \to \infty} \frac{1}{n}\log\mu_\iii^{\xi_i}(\cP_0^{n-1}(\iii))
  \end{equation*}
  for $\mu$-almost all $\iii \in \Sigma'$. Furthermore, if $X$ is a proximal box-like self-affine
  set with rotations and $\nu \in \MM_\sigma(\Sigma_Z)$, then for each $i \in \{0,1\}$ we also have
  \begin{equation*}
    \overline{h}_i(\nu,\iii) = -\lim_{n \to \infty}
    \frac{1}{n}\log\nu_\iii^{\overline{\xi}_i}(\overline{\cP}_0^{n-1}(\iii))
  \end{equation*}
  for $\nu$-almost all $\iii \in \Sigma_Z$. 
\end{lemma}

Recalling Bowen \cite[\S 2]{Bowen2008}, the \emph{Kolmogorov--Sinai entropy} of $\mu$ is
\begin{align*}
  h(\mu) 
  &= \sup\biggl\{ \lim_{n \to
    \infty}\frac{1}{n}H_\mu\biggl(\bigvee_{i=0}^{n-1}\sigma^{-i}(\cQ)\biggr) : \cQ
  \text{ is a finite Borel partition of }\Sigma\biggr\} \\ 
  &= \lim_{n \to \infty}\frac{1}{n}H_\mu(\cP_0^{n-1}).
\end{align*}
By the Shannon--McMillan--Breiman Theorem (see, e.g., Petersen \cite[\S3, Theorem 2.3]{Petersen1983}), we have
\begin{equation*}
  h(\mu) = -\lim_{n \to \infty} \frac{1}{n}\log\mu(\cP_0^{n-1}(\iii))
\end{equation*}
for $\mu$-almost all $\iii \in \Sigma'$ whenever $\mu$ is ergodic. D.-J.~Feng \cite[Lemma 4.6]{Feng2023}
also gives the following lemma.

\begin{lemma} \label{thm:feng_entropy_kolmogorov_sinai}
  If $X$ is a planar self-affine set and $\mu \in \MM_\sigma(\Sigma)$ is ergodic with
  $\lambda_2(\A,\mu) < \lambda_1(\A,\mu)$, then we have
  \begin{equation*}
    h_0(\mu,\iii) = -\lim_{n \to \infty} \frac{1}{n}\log\mu(\cP_0^{n-1}(\iii))
  \end{equation*}
  for $\mu$-almost all $\iii \in \Sigma'$. 
\end{lemma}

The lemma identifies $h_0(\mu)$ with the Kolmogorov--Sinai entropy $h(\mu)$. More generally, we use
$h_0(\mu)=h(\mu)$ and $\overline h_0(\nu)=h(\nu)$ for all invariant measures on the respective
shifts. Let $\delta_\iii$ be the Dirac measure at $\iii$. The following result strengthens the
vanishing criterion in D.-J.~Feng~\cite[Corollary~4.7]{Feng2023} to a pointwise equivalence.

\begin{lemma} \label{thm:feng_entropy_zero}
  If $X$ is a planar self-affine set and $\mu\in\MM_\sigma(\Sigma)$ is ergodic with
  $\lambda_2(\A,\mu)<\lambda_1(\A,\mu)$, then, for each $i\in\{0,1,2\}$,
  \begin{equation*}
    h_i(\mu,\iii)=0\quad\Leftrightarrow\quad\mu_\iii^{\xi_i}=\delta_\iii
  \end{equation*}
  for $\mu$-almost every $\iii\in\Sigma'$. If $X$ is a proximal box-like self-affine set
  with rotations and $\nu\in\MM_\sigma(\Sigma_Z)$, then, for each $i\in\{0,1\}$,
  \begin{equation*}
    \overline h_i(\nu,\iii)=0\quad\Leftrightarrow\quad
    \nu_\iii^{\overline\xi_i}=\delta_\iii
  \end{equation*}
  for $\nu$-almost every $\iii\in\Sigma_Z$.
\end{lemma}

\begin{proof}
  We first prove the equivalence for $\nu$ on $\Sigma_Z$. Fix $i\in\{0,1\}$, and write
  $f=I_\nu(\overline\cP\,|\,\widehat{\overline\xi}_i)$ and $H=\E_\nu(f\,|\,\overline\cI)=\overline
  h_i(\nu,\,\cdot\,)$. We have $f\ge0$ and $f\in L^1(\nu)$. D.-J.~Feng~\cite[Lemma~4.6]{Feng2023}
  gives, on a common full-measure set,
  \begin{equation*}
    -\log\nu_\iii^{\overline\xi_i}(\overline\cP_0^{n-1}(\iii)) =\sum_{k=0}^{n-1}f(\sigma^k(\iii))
  \end{equation*}
  for every $n\ge1$. The set $E=\{H=0\}$ is invariant modulo $\nu$, and
  $\int_E f\dd\nu=\int_E H\dd\nu=0$. Thus, for $\nu$-almost every $\iii\in E$, we have
  $f(\sigma^k(\iii))=0$ for every $k\ge0$. The displayed identity shows that
  $\nu_\iii^{\overline\xi_i}(\overline\cP_0^{n-1}(\iii))=1$ for every $n$. Since the
  conditional measure is supported on an atom with fixed past, intersecting these cylinders
  gives $\nu_\iii^{\overline\xi_i}(\{\iii\})=1$. Conversely, a Dirac conditional measure
  gives zero in the limit in \cref{thm:feng_entropy_as_limit}, so $H(\iii)=0$.
  The same argument proves the equivalence for $\mu$.
\end{proof}

Together with \cref{thm:finer_partition_dirac_mass}, this criterion relates vanishing
entropies for nested partitions.

To finish this section, let us show that the partitions $\xi_1$ of $\Sigma'$ and $\overline{\xi}_1$
of $\Sigma_Z$ are connected. Recall from \cref{eq:sigma_tau_omega} that
\begin{align*}
  \Sigma^\tau &= \{\iii \in \Sigma' : V_\iii^1 \text{ is the $y$-axis}\}, \\ 
  \Sigma^\omega &= \{\iii \in \Sigma' : V_\iii^1 \text{ is the $x$-axis}\}.
\end{align*}
Recall from \cref{sec:subshift} the definitions of the injective maps $\tau, \omega \colon \Sigma
\to \Sigma_Z$ and their inverses $\overline{\tau}$ and $\overline{\omega}$ defined on $\tau(\Sigma)$
and $\omega(\Sigma)$, respectively. By \cref{eq:tau_omega_partition}, we have
\begin{equation*}
  \Sigma_Z = \tau(\Sigma) \cup \omega(\Sigma) \quad\text{and}\quad \tau(\Sigma) \cap \omega(\Sigma)
  = \emptyset.
\end{equation*}
Recall also that, by \cref{eq:commutative_diagram},
\begin{align*}
  \overline{\pi}(\tau(\iii)) &= \proj_x(\pi(\iii)), \\ 
  \overline{\pi}(\omega(\iii)) &= \proj_y(\pi(\iii))
\end{align*}
for all $\iii\in\Sigma$.

\begin{lemma} \label{thm:partitions_xi1}
  If $X$ is a box-like self-affine set with rotations and
  $\mu\in\MM_\sigma(\Sigma)$ is ergodic with $\lambda_2(\A,\mu)<\lambda_1(\A,\mu)$, then
  \begin{equation*}
    \xi_1(\iii) =
    \begin{cases}
      \overline{\tau}(\overline{\xi}_1(\tau(\iii))) \subseteq \Sigma^\tau, &\text{if } \iii \in
      \Sigma^\tau, \\ 
      \overline{\omega}(\overline{\xi}_1(\omega(\iii))) \subseteq \Sigma^\omega, &\text{if } \iii
      \in \Sigma^\omega,
    \end{cases}
  \end{equation*}
  and
  \begin{equation*}
    \xi^\bot_1(\iii) =
    \begin{cases}
      \overline{\omega}(\overline{\xi}_1(\omega(\iii))) \subseteq \Sigma^\tau, &\text{if } \iii \in
      \Sigma^\tau, \\ 
      \overline{\tau}(\overline{\xi}_1(\tau(\iii))) \subseteq \Sigma^\omega, &\text{if } \iii \in
      \Sigma^\omega.
    \end{cases}
  \end{equation*}
  for every $\iii\in\Sigma'$.
\end{lemma}

\begin{proof}
  We may assume that $\iii \in \Sigma^\tau$ since the proof is symmetric if $\iii \in
  \Sigma^\omega$. The transition rule for $Z$ shows that the past determines whether a sequence
  belongs to $\tau(\Sigma)$ or $\omega(\Sigma)$, so the inverse maps below are defined on the
  whole partition atoms. Let us first show that $\xi_1(\iii) =
  \overline{\tau}(\overline{\xi}_1(\tau(\iii)))$.
  To prove $\overline{\tau}(\overline{\xi}_1(\tau(\iii))) \subseteq \xi_1(\iii)$, let $\jjj \in
  \overline{\tau}(\overline{\xi}_1(\tau(\iii)))$. Since $\tau(\jjj) \in
  \overline{\xi}_1(\tau(\iii))$, we have $\tau(\jjj)^- = \tau(\iii)^-$ and
  $\overline{\pi}(\tau(\jjj)) = \overline{\pi}(\tau(\iii))$. Hence, $\jjj^- = \iii^-$ and
  $\proj_x(\pi(\jjj)) = \proj_x(\pi(\iii))$.
  By our choice of $\Sigma'$ in \cref{sec:lyapunov_exponents}, we have $\jjj\in\Sigma'$ and
  $V_\jjj^1=V_\iii^1$, so $\jjj\in\Sigma^\tau$. Since $V_\iii^1$ is the $y$-axis, we have
  $\pi(\jjj)-\pi(\iii) \in V_\iii^1$ and hence $\jjj \in \xi_1(\iii)$.
  Since the reasoning reverses, we have shown the claim.

  Let us then prove that $\xi^\bot_1(\iii) = \overline{\omega}(\overline{\xi}_1(\omega(\iii)))$. To
  show that $\overline{\omega}(\overline{\xi}_1(\omega(\iii))) \subseteq \xi^\bot_1(\iii)$, let
  $\jjj \in \overline{\omega}(\overline{\xi}_1(\omega(\iii)))$. Since $\omega(\jjj) \in
  \overline{\xi}_1(\omega(\iii))$, we have $\omega(\jjj)^- = \omega(\iii)^-$ and
  $\overline{\pi}(\omega(\jjj)) = \overline{\pi}(\omega(\iii))$. Therefore, $\jjj^- = \iii^-$ and
  $\proj_y(\pi(\jjj)) = \proj_y(\pi(\iii))$.
  As above, $\jjj\in\Sigma^\tau$, and $\pi(\jjj)-\pi(\iii) \,\bot\, V_\iii^1$, so
  $\jjj \in \xi^\bot_1(\iii)$. Since the reasoning reverses, we have finished the proof.
\end{proof}

The proof also shows that $\Sigma^\tau$ and $\Sigma^\omega$ are unions of past atoms, and
that each of $\tau(\Sigma^\tau)$, $\tau(\Sigma^\omega)$, $\omega(\Sigma^\tau)$ and
$\omega(\Sigma^\omega)$ is a union of atoms of $\overline\xi_0$, hence also of
$\overline\xi_1$. The following technical lemma will help us compare the entropies of ergodic
measures on $\Sigma$ and $\Sigma_Z$.

\begin{lemma} \label{thm:partition_tau_conj}
  If $X$ is a box-like self-affine set with rotations,
  $\mu\in\MM_\sigma(\Sigma)$ is ergodic with $\lambda_2(\A,\mu)<\lambda_1(\A,\mu)$, and
  $\iii\in\Sigma'$, then
  \begin{equation*}
    (\xi_1 \lor \cP_0^{n-1})(\iii) = 
    \begin{cases}
      \overline{\tau}((\overline{\xi}_1 \lor \overline{\cP}_0^{n-1})(\tau(\iii))), &\text{if } \iii
      \in \Sigma^\tau, \\ 
      \overline{\omega}((\overline{\xi}_1 \lor \overline{\cP}_0^{n-1})(\omega(\iii))), &\text{if }
      \iii \in \Sigma^\omega,
    \end{cases}
  \end{equation*}
  for every $n\in\N^+$.
\end{lemma}

\begin{proof}
  We prove the claim for $\iii = (i_k)_{k \in \Z} \in \Sigma^\tau$ and note that the proof in the
  case $\iii \in \Sigma^\omega$ is analogous. By \cref{eq:def_of_tau_omega}, for any
  $\jjj\in\Sigma$ and $n\ge1$, the equalities $i_k=j_k$ for $0\le k<n$ hold if and only if
  $\tau_k(\iii)=\tau_k(\jjj)$ for $0\le k<n$. Consequently,
  \begin{equation*}
    \overline\tau\bigl(\overline\cP_0^{n-1}(\tau(\iii))\cap\tau(\Sigma)\bigr)
    =\cP_0^{n-1}(\iii).
  \end{equation*}
  Intersecting this identity with the corresponding partition atom in
  \cref{thm:partitions_xi1} proves the assertion for every $n$.
\end{proof}

\subsection{Affinity dimension and equilibrium state} \label{sec:affinity_dim}
Let $\Sigma_*^+ = \bigcup_{n \in \N^+} \Sigma_n^+$, where $\Sigma_n^+ = \{1,\ldots,N\}^n$
for all $n \in \N^+$, be the collection of finite words.
For each $\iii \in \Sigma_*^+$, we let the \emph{length} $|\iii|$ be
the unique $n \in \N^+$ for which $\iii \in \Sigma_n^+$. If $\mathsf{A} = (A_1,\ldots,A_N)
\in \GL^N$, then we write $A_\iii = A_{i_0} \cdots A_{i_{n-1}}$ for all $\iii = i_0 \cdots i_{n-1}
\in \Sigma_n^+$ and $n \in \N^+$. For each $\A \in \GL^N$ and $s \ge 0$ we define the
\emph{singular value pressure} by setting
\begin{equation*}
  P(\A,s) = \lim_{n \to \infty} \frac{1}{n}\log\sum_{\iii \in \Sigma_n^+} \fii^s(A_\iii),
\end{equation*}
where $\fii^s$ is the singular value function defined in \cref{sec:Singular_value_function}.
As the singular value function is sub-multiplicative, the sequence $(\log\sum_{\iii \in \Sigma_n^+}
\fii^s(A_\iii))_{n \in \N}$ is sub-additive and hence, the limit above exists by Fekete's lemma. It
is also easy to see that the pressure function $s \mapsto P(\A,s)$ is continuous and convex between
consecutive integers. If $\max_{i \in \{1,\ldots,N\}}\|A_i\|<1$, then it is also strictly decreasing
with $P(\A,0) \ge 0$ and $\lim_{s \to \infty}P(\A,s) = -\infty$, see K\"aenm\"aki--Vilppolainen
\cite[Lemma 2.1]{KaenmakiVilppolainen2010}. In this case, we may thus define the \emph{affinity
dimension} by setting $\dimaff(\A)$ to be the unique $s \ge 0$ for which $P(\A,s)=0$.

We are interested in understanding the geometry of self-affine sets. Relying on
\cref{eq:self_affine_set_def}, the self-affine set $X$ can naturally be covered by the sets
$\fii_\iii(B)$, where $\fii_\iii = \fii_{i_0} \circ \cdots \circ \fii_{i_{n-1}}$ for all $\iii = i_0
\cdots i_{n-1} \in \Sigma_n^+$ and $B$ is a (closed) ball containing $X$. Observe that each ellipse
$\fii_\iii(B)$ can be covered by one ball of radius $\alpha_1(A_\iii)\diam(B)$ or by
$\lceil \alpha_1(A_\iii)/\alpha_2(A_\iii) \rceil$ many balls of radius
$\alpha_2(A_\iii)\diam(B)$. This motivates
us to study the limiting behaviour of the sums $\sum_{\iii \in \Sigma_n^+} \fii^s(A_\iii)$ and
indeed, it is straightforward to see that
\begin{equation} \label{eq:hausdorff_affinity_upper_bound}
  \dimh(X) \le \min\{2,\dimaff(\A)\},
\end{equation}
see Falconer \cite[Theorem 5.4]{Falconer1988}. While the affinity dimension only depends on $\A$, we
may write $\dimaff(X)=\dimaff(\A)$ for the associated attractor.

For the definitions of the Minkowski (or box-counting) dimension $\dimm$ and the packing dimension
$\dimp$ appearing in \cref{thm:fraser}, we refer to Falconer \cite[\S 2]{Falconer1997}. If $X$ is a
box-like self-affine set, then the singular value function admits a modification, introduced by
Fraser \cite[\S 2.1]{Fraser2012}, which takes the coordinate projections of $X$ into account. In
this case, each $A_\iii$ for $\iii \in \Sigma_*^+$ is either diagonal or anti-diagonal (recall the
proof of \cref{thm:box_like_sa}) and hence, $A_\iii([-1,1]^2)$ is a rectangle with sides parallel
to the coordinate axes. The modified singular value function weights the longer side of this
rectangle by the Minkowski dimension of the projection of $X$ onto the coordinate axis parallel
to it, and the \emph{modified affinity dimension} $\dimmaff(X)$ is the zero of the associated
pressure. We refer to \cite[\S 2.1]{Fraser2012} for the precise definition. It satisfies
$\udimm(X) \le \dimmaff(X) \le \dimaff(X)$, with the first inequality given by \cref{thm:fraser}.
The second inequality follows by comparing the two pressures at $s = \dimaff(X)$. The modified
singular value function is bounded above by the singular value function whenever $s$ is at least
the Minkowski dimensions of both coordinate projections. These dimensions do not exceed
$\udimm(X)$, which is at most $\dimaff(X)$ by the proof of \cite[Theorem~5.4]{Falconer1988}.
Under the rectangular open set condition, it agrees with the Minkowski and packing dimensions of
$X$. We note that, in contrast to the singular value function, the modified singular value
function depends not only on $\A$ but, through the Minkowski dimensions of the coordinate
projections, also on the set $X$ itself. In particular, the modified affinity dimension is not
determined by $\A$ alone, and there is no analogue of the convention
$\dimaff(X) = \dimaff(\A)$.

Write $\iii|_n = i_0\cdots i_{n-1}$ for all $\iii = (i_n)_{n \in \Z} \in \Sigma$ and $n \in
\N^+$. If $\mu \in \MM_\sigma(\Sigma)$ and $\mathsf{A} = (A_1,\ldots,A_N) \in \GL^N$, then we
define the \emph{$s$-energy} of $\mu$ by setting
\begin{equation} \label{eq:energy_def}
  \begin{split}
    \Lambda(\A,\mu,s) 
    &= \lim_{n \to \infty} \frac{1}{n} \int_{\Sigma} \log\fii^s(A_{\iii|_n})
    \dd\mu(\iii) \\
    &=
    \begin{cases}
      s\lambda_1(\A,\mu), &\text{if } 0 \le s \le 1, \\ 
      \lambda_1(\A,\mu) + (s-1)\lambda_2(\A,\mu), &\text{if } 1 < s \le 2, \\ 
      \tfrac{s}{2}(\lambda_1(\A,\mu)+\lambda_2(\A,\mu)), &\text{if } s > 2.
    \end{cases}
  \end{split}
\end{equation}
Since the singular value function is sub-multiplicative and $\mu$ is invariant, the limit above
exists by Fekete's lemma. An application of Jensen's inequality yields
\begin{equation} \label{eq:variational}
  P(\A,s) \ge h_0(\mu) + \Lambda(\A,\mu,s)
\end{equation}
for all $\mu \in \MM_\sigma(\Sigma)$ and $s \ge 0$. A measure $\mu \in \MM_\sigma(\Sigma)$ is called
an \emph{$s$-equilibrium state} for $\mathsf{A}$ if it satisfies
\begin{equation} \label{eq:eq_state_def}
  P(\A,s) = h_0(\mu) + \Lambda(\A,\mu,s).
\end{equation}
K\"{a}enm\"{a}ki \cite[Theorems 2.6 and 4.1]{Kaenmaki2004} proved that for each $\mathsf{A} \in
\GL^N$ and $s \ge 0$ there exists an ergodic $s$-equilibrium state for $\mathsf{A}$. The example of
K\"aenm\"aki--Vilppolainen \cite[Example 6.2]{KaenmakiVilppolainen2010} shows that such an
equilibrium state is not necessarily unique, and the result of Feng--K\"aenm\"aki \cite[Theorem
1.7]{FengKaenmaki2011} shows that, in general, there are at most two ergodic $s$-equilibrium states
for $\A$. Moreover, by the following proposition, which is \cite[Proposition
1.2]{FengKaenmaki2011}, only reducible tuples $\A$ can have two $s$-equilibrium states.

\begin{proposition}[Feng--K\"aenm\"aki] \label{thm:irred_eq_state}
  If $\A = (A_1,\ldots,A_N) \in \GL^N$ is irreducible, then for each $s>0$ there is a unique
  $s$-equilibrium state $\mu \in \MM_\sigma(\Sigma)$. This measure is ergodic, and there exists
  a constant $C=C(\A,s)\ge1$ such that
  \begin{equation*}
    C^{-1}\exp(-|\iii|P(\A,s)) \fii^s(A_\iii) \le \mu([\iii]) \le C\exp(-|\iii|P(\A,s)) \fii^s(A_\iii)
  \end{equation*}
  for all $\iii \in \Sigma_*^+$.
\end{proposition}

Feng--K\"aenm\"aki state their proposition for the potentials $\iii \mapsto \|A_\iii\|^q$ with
$q>0$, which is $\fii^s(A_\iii)$ when $q=s\le1$. For $1<s<2$ one applies it to the rescaled tuple
$\A' = (A_1',\ldots,A_N')$, where $A_i' = |\det A_i|^{(s-1)/(2-s)}A_i$, with $q=2-s$, since
$\|A_\iii'\|^{2-s} = \fii^s(A_\iii)$ for all $\iii \in \Sigma_*^+$ and $\A'$ is irreducible
whenever $\A$ is.

Let $\iii\jjj$ denote the concatenation of two finite words $\iii$ and $\jjj$. We say that a Borel
probability measure $\mu$ is \emph{sub-multiplicative} if there is $C \ge 1$ such that
$\mu([\iii\jjj]) \le C\mu([\iii])\mu([\jjj])$ for all $\iii, \jjj \in \Sigma_*^+$. It follows from
\cref{thm:irred_eq_state} and the sub-multiplicativity of the singular value function that the
$s$-equilibrium state $\mu$ satisfies
\begin{equation*}
  \mu([\iii\jjj]) \le C\exp(-(|\iii|+|\jjj|)P(\A,s)) \fii^s(A_\iii)\fii^s(A_\jjj) \le
  C^3\mu([\iii])\mu([\jjj])
\end{equation*}
for all $\iii, \jjj \in \Sigma_*^+$ and hence $\mu$ is sub-multiplicative. Note that if $s =
\dimaff(\A)$, then the $s$-equilibrium state $\mu \in \MM_\sigma(\Sigma)$ gives each cylinder a mass
which is approximately the value of the corresponding singular value function, that is,
\begin{equation*}
  C^{-1} \fii^s(A_\iii) \le \mu([\iii]) \le C \fii^s(A_\iii)
\end{equation*}
for all $\iii \in \Sigma_*^+$. If $\mu \in \MM_\sigma(\Sigma)$ and $\A = (A_1,\ldots,A_N) \in \GL^N$
is such that $\max_{i \in \{1,\ldots,N\}}\|A_i\|<1$, then the unique nonnegative real number $s$ for which
$h_0(\mu) + \Lambda(\A,\mu,s) = 0$ is called the \emph{Lyapunov dimension} of $\mu$, denoted by
$\diml(\A,\mu)$. By \cref{eq:energy_def}, we have
\begin{equation} \label{eq:diml_min2}
  \diml(\A,\mu) =
  \begin{dcases}
    -\frac{h_0(\mu)}{\lambda_1(\A,\mu)},
    &\text{if } 0 \le h_0(\mu) \le -\lambda_1(\A,\mu), \\
    1-\frac{h_0(\mu)+\lambda_1(\A,\mu)}{\lambda_2(\A,\mu)},
    &\text{if } 0 < h_0(\mu) + \lambda_1(\A,\mu) \le -\lambda_2(\A,\mu), \\
    -\frac{2h_0(\mu)}{\lambda_1(\A,\mu)+\lambda_2(\A,\mu)},
    &\text{if } h_0(\mu) > -\lambda_1(\A,\mu)-\lambda_2(\A,\mu).
  \end{dcases}
\end{equation}

Recalling \cref{eq:variational}, we see that $\diml(\A,\mu) \le \dimaff(\A)$ for all $\mu \in
\MM_\sigma(\Sigma)$. If $s=\dimaff(\A)$ and $\mu \in \MM_\sigma(\Sigma)$ is an $s$-equilibrium state
for $\A$, then, by \cref{eq:eq_state_def}, we have $\diml(\A,\mu) = \dimaff(\A)$. 

\subsection{Normalisation and simple Lyapunov spectrum} \label{sec:normalisation}
We now justify the normalisation described in \cref{sec:irreducibility}. Let $\Phi$ be a
proximal box-like system with rotations in fixed coordinates, and let $\Phi^{(2)}$ and
$\A^2$ be its second iterates and their linear parts. Write
$\Sigma^{(2)}=(\{1,\ldots,N\}^2)^\Z$, and denote its left shift and canonical projection by
$\sigma^{(2)}$ and $\pi^{(2)}$, respectively. The grouping map
$\beta\colon\Sigma\to\Sigma^{(2)}$ is defined by
\begin{equation*}
  \beta((i_n)_{n\in\Z})=((i_{2n},i_{2n+1}))_{n\in\Z}.
\end{equation*}
For $s>0$, let $\mu$ be the $s$-equilibrium state of $\A$ and write
$\mu^{(2)}=\beta_*\mu$. Let $\Psi$ and $\Psi^{(2)}$ be the projective systems of $\Phi$ and
$\Phi^{(2)}$, respectively, as defined in \cref{sec:projective_system}.

\begin{lemma} \label{thm:two-block-normalisation}
  The system $\Phi^{(2)}$ has attractor $X$, and $\A^2$ contains both matrix types and a proximal
  diagonal generator. The strong separation, strong open set and projective exponential separation
  conditions each pass from $\Phi$ to $\Phi^{(2)}$. Furthermore, $\mu^{(2)}$ is the unique ergodic
  $s$-equilibrium state of $\A^2$, and
  \begin{equation*}
    P(\A^2,t) = 2P(\A,t), \quad h_0(\mu^{(2)}) = 2h_0(\mu),
    \quad\text{and}\quad \lambda_i(\A^2,\mu^{(2)}) = 2\lambda_i(\A,\mu)
  \end{equation*}
  for all $t\ge0$ and $i\in\{1,2\}$. Moreover, $\pi^{(2)}_*\mu^{(2)}=\pi_*\mu$.
\end{lemma}

\begin{proof}
  The identity $X=\bigcup_{i,j}\fii_i\circ\fii_j(X)$ gives the attractor. Products of two
  diagonal generators and of a diagonal and an anti-diagonal generator show that both matrix
  types occur. By \cref{thm:box_like_proximal}, either $\A^2$ already contains a proximal
  diagonal generator or $\A$ contains one, whose square belongs to $\A^2$. Hence $\A^2$ is
  irreducible: a common invariant line would have to be an axis, whereas an anti-diagonal
  generator exchanges the axes. Disjointness of the first-level images passes to the
  second-level images by injectivity. A feasible open set for $\Phi$ is also feasible for
  $\Phi^{(2)}$, so both planar separation conditions are preserved.

  The map $\beta$ is a homeomorphism satisfying
  $\beta\circ\sigma^2=\sigma^{(2)}\circ\beta$ and $\pi^{(2)}\circ\beta=\pi$. A word of length
  $n$ over the paired alphabet corresponds bijectively to a word of length $2n$ over the
  original alphabet, with the same matrix product. The pressure and Lyapunov identities
  follow by taking the even subsequences in their defining limits. The same correspondence
  identifies the block entropies of $\mu^{(2)}$ at length $n$ with those of $\mu$ at length
  $2n$, and therefore gives $h_0(\mu^{(2)})=2h_0(\mu)$. Hence,
  \begin{equation*}
    h_0(\mu^{(2)})+\Lambda(\A^2,\mu^{(2)},s)
      =2\bigl(h_0(\mu)+\Lambda(\A,\mu,s)\bigr)=P(\A^2,s).
  \end{equation*}
  Thus $\mu^{(2)}$ is an equilibrium state. Irreducibility of $\A^2$ and
  \cref{thm:irred_eq_state} give uniqueness and ergodicity.

  To prove preservation of projective separation, put $\eps_i=0$ for a diagonal generator and $\eps_i=1$ for an anti-diagonal generator. For $q\in\{0,1\}$, write $g_{i,q}=g_{i+qN}$. By \cref{eq:phi_def_by_g}, the projective map of $\fii_i\circ\fii_j$ corresponding to coordinate $q$ is
  \begin{equation*}
    g^{(2)}_{(i,j),q}=g_{i,q}\circ g_{j,(q+\eps_i)\bmod2}.
  \end{equation*}
  The index $((i,j),q)$ determines a distinct two-letter word in the alphabet of $\Psi$.
  Consequently, distinct words of length $n$ for $\Psi^{(2)}$ give distinct words of length
  $2n$ for $\Psi$, and $\Delta_n^{\Psi^{(2)}}\ge\Delta_{2n}^\Psi$.

  By \cref{sec:exponential_separation}, we have $\Delta_{m+1}^\Psi\le\Delta_m^\Psi$ for every $m$.
  Suppose $\Delta_{n_k}^\Psi\ge c^{n_k}$ along an increasing sequence, where we may take $0<c<1$.
  Set $m_k=\lfloor n_k/2\rfloor$ and pass to an increasing subsequence with $m_k\ge1$. Then
  \begin{equation*}
    \Delta_{m_k}^{\Psi^{(2)}}\ge\Delta_{2m_k}^\Psi
      \ge\Delta_{n_k}^\Psi\ge c^{n_k}
      \ge c^{2m_k+1}\ge(c^3)^{m_k}
  \end{equation*}
  as required.
\end{proof}

Thus the affinity dimension, Lyapunov dimension of the equilibrium
state, and the projected planar measure itself are preserved. Hence both coordinate projections
are unchanged by the grouping.
Similarly, invertible linear changes of coordinates $T$ preserve pressure, Lyapunov exponents
and the dimensions of planar sets and measures.

We say that a Borel probability measure $\mu$ on $\Sigma$ is \emph{Bernoulli} if there is a
probability vector $(p_1,\ldots,p_N)$ such that
\begin{equation*}
  \mu([i_k\cdots i_\ell]) = p_{i_k}\cdots p_{i_\ell}
\end{equation*}
for all $i_j \in \{1,\ldots,N\}$ and all $k \le \ell$ in $\Z$. Since the condition is imposed at
every position, a Bernoulli measure is $\sigma$-invariant, and it is well-known that it is ergodic,
see e.g.\ Walters~\cite{Walters1982}. Now fix $0<s<2$, suppose that $\A$ is proximal and box-like
with rotations, and let $\mu$ be its equilibrium state. If $\A$ contains a proximal diagonal
generator, Morris--Shmerkin~\cite[Lemma~3.5]{MorrisShmerkin2019} shows that $\mu$ is not
Bernoulli, and \cite[Lemma~3.4]{MorrisShmerkin2019} then gives distinct Lyapunov exponents.
Otherwise we apply these results to $\A^2$ and $\mu^{(2)}$ from
\cref{thm:two-block-normalisation}. Grouping a Bernoulli measure in pairs gives a Bernoulli
measure, so $\mu^{(2)}$ being non-Bernoulli implies the same for $\mu$. The scaling of
Lyapunov exponents in that lemma also transfers their strict inequality to $\mu$. Thus
\begin{equation} \label{eq:simple_lyapunov}
  \lambda_2(\A,\mu) < \lambda_1(\A,\mu)
\end{equation}
and the assumptions used to define the Oseledets partitions in
\cref{sec:lyapunov_exponents,sec:entropies} hold for these equilibrium states.

\begin{lemma} \label{thm:bernoulli_equilibrium}
  For $s\ge2$, we have $P(\A,s)=\log\sum_{i=1}^N|\det A_i|^{s/2}$, and the unique
  $s$-equilibrium state is Bernoulli with weights
  \begin{equation*}
    p_i=\frac{|\det A_i|^{s/2}}{\sum_{j=1}^N|\det A_j|^{s/2}}.
  \end{equation*}
\end{lemma}

\begin{proof}
  The pressure formula follows from $\fii^s(A_\iii)=|\det A_\iii|^{s/2}$ and multiplicativity
  of the determinant. If $\nu\in\MM_\sigma(\Sigma)$ has one-symbol weights $q_i=\nu([i])$,
  then its entropy rate is at most $-\sum_i q_i\log q_i$, with equality precisely when $\nu$
  is Bernoulli. Hence Gibbs' inequality gives
  \begin{equation*}
    h_0(\nu)+\Lambda(\A,\nu,s)
    \le -\sum_{i=1}^N q_i\log q_i+\frac{s}{2}\sum_{i=1}^N q_i\log|\det A_i|
    \le P(\A,s),
  \end{equation*}
  with equality throughout precisely when $\nu$ is Bernoulli with weights $q_i=p_i$.
\end{proof}

Suppose that $\A \in \GL^N$ is proximal and box-like with rotations and let $\mu \in
\MM_\sigma(\Sigma)$ be the unique ergodic $s$-equilibrium state with $0<s<2$.
Recalling \cref{eq:simple_lyapunov}, let $\Sigma' \subseteq \Sigma$ be the set of full $\mu$-measure
as in \cref{sec:lyapunov_exponents}. Recall from \cref{eq:sigma_tau_omega} that
\begin{align*}
  \Sigma^\tau &= \{\iii \in \Sigma' : V_\iii^1 \text{ is the $y$-axis}\}, \\
  \Sigma^\omega &= \{\iii \in \Sigma' : V_\iii^1 \text{ is the $x$-axis}\}.
\end{align*}
Although $\mu(\Sigma^\tau) + \mu(\Sigma^\omega) = 1$ by \cref{eq:tau_omega_partition2}, the
following lemma guarantees that both sets are non-negligible with respect to $\mu$.

\begin{lemma} \label{thm:sigma_tau_omega_positive_eq_measure}
  If $X$ is a proximal box-like self-affine set with rotations and $\mu \in \MM_\sigma(\Sigma)$ is
  the unique ergodic $s$-equilibrium state with $0<s<2$, then
  \begin{equation*}
    \mu(\Sigma^\tau) > 0 \quad\text{and}\quad \mu(\Sigma^\omega) > 0.
  \end{equation*}
\end{lemma}

\begin{proof}
  Recall that $\ell$ is the largest integer such that $A_i$ is a diagonal matrix for all
  $1\leq i\leq \ell$.
  We first show 
  \begin{equation} \label{eq:sigma_tau_omega_positive_eq_measure1}
    \sigma^{k-1}(\iii) \in \bigcup_{j_0=\ell+1}^N [j_0] \quad\Rightarrow\quad \sigma^{k}(\iii)
    \in \Sigma^\tau \text{ and } \sigma^{k-1}(\iii) \in \Sigma^\omega \text{ or vice versa}
  \end{equation}
  for all $\iii \in \Sigma'$ and $k \in \N$. Notice first that if $\iii = (i_n)_{n \in \Z}$ is such
  that $\sigma^{k-1}(\iii) \in \bigcup_{j_0=\ell+1}^N [j_0]$, then $i_{k-1} \ge \ell+1$.
  Thus, $A_{i_{k-1}}$ is anti-diagonal.
  By \cref{eq:V_lyapunov},
  \begin{align*}
    V_{\sigma^k(\iii)}^1 
      &= \{0\} \cup \biggl\{ v \in \R^2 \setminus \{0\} :
	\lim_{n\to\infty}\frac{1}{n} \log| A_{i_{k-n}} \cdots A_{i_{k-2}} A_{i_{k-1}} v| =
      \lambda_{2}(\A,\mu) \biggr\}, \\ 
      V_{\sigma^{k-1}(\iii)}^1 
      &= \{0\} \cup \biggl\{ v \in \R^2 \setminus \{0\} :
	\lim_{n\to\infty}\frac{1}{n} \log| A_{i_{k-n}} \cdots A_{i_{k-2}} v| =
      \lambda_{2}(\A,\mu) \biggr\}.
  \end{align*}
  Since $A_{i_{k-1}}v\bot v$ for every nonzero vector $v$ on either coordinate axis, and all
  matrices in the product are diagonal or anti-diagonal, we see that $i_{k-1} \ge \ell+1$ implies
  $V_{\sigma^k(\iii)}^1 \bot V_{\sigma^{k-1}(\iii)}^1$ giving
  \cref{eq:sigma_tau_omega_positive_eq_measure1}.

  Since $\Sigma^\tau$ and $\Sigma^\omega$ are disjoint, \cref{eq:sigma_tau_omega_positive_eq_measure1}
  places exactly one of $\sigma^k(\iii)$ and $\sigma^{k+1}(\iii)$ in $\Sigma^\tau$ whenever
  $\sigma^k(\iii) \in \bigcup_{j_0=\ell+1}^N [j_0]$, so that
  \begin{equation*}
    \mathds{1}_{\Sigma^\tau}(\sigma^k(\iii)) + \mathds{1}_{\Sigma^\tau}(\sigma^{k+1}(\iii)) \ge
    \mathds{1}_{\bigcup_{j_0=\ell+1}^N [j_0]}(\sigma^k(\iii))
  \end{equation*}
  for all $\iii \in \Sigma'$ and $k \in \N$. Summing over $k \in \{0,\ldots,n-1\}$, each $k$ with
  $\sigma^k(\iii) \in \Sigma^\tau$ is counted at most twice on the left-hand side.

  Now by \cref{thm:irred_eq_state}, the $s$-equilibrium state $\mu$ is ergodic and fully
  supported. Hence, Birkhoff's pointwise ergodic theorem gives
  \begin{equation}
    \mu(\Sigma^\tau) 
	  = \lim_{n \to \infty} \frac{1}{n} \#\{k \in \{0,\ldots,n-1\} : \sigma^k(\iii)
	  \in \Sigma^\tau\}
  \end{equation}
  for $\mu$-almost all $\iii \in \Sigma'$.
  Combining this limit with the preceding counting inequality gives
  \begin{align*}
  \mu(\Sigma^\tau) 
	&\ge \frac12 \lim_{n \to \infty} \frac{1}{n} \# \biggl\{k \in \{0,\ldots,n-1\} : \sigma^k(\iii)
	\in \bigcup_{j_0=\ell+1}^N [j_0]\biggr\} \\ 
	&= \frac12 \mu\biggl( \bigcup_{j_0=\ell+1}^N [j_0] \biggr) > 0
  \end{align*}
  for $\mu$-almost all $\iii \in \Sigma'$, where the $\tfrac12$ factor was introduced by the argument
  above. Since $\mu(\Sigma^\omega) > 0$ can be shown similarly, we
  have finished the proof.
\end{proof}

\subsection{Gibbs measures} \label{sec:gibbs_measures}
Let $\Sigma_Z$ be the sub-shift of finite type as in \cref{sec:subshift} and let $f \colon
\Sigma_Z \to \R$ be a H\"older continuous potential. Recall that $\iii|_n = i_0 \cdots i_{n-1}$ for
all $\iii = (i_n)_{n \in \Z} \in \Sigma$ and $n \in \N$. The \emph{Birkhoff sum} of $f$ is
$S_nf = \sum_{k=0}^{n-1}f\circ\sigma^k$ and the \emph{pressure} of $f$ is
\begin{equation*}
  P_Z(f) = \lim_{n \to \infty} \frac{1}{n}\log\sum_{\iii \in \Sigma_Z \,:\, \sigma^n(\iii)=\iii}
  \exp\Bigl(\sup_{\jjj \in [\iii|_n]} S_nf(\jjj)\Bigr).
\end{equation*}
Write $\MM_\sigma(\Sigma_Z)$ for the collection of all $\sigma$-invariant Borel probability measures
on $\Sigma_Z$. A measure $\nu \in \MM_\sigma(\Sigma_Z)$ is a \emph{Gibbs measure} if
\begin{equation*}
  C^{-1}\exp(-nP_Z(f)+S_nf(\iii)) \le \nu([\iii|_n]) \le C\exp(-nP_Z(f)+S_nf(\iii))
\end{equation*}
for all $\iii \in \Sigma_Z$ and $n \in \N^+$. Since
\begin{equation*}
  \nu([\iii|_n]) \le C^2\exp(P_Z(f)+\|f\|_\infty) \nu([\iii|_{n+1}])
\end{equation*}
for all $\iii \in \Sigma_Z$ and $n \in \N^+$, and since the invariance of $\nu$ gives the same
bound with the same constant for the cylinders extended also to the past, which are the balls of
$\Sigma_Z$, we see that the Gibbs measure is \emph{doubling}: the
measure of any ball is at most a constant times the measure of the corresponding ball with half the
radius. Recall that, since the sub-shift $\Sigma_Z$ is mixing, for each H\"older continuous
potential $f$ there exists a unique ergodic Gibbs measure 
\cite[Theorem 4.7 and Corollary 4.2]{SarigNotes}.

Let the iterated function system $\Psi = (g_1,\ldots,g_{2N})$ be as in \cref{sec:projective_system}
and fix $0 < s < 2$. In our applications, $s$ is the affinity dimension $\dimaff(\A)$. We define
locally constant potentials $f_{1,s}, f_{2,s} \colon \Sigma_Z \to \R$ by setting
\begin{equation*}
  f_{1,s}(\iii) =
  \begin{cases}
    s\log |g_{i_0}'(0)| = s\log |a_{i_0}|, &\text{if } i_0 \in \{1,\ldots,N\}, \\ 
    s\log |g_{i_0}'(0)| = s\log |b_{i_0-N}|, &\text{if } i_0 \in \{N+1,\ldots,2N\},
  \end{cases}
\end{equation*}
and
\begin{equation*}
  f_{2,s}(\iii) =
  \begin{cases}
    s\log |g_{i_0+N}'(0)| = s\log |b_{i_0}|, &\text{if } i_0 \in \{1,\ldots,N\}, \\ 
    s\log |g_{i_0-N}'(0)| = s\log |a_{i_0-N}|, &\text{if } i_0 \in \{N+1,\ldots,2N\},
  \end{cases}
\end{equation*}
for all $\iii = (i_n)_{n \in \Z} \in \Sigma_Z$ when $0 < s < 1$, and similarly
\begin{align*}
  f_{1,s}(\iii) 
    &=
    \begin{cases}
      \log |g_{i_0}'(0)| + (s-1)\log |g_{i_0+N}'(0)|, &\text{if } i_0 \in \{1,\ldots,N\}, \\ 
      \log |g_{i_0}'(0)| + (s-1)\log |g_{i_0-N}'(0)|, &\text{if } i_0 \in \{N+1,\ldots,2N\},
    \end{cases} \\
    &=
    \begin{cases}
      \log |a_{i_0}| + (s-1)\log |b_{i_0}|, &\text{if } i_0 \in \{1,\ldots,N\}, \\ 
      \log |b_{i_0-N}| + (s-1)\log |a_{i_0-N}|, &\text{if } i_0 \in \{N+1,\ldots,2N\},
    \end{cases}
\end{align*}
and
\begin{align*}
  f_{2,s}(\iii) 
    &=
    \begin{cases}
      \log |g_{i_0+N}'(0)| + (s-1)\log |g_{i_0}'(0)|, &\text{if } i_0 \in \{1,\ldots,N\}, \\ 
      \log |g_{i_0-N}'(0)| + (s-1)\log |g_{i_0}'(0)|, &\text{if } i_0 \in \{N+1,\ldots,2N\},
    \end{cases} \\ 
    &= 
    \begin{cases}
      \log |b_{i_0}| + (s-1)\log |a_{i_0}|, &\text{if } i_0 \in \{1,\ldots,N\}, \\ 
      \log |a_{i_0-N}| + (s-1)\log |b_{i_0-N}|, &\text{if } i_0 \in \{N+1,\ldots,2N\},
    \end{cases}
\end{align*}
for all $\iii = (i_n)_{n \in \Z} \in \Sigma_Z$ when $1 \le s < 2$. Let $\tau, \omega \colon \Sigma
\to \Sigma_Z$ and their inverses $\overline{\tau}$, $\overline{\omega}$
be as in \cref{sec:subshift}. Recalling that $g_{\iii|_n} = g_{i_0} \circ \cdots \circ g_{i_{n-1}}$
for all $\iii = (i_k)_{k \in \Z} \in \Sigma_Z$ and $n \in \N^+$, we have
\begin{equation*}
  \exp S_n f_{1,1}(\tau(\iii)) = |g_{\tau(\iii)|_n}'(0)| \quad\text{ and }\quad
  \exp S_n f_{1,1}(\omega(\iii)) = |g_{\omega(\iii)|_n}'(0)|
\end{equation*}
for all $\iii \in \Sigma$ and hence, by \cref{eq:tau_omega_partition},
\begin{equation} \label{eq:lyapunov_exp_for_psi}
  |g_{\iii|_n}'(0)| = \exp S_n f_{1,1}(\iii)
\end{equation}
for all $\iii \in \Sigma_Z$. It follows that $\fii_{\iii|_n}([-1,1]^2)$ is a rectangle with sides
parallel to the axes such that the horizontal side has length $2\exp(S_nf_{1,1}(\tau(\iii)))$ and
the vertical side has length $2\exp(S_nf_{2,1}(\tau(\iii)))$. Therefore,
\begin{equation*}
  \fii^s(A_{\iii|_n}) = \max\{\exp(S_nf_{1,s}(\tau(\iii))), \exp(S_nf_{2,s}(\tau(\iii)))\}
\end{equation*}
for all $\iii \in \Sigma$ and $n \in \N^+$. Note also that $f_{1,s}(\tau(\iii)) =
f_{2,s}(\omega(\iii))$ and $f_{2,s}(\tau(\iii)) = f_{1,s}(\omega(\iii))$.

Let $\mu_1, \mu_2 \in \MM_\sigma(\Sigma_Z)$ be the Gibbs measures associated to the potentials
$f_{1,s}$ and $f_{2,s}$, respectively. The measures $\mu_1$ and $\mu_2$ are called the
\emph{$s$-Gibbs measures}. If $0 < s \le 1$, then \cref{eq:lyapunov_exp_for_psi} shows that
\begin{equation*}
  C^{-1}\exp(-nP_Z(f_{1,s}))|g_{\iii|_n}'(0)|^s \le \mu_1([\iii|_n]) \le
  C\exp(-nP_Z(f_{1,s}))|g_{\iii|_n}'(0)|^s
\end{equation*}
for all $\iii \in \Sigma_Z$ and $n \in \N$, and similarly for $\mu_2$.

 Notice that, by symmetry, we
have $P_Z(f_{1,s}) = P_Z(f_{2,s})$. Observe also that $\mu_1 \ne \mu_2$. Indeed, by
our normalisation, there is $i_0 \in \{1,\ldots,\ell\}$ with $|a_{i_0}| \ne |b_{i_0}|$ and, since
$Z(i_0,i_0)=1$, the constant sequence $\jjj = \cdots i_0i_0i_0 \cdots$ lies in $\Sigma_Z$. A direct
calculation shows that
\begin{equation*}
  S_nf_{1,s}(\jjj) - S_nf_{2,s}(\jjj) = c_s n(\log|a_{i_0}| - \log|b_{i_0}|)
\end{equation*}
for all $n \in \N^{+}$, where $c_s = s$ for $0 < s < 1$ and $c_s = 2-s$ for $1 \le s < 2$. As
$|a_{i_0}| \ne |b_{i_0}|$ and $0 < s < 2$, this difference is unbounded in $n$. If $\mu_1 = \mu_2$,
then, since $P_Z(f_{1,s}) = P_Z(f_{2,s})$, the Gibbs bounds for $\mu_1([\jjj|_n])$ and
$\mu_2([\jjj|_n])$ would force it to be bounded, a contradiction. This is where the assumption
$|a_{i_0}| \ne |b_{i_0}|$ for some $i_0 \in \{1,\ldots,\ell\}$ enters. In particular, being
distinct ergodic measures, $\mu_1$ and $\mu_2$ are mutually singular. Since, for all $\iii \in
\Sigma$ and $n \in \N^+$,
\begin{equation} \label{eq:tau_cylinder}
  \begin{split}
    [\tau(\iii)|_n] &= \{\jjj = (j_k)_{k \in \Z} \in \Sigma_Z : j_k=\tau_k(\iii) \text{ for all
    }k\in\{0,\ldots,n-1\}\} \\ 
  &= \tau(\{\jjj = (j_k)_{k \in \Z} \in \Sigma : j_k=i_k \text{ for all }k\in\{0,\ldots,n-1\}\}) =
  \tau([\iii|_n])
\end{split}
\end{equation}
and $\overline{\tau}_*\mu_k([\iii|_n]) = \mu_k([\tau(\iii)|_n]) > 0$, we see that
$\overline{\tau}_*\mu_k(\Sigma)>0$ for both $k \in \{1,2\}$.
Notice also that $\overline{\tau}_*\mu_1 = \overline{\omega}_*\mu_2$ and $\overline{\tau}_*\mu_2 =
\overline{\omega}_*\mu_1$. The $s$-equilibrium state is defined in
\cref{sec:affinity_dim}, the Lyapunov exponents $\lambda_i$ and $\overline{\lambda}$ are defined in
\cref{sec:lyapunov_exponents} and the pressure $P(\A,s)$ is defined in \cref{sec:affinity_dim}.

\begin{proposition}[Fraser--Jordan--Jurga {\cite[Lemma 1.3 and Corollaries 1.5 and 1.6]{FraserJordanJurga2020}}] \label{thm:eq_state_is_sum}
  If $X$ is a proximal box-like self-affine set with rotations whose tuple contains a proximal
  diagonal generator and $\mu \in \MM_\sigma(\Sigma)$ is the unique ergodic $s$-equilibrium state
  with $0<s<2$, then
  \begin{align*}
    P(\A,s) &= P_Z(f_{1,s}) = P_Z(f_{2,s}), \\ 
    \lambda_1(\A,\mu) &= \overline{\lambda}(\A,\mu_1) = \int_{\Sigma_Z} f_{1,1} \dd\mu_1 =
    \int_{\Sigma_Z} f_{2,1} \dd\mu_2, \\ 
    \lambda_2(\A,\mu) &= \overline{\lambda}(\A,\mu_2) = \int_{\Sigma_Z} f_{1,1} \dd\mu_2 =
    \int_{\Sigma_Z} f_{2,1} \dd\mu_1,
  \end{align*}
  and
  \begin{align*}
    \mu &= \overline{\tau}_*\mu_1 + \overline{\tau}_*\mu_2 = \overline{\omega}_*\mu_1 +
    \overline{\omega}_*\mu_2 \\ 
	&= \overline{\tau}_*\mu_1 + \overline{\omega}_*\mu_1 = \overline{\tau}_*\mu_2 +
	\overline{\omega}_*\mu_2, 
  \end{align*}
  where $\mu_1, \mu_2 \in \MM_\sigma(\Sigma_Z)$ are the $s$-Gibbs measures associated to the
  potentials $f_{1,s}$ and $f_{2,s}$, respectively.
\end{proposition}

It follows immediately from \cref{thm:eq_state_is_sum} that $\overline{\tau}_*\mu_1$,
$\overline{\tau}_*\mu_2$, $\overline{\omega}_*\mu_1$, and $\overline{\omega}_*\mu_2$ are absolutely
continuous with respect to $\mu$.
We extend the definitions of $f_{1,s}$ and $f_{2,s}$ to $s=2$ by the same formulae.
At $s=2$, the potentials coincide, $f_{1,2}=f_{2,2}$, and hence so do the Gibbs measures.

\subsection{Dimension drop from coordinate projections} \label{sec:dimension-drop}

We can now justify the dimension-drop assertion in \cref{sec:intro_fraser}. Let $X$ be a proximal
box-like self-affine set with rotations and $s=\dimaff(X)<2$. By
Fraser~\cite[Lemma~2.7]{Fraser2012}, its two coordinate projections have a common Minkowski
dimension $t$. If $t<\min\{1,s\}$, then the modified singular value function at $s$ satisfies
\begin{equation*}
  \alpha_1(A)^t\alpha_2(A)^{s-t}
  = \fii^s(A)\biggl(\frac{\alpha_2(A)}{\alpha_1(A)}\biggr)^{\min\{1,s\}-t}
\end{equation*}
for every product $A$ of matrices in $\A$. The factor on the right decays exponentially along
words typical for the $s$-equilibrium state $\mu$, since
$\lambda_2(\A,\mu)<\lambda_1(\A,\mu)$ by \cref{eq:simple_lyapunov}.

By large deviations for the two Gibbs measures $\mu_1$ and $\mu_2$, whose potentials are
locally constant, and the decomposition in \cref{thm:eq_state_is_sum}, there is $\eps>0$ such
that the words of length $n$ on which this factor exceeds $e^{-\eps n}$ carry exponentially
small $\mu$-measure. Since $P(\A,s)=0$, the Gibbs property of $\mu$ from
\cref{thm:irred_eq_state} makes the pressure of the modified function at $s$ negative. Hence
$\dimmaff(X)<\dimaff(X)$, and \cref{thm:fraser} gives
$\dimh(X)\le\udimm(X)<\dimaff(X)$ without any separation condition.

\section{Hausdorff dimension} \label{sec:hausdorff_dimension}
\subsection{Ledrappier--Young formula} \label{sec:ledrappier_young_formula}
If $\mu$ is a nonzero finite Borel measure on $\R^d$, then the \emph{upper and lower pointwise
dimensions} of $\mu$ at $x \in \R^d$ are
\begin{align*}
  \udimloc(\mu,x) &= \limsup_{r \downarrow 0} \frac{\log \mu(B(x,r))}{\log r}, \\
  \ldimloc(\mu,x) &= \liminf_{r \downarrow 0} \frac{\log \mu(B(x,r))}{\log r},
\end{align*}
respectively. If the limit above exists, then the common value is denoted by $\dimloc(\mu,x)$. We
use the convention that writing $\dimloc(\mu,x)$ implicitly means that the limit exists. If there
exists a constant $s$ such that $\udimloc(\mu,x) = \ldimloc(\mu,x) = s$ for $\mu$-almost all $x \in
\R^d$, then we write $\dim(\mu) = s$ and say that $\mu$ is \emph{exact-dimensional}. We say that
$\mu$ is \emph{double-dimensional} if its pointwise dimension exists at $\mu$-almost every point and
takes exactly two values, each on a set of positive $\mu$-measure. The \emph{upper packing
dimension} and \emph{lower Hausdorff dimension} of $\mu$ are
\begin{align*}
  \udimp(\mu) &= \esssup_{x \sim \mu} \udimloc(\mu,x), \\
  \ldimh(\mu) &= \essinf_{x \sim \mu} \ldimloc(\mu,x),
\end{align*}
respectively. Equivalently, we have
\begin{equation} \label{eq:dimension_measure_def}
  \begin{split}
    \udimp(\mu) &= \inf\{\dimp(A):A\text{ is a Borel set such that }\mu(\R^d \setminus A)=0\}, \\ 
    \ldimh(\mu) &= \inf\{\dimh(A):A\text{ is a Borel set such that }\mu(A)>0\},
  \end{split}
\end{equation}
where $\dimp(A)$ and $\dimh(A)$ denote the packing and Hausdorff dimensions of a set $A$,
respectively, see, for example, Falconer \cite[\S 10]{Falconer1997}. Since $\mu \ll \mu + \nu$,
where $\nu$ is another nonzero finite Borel measure on $\R^d$, it follows from
\cref{eq:dimension_measure_def} that $\udimp(\mu + \nu) \geq \udimp(\mu)$ and $\ldimh(\mu + \nu)
\leq \ldimh(\mu)$. By selecting Borel sets of full measure with packing dimensions arbitrarily close
to $\udimp(\mu)$ and $\udimp(\nu)$, respectively, and a Borel set of positive measure with Hausdorff
dimension arbitrarily close to $\ldimh(\mu + \nu)$, we deduce that
\begin{equation} \label{eq:haus_dim_sum}
  \begin{split}
    \udimp(\mu+\nu) &= \max\{\udimp(\mu), \udimp(\nu)\}, \\ 
    \ldimh(\mu+\nu) &= \min\{\ldimh(\mu), \ldimh(\nu)\}.
  \end{split}
\end{equation}
Notice that $\mu$ is exact-dimensional if and only if $\udimp(\mu) = \ldimh(\mu) = \dim(\mu)$, and
double-dimensional if and only if the pointwise dimension exists almost everywhere, takes only the
values $\udimp(\mu)$ and $\ldimh(\mu)$, and $\ldimh(\mu) < \udimp(\mu)$. We use the convention that
writing $\dim(\mu)$ implicitly means that the measure $\mu$ is exact-dimensional.

Adding two exact-dimensional measures of different dimensions trivially gives a double-dimensional
measure, and we record some simple facts below. We omit its straightforward proof.
\begin{lemma} \label{thm:exact-dimensional-mixture}
  Let $\nu_1$ and $\nu_2$ be nonzero finite exact-dimensional Borel measures on
  $\R^d$, and write $\nu=\nu_1+\nu_2$. For each $k\in\{1,2\}$, we have
  \begin{equation*}
    \dimloc(\nu,x)=\dim(\nu_k)
  \end{equation*}
  for $\nu_k$-almost every $x$. In particular, each value $\dim(\nu_k)$ occurs on a set of positive
  $\nu$-measure, and $\nu$ is double-dimensional if $\dim(\nu_1) \ne \dim(\nu_2)$ and
  exact-dimensional otherwise.
\end{lemma}

The Lyapunov exponents $\lambda_i$ and $\overline{\lambda}$ are defined in \cref{sec:lyapunov_exponents}
and the entropies $h_i$ and $\overline{h}_i$ are defined in \cref{sec:entropies}. Recall also that
if $\mu \in \MM_\sigma(\Sigma)$ and $\nu \in \MM_\sigma(\Sigma_Z)$, then the canonical projections
$\pi_*\mu$ and $\overline{\pi}_*\nu$ are Borel probability measures on $\R^2$ and $\R$,
respectively. 

\begin{theorem}[D.-J.~Feng {\cite[Theorems 1.2 and 1.3]{Feng2023}}] \label{thm:feng1}
  If $X$ is a planar self-affine set and $\mu \in \MM_\sigma(\Sigma)$ is ergodic with
  $\lambda_2(\A,\mu)<\lambda_1(\A,\mu)$, then
  \begin{equation} \label{eq:feng1_1}
    \dim(\pi_*\mu) = \frac{h_1(\mu)-h_0(\mu)}{\lambda_1(\A,\mu)} +
    \frac{h_2(\mu)-h_1(\mu)}{\lambda_2(\A,\mu)}.
  \end{equation}
  Furthermore, if $X$ is a proximal box-like self-affine set with rotations and $\nu \in
  \MM_\sigma(\Sigma_Z)$ is ergodic, then we also have
  \begin{equation*}
    \dim(\overline{\pi}_*\nu) =
    \frac{\overline{h}_1(\nu)-\overline{h}_0(\nu)}{\overline{\lambda}(\A,\nu)}.
  \end{equation*}
\end{theorem}

Rossi \cite[Theorem 2.2]{Rossi2014} proved that for each ergodic measure $\mu \in
\MM_\sigma(\Sigma)$, it holds that
\begin{equation} \label{eq:hausdorff_ergodic_lyapunov}
  \udimh(\pi_*\mu) \le \min\{2,\diml(\A,\mu)\}.
\end{equation}
Therefore, in order to prove that $\dimh(X) = \min\{2,\dimaff(\A)\}$, it suffices to show that if
$s=\dimaff(\A)$, then the $s$-equilibrium state $\mu$ satisfies $\dim(\pi_*\mu) \ge
\min\{2,\diml(\A,\mu)\}$. Under the projective exponential separation condition,
\cref{thm:main_extended} establishes this lower bound up to the threshold of \cref{thm:main2} and,
beyond it, under the strong open set condition, and so proves the claim on $\pi_*\mu$ in
\cref{thm:main2}. The proof of the dimension formula in \cref{thm:main2} instead uses Bernoulli
measures on finite diagonal subsystems, see \cref{sec:dimension_result_proof}. For a proximal
box-like tuple with rotations and $0<s<2$, the $s$-equilibrium state is not a Bernoulli measure, see
\cref{sec:normalisation}. Were it a Bernoulli measure, we could have applied the result of
B\'ar\'any--K\"aenm\"aki \cite[Theorem 2.3]{BaranyKaenmaki2017} to derive \cref{thm:feng1}. However,
as this does not hold, we rely on D.-J.~Feng's general result \cite{Feng2023}, which applies to the
ergodic $s$-equilibrium state.

For the rest of this discussion, fix $0<s<2$ and a proximal box-like tuple with rotations. Recall
from \cref{sec:affinity_dim} that its unique ergodic $s$-equilibrium state is sub-multiplicative.
D.-J.~Feng \cite[Theorem 1.6(iii)]{Feng2023} showed that if the measure is sub-multiplicative, then
on a positive measure set the first summand in \cref{eq:feng1_1} is the pointwise dimension of the
measure projected onto the orthogonal complement $V^\bot$ of $V$ and the second summand is the
pointwise dimension of the measure conditioned on $V$, for almost all $V$ with respect to the
Furstenberg measure. Recall from \cref{sec:lyapunov_exponents} and \cref{eq:simple_lyapunov} that if
$X$ is a proximal box-like self-affine set with rotations, then the support of the Furstenberg
measure $\Pi_*\mu$ is the union of coordinate axes. We will study dimensions of the projected
$s$-equilibrium state in detail in \cref{sec:dimension_of_projections}.

\subsection{Kolmogorov--Sinai entropy of the equilibrium state} \label{sec:entropy_eq_state}
For the normalised tuple and $0<s<2$, we show that the Kolmogorov--Sinai entropy of the
$s$-equilibrium state $\mu$
defined in \cref{sec:affinity_dim} equals the entropies of the $s$-Gibbs measures $\mu_1$ and
$\mu_2$ defined in \cref{sec:gibbs_measures}. 

\begin{proposition} \label{thm:entropy_eq_gibbs}
  If $X$ is a proximal box-like self-affine set with rotations whose tuple contains a proximal
  diagonal generator and $\mu \in \MM_\sigma(\Sigma)$ is the unique ergodic $s$-equilibrium state
  with $0<s<2$, then
  \begin{equation*}
    h_0(\mu) = \overline{h}_0(\mu_1) = \overline{h}_0(\mu_2).
  \end{equation*}
\end{proposition}

\begin{proof}
  Define $q\colon\Sigma_Z\to\Sigma$ by replacing each symbol $j>N$ by $j-N$ and leaving the
  other symbols unchanged. This map is continuous and commutes with $\sigma$. Its restrictions
  to $\tau(\Sigma)$ and $\omega(\Sigma)$ are $\overline{\tau}$ and $\overline{\omega}$,
  respectively, so every $\iii\in\Sigma$ has exactly two lifts, $\tau(\iii)$ and $\omega(\iii)$.
  By \cref{eq:tau_omega_partition,thm:eq_state_is_sum}, we have
  \begin{equation*}
    q_*\mu_k=\overline{\tau}_*\mu_k+\overline{\omega}_*\mu_k=\mu
  \end{equation*}
  for both $k\in\{1,2\}$.

  Recall the coordinate partitions $\cP$ and $\overline{\cP}$ from \cref{sec:entropies}.
  For a base block $i_0\cdots i_{n-1}$, the choice of its initial lifted symbol
  $j_0\in\{i_0,i_0+N\}$ determines all remaining lifted symbols by the transition rule for $Z$.
  Hence each atom of $q^{-1}(\cP_0^{n-1})$ is a union of at most two atoms of
  $\overline{\cP}_0^{n-1}$. Since $q_*\mu_k=\mu$, the entropy chain rule gives
  \begin{equation*}
    0\le H_{\mu_k}(\overline{\cP}_0^{n-1})-H_\mu(\cP_0^{n-1})\le\log2
  \end{equation*}
  for every $n\in\N^+$ and $k\in\{1,2\}$. Dividing by $n$ and passing to the limit gives equality
  of the Kolmogorov--Sinai entropies. By the definitions in \cref{sec:entropies}, $h_0$ and
  $\overline{h}_0$ are the one-symbol entropies conditional on the past, so they equal these
  entropy rates. Thus, for each $k\in\{1,2\}$,
  \begin{equation*}
    h_0(\mu)
    =\lim_{n\to\infty}\frac{1}{n}H_\mu(\cP_0^{n-1})
    =\lim_{n\to\infty}\frac{1}{n}H_{\mu_k}(\overline{\cP}_0^{n-1})
    =\overline{h}_0(\mu_k). \qedhere
  \end{equation*}
\end{proof}

\subsection{Dimension of projections} \label{sec:dimension_of_projections}
We first consider a proximal box-like tuple $\A$ with rotations containing a proximal
diagonal generator, in the coordinates of \cref{sec:subshift}. Fix
$0<s<2$ and let $\mu \in \MM_\sigma(\Sigma)$ be the unique ergodic $s$-equilibrium state as
in \cref{sec:affinity_dim} and $\mu_1, \mu_2 \in \MM_\sigma(\Sigma_Z)$ be the $s$-Gibbs measures as
in \cref{sec:gibbs_measures}. We write
\begin{equation*}
  \nu_x = (\proj_x)_*(\pi_*\mu)
  \quad\text{and}\quad
  \nu_y = (\proj_y)_*(\pi_*\mu)
\end{equation*}
for the coordinate projections. For $k\in\{1,2\}$, we also write
\begin{equation*}
  \nu_k = \overline{\pi}_*\mu_k, \quad \nu_{x,k} = \overline{\pi}_*(\mu_k|_{\tau(\Sigma)}),
  \quad\text{and}\quad \nu_{y,k} = \overline{\pi}_*(\mu_k|_{\omega(\Sigma)}).
\end{equation*}
The decompositions of $\mu$ in \cref{thm:eq_state_is_sum}, together with
\cref{eq:commutative_diagram,eq:inverse-coding-restrictions}, give
\begin{equation} \label{eq:proj_eq_is_sum}
  \nu_x = \nu_{x,1}+\nu_{x,2}
  \quad\text{and}\quad
  \nu_y = \nu_{y,1}+\nu_{y,2}.
\end{equation}
The Gibbs bounds give positive mass to every cylinder, so all four components are nonzero.
Moreover, $\nu_{x,k}\le\nu_k$ and $\nu_{y,k}\le\nu_k$.

Applying Feng--Hu~\cite[Theorem~2.8]{FengHu2009} to the one-sided versions of the ergodic measures
$\mu_k$ shows that $\nu_k$ is exact-dimensional without any separation assumption. Differentiation
of measures then shows that the restricted projected measures have the same dimensions:
\begin{equation} \label{eq:restricted-projection-dimensions}
  \dim(\nu_{x,k}) = \dim(\nu_{y,k}) = \dim(\nu_k)
\end{equation}
for both $k\in\{1,2\}$.

By \cref{eq:proj_eq_is_sum,thm:exact-dimensional-mixture}, the pointwise dimension of both
coordinate projections therefore exists almost everywhere and takes only the values $\dim(\nu_1)$
and $\dim(\nu_2)$, each on a set of positive measure, without any separation assumption. The
following theorem identifies these values under projective exponential separation.

\begin{theorem} \label{thm:eq_state_proj_dim}
  If $X$ is a proximal box-like self-affine set with rotations satisfying the projective exponential
  separation condition and $\mu \in \MM_\sigma(\Sigma)$ is the unique ergodic $s$-equilibrium state
  with $s>0$, then the pointwise dimension of each of $(\proj_x)_*(\pi_*\mu)$ and
  $(\proj_y)_*(\pi_*\mu)$ exists almost everywhere and takes at most two values, namely the upper
  packing and the lower Hausdorff dimension of the measure, each on a set of positive measure, with
  \begin{align*}
    \udimp((\proj_x)_*(\pi_*\mu)) 
    &= \udimp((\proj_y)_*(\pi_*\mu)) = \min\biggl\{ 1,
    -\frac{h_0(\mu)}{\lambda_1(\A,\mu)} \biggr\}, \\
    \ldimh((\proj_x)_*(\pi_*\mu))
    &= \ldimh((\proj_y)_*(\pi_*\mu)) = \min\biggl\{ 1,
    -\frac{h_0(\mu)}{\lambda_2(\A,\mu)} \biggr\}.
  \end{align*}
\end{theorem}

\begin{proof}
  Suppose first that $s\ge2$. By \cref{thm:bernoulli_equilibrium}, $\mu$ is Bernoulli with weights
  $p_i$. The measure $\widehat\mu=(\tau_*\mu+\omega_*\mu)/2$ on $\Sigma_Z$ chooses the initial
  orientation uniformly and then chooses each planar symbol independently with weights $p_i$. It is
  a stationary mixing Markov measure, since diagonal symbols preserve orientation and anti-diagonal
  symbols switch it, and both types have positive probability. Each planar word has exactly two
  orientation lifts, so $\overline h_0(\widehat\mu)=h_0(\mu)$. Birkhoff's theorem applied to the
  logarithms of the scalar contractions along both lifts, together with
  \cref{thm:matrix_products_g_derivatives}, gives
  \begin{equation*}
    \lambda_1(\A,\mu)=\lambda_2(\A,\mu)
    =\overline\lambda(\A,\widehat\mu)
    =\frac12\sum_{i=1}^N p_i\log|\det A_i|.
  \end{equation*}
  By \cref{eq:commutative_diagram}, $\overline\pi_*\widehat\mu=(\nu_x+\nu_y)/2$.
  Jordan--Rapaport~\cite[Theorem~1.1]{JordanRapaport2021} therefore gives this average exact
  dimension $\min\{1,-h_0(\mu)/\lambda_1(\A,\mu)\}$. Both coordinate projections are
  absolutely continuous with respect to the average, so differentiation of measures gives them
  the same exact dimension.

  We may now assume $0<s<2$. By \cref{thm:two-block-normalisation}, grouping preserves both coordinate projections and the ratios of entropy to Lyapunov exponents. We may therefore assume that $\A$ contains a proximal diagonal generator and use the components constructed above. Fix $z\in\{x,y\}$. We prove the assertion for $\nu_z$. Write
  \begin{equation*}
    s_1 = \min\biggl\{ 1, -\frac{h_0(\mu)}{\lambda_1(\A,\mu)} \biggr\} \quad\text{and}\quad
    s_2 = \min\biggl\{ 1, -\frac{h_0(\mu)}{\lambda_2(\A,\mu)} \biggr\}.
  \end{equation*}
  Since $X$ satisfies the projective exponential separation condition, Jordan--Rapaport
  \cite[Theorem~1.1]{JordanRapaport2021} together with
  \cref{thm:entropy_eq_gibbs,thm:eq_state_is_sum} show that
  \begin{equation} \label{eq:gibbs_proj_dim}
    \begin{split}
      \dim(\nu_1)
      &= \min\biggl\{ 1,
      -\frac{\overline{h}_0(\mu_1)}{\overline{\lambda}(\A,\mu_1)} \biggr\} = s_1, \\
      \dim(\nu_2)
      &= \min\biggl\{ 1,
      -\frac{\overline{h}_0(\mu_2)}{\overline{\lambda}(\A,\mu_2)} \biggr\} = s_2.
    \end{split}
  \end{equation}
  By \cref{eq:gibbs_proj_dim,eq:restricted-projection-dimensions}, we have $\dim(\nu_{z,k})=s_k$ for $k\in\{1,2\}$. Since $\nu_z=\nu_{z,1}+\nu_{z,2}$ by \cref{eq:proj_eq_is_sum}, applying \cref{thm:exact-dimensional-mixture} gives
  \begin{equation*}
    \dimloc(\nu_z,t)=s_k
  \end{equation*}
  for $\nu_{z,k}$-almost every $t$. Thus the pointwise dimension of $\nu_z$ exists almost everywhere and takes only the values $s_1$ and $s_2$, each on a set of positive $\nu_z$-measure. Since $\lambda_2(\A,\mu)<\lambda_1(\A,\mu)<0$ by \cref{eq:simple_lyapunov}, we have $s_2\le s_1$. The definitions of upper packing and lower Hausdorff dimension now give
  \begin{equation*}
    \udimp(\nu_z)=s_1 \quad\text{and}\quad \ldimh(\nu_z)=s_2. \qedhere
  \end{equation*}
\end{proof}

To prove \cref{thm:main1}, we compare the two values of the pointwise dimension given by \cref{thm:eq_state_proj_dim}. Since $s=\dimaff(X)$, the definition of the equilibrium state \cref{eq:eq_state_def} expresses the entropy in terms of $s$ and the Lyapunov exponents, and this allows us to determine when the two values differ.

\begin{proof}[Proof of \cref{thm:main1}]
  By \cref{thm:eq_state_proj_dim}, the pointwise dimension of each coordinate projection exists
  almost everywhere and takes only the values $\min\{1,-h_0(\mu)/\lambda_1(\A,\mu)\}$ and
  $\min\{1,-h_0(\mu)/\lambda_2(\A,\mu)\}$, each on a set of positive measure. The projections are
  therefore double-dimensional if these values differ, and exact-dimensional otherwise. If $s\ge2$,
  the proof of \cref{thm:eq_state_proj_dim} gives $\lambda_1(\A,\mu)=\lambda_2(\A,\mu)$, and
  \cref{eq:eq_state_def,eq:energy_def} give $h_0(\mu)=-s\lambda_1(\A,\mu)$. Both values are
  therefore one, as asserted, since the threshold is one.

  We may now assume $0<s<2$ and write $\theta=\lambda_1(\A,\mu)/\lambda_2(\A,\mu)$, so that
  $0<\theta<1$ by \cref{eq:simple_lyapunov}. At the affinity parameter, \cref{eq:eq_state_def} gives
  \begin{equation*}
    h_0(\mu) =
    \begin{cases}
      -s\lambda_1(\A,\mu), &\text{if } 0<s\le1, \\
      -\lambda_1(\A,\mu)-(s-1)\lambda_2(\A,\mu), &\text{if } 1<s<2.
    \end{cases}
  \end{equation*}
  Dividing $h_0(\mu)$ by $-\lambda_1(\A,\mu)$ and $-\lambda_2(\A,\mu)$, we see that
  $-h_0(\mu)/\lambda_1(\A,\mu)=s$ and $-h_0(\mu)/\lambda_2(\A,\mu)=s\theta$ if $0<s\le1$, and
  $-h_0(\mu)/\lambda_1(\A,\mu)=1+(s-1)/\theta>1$ and $-h_0(\mu)/\lambda_2(\A,\mu)=s-1+\theta$ if
  $1<s<2$. Hence the first value is $\min\{1,s\}$, and the second is $s\theta$ if $0<s\le1$ and
  $\min\{1,s-1+\theta\}$ if $1<s<2$. If $0<s\le1$, the values $s$ and $s\theta$ differ, and
  $s<2-\theta$. If $1<s<2$, the first value is one, and the second is less than one precisely when
  $s<2-\theta$. Thus the two values differ precisely when $s<2-\theta$, and otherwise both are one.
\end{proof}

\begin{example} \label{ex:projection_double_dimensional}
  For every $0<s_2<s_1\le1$, we construct a proximal box-like self-affine set $X$ with
  rotations, defined by four maps and satisfying strong separation and projective exponential
  separation, such that $\dimaff(X)=s_1$ and both coordinate projections of its
  $s_1$-equilibrium state are double-dimensional with pointwise dimension values $s_1$ and $s_2$.
  Write
  \begin{equation*}
    H(p)=-p\log p-(1-p)\log(1-p),
    \qquad J(p)=-p\log(1-p)-(1-p)\log p,
  \end{equation*}
  and choose $p\in(1/2,1)$ such that
  \begin{equation*}
    \frac{\log2+H(p)}{\log2+J(p)}=\frac{s_2}{s_1}.
  \end{equation*}
  This choice is unique: $H$ decreases and $J$ increases on $(1/2,1)$, and the ratio decreases
  continuously from one to zero. Set $a=(p/2)^{1/s_1}$ and $b=((1-p)/2)^{1/s_1}$, and define
  \begin{align*}
    \fii_i(x,y)&=(ax+t_i,by+t_{i+4}), \quad i\in\{1,2\},\\
    \fii_i(x,y)&=(ay+t_i,bx+t_{i+4}), \quad i\in\{3,4\}.
  \end{align*}
  Since $0<b<a<1/2$, the diagonal generators are proximal and the anti-diagonal generators
  interchange their invariant axes, so the tuple is irreducible.

  To choose the translations, put
  \begin{equation*}
    R=\frac{2}{1-2a} \quad\text{and}\quad B=(0,R,0,R,0,0,R,R)+[0,1]^8.
  \end{equation*}
  For every $t\in B$, the four maps send $[0,2R]^2$ into itself. Indeed,
  $2aR+1=R-1$, so each coordinate of an image rectangle lies in $[0,R-1]$ or
  $[R,2R-1]$, according to its translation interval. The four rectangles occupy different
  pairs of these intervals and are therefore separated by gaps of at least one.
  Thus every $t\in B$ gives strong separation.

  We use Lebesgue measure $\LL^8$ on the translation cube $B\subset\R^8$ to show that some $t$ also
  gives projective exponential separation. Write $g_i^t(x)=ax+t_i$ and
  $g_{i+4}^t(x)=bx+t_{i+4}$ for $1\le i\le4$. For distinct words
  $\iii=i_1\cdots i_n$ and $\jjj=j_1\cdots j_n$ in $\{1,\ldots,8\}^n$, the difference
  $L_{\iii,\jjj}(t)=g_\iii^t(0)-g_\jjj^t(0)$ is linear in $t$. Let
  $k=\min\{m:i_m\ne j_m\}$ be the first index where the words differ. Since $i_m=j_m$
  for $m<k$, we have
  \begin{equation*}
    L_{\iii,\jjj}(t)
    =(g_{i_1\cdots i_{k-1}}^t)'(0)
      \bigl(g_{i_k\cdots i_n}^t(0)-g_{j_k\cdots j_n}^t(0)\bigr),
  \end{equation*}
  where the first composition is the identity when $k=1$. The coefficient of $t_{i_k}$ in
  $g_{i_k\cdots i_n}^t(0)$ is at least one, whereas its coefficient in
  $g_{j_k\cdots j_n}^t(0)$ is at most $a+a^2+\cdots=a/(1-a)$, since $j_k\ne i_k$.
  Their difference is therefore at least $\kappa=(1-2a)/(1-a)>0$. As
  $(g_{i_1\cdots i_{k-1}}^t)'(0)\ge b^{k-1}\ge b^{n-1}$, the coefficient of $t_{i_k}$
  in $L_{\iii,\jjj}$ is at least $\kappa b^{n-1}$.

  Now fix the other seven translations. As a function of $t_{i_k}$, the difference
  $L_{\iii,\jjj}$ has slope at least $\kappa b^{n-1}$, so
  $|L_{\iii,\jjj}(t)|<\eps$ can hold only on an interval of length at most
  $2\eps/(\kappa b^{n-1})$. Each parameter interval in $B$ has length one. Integrating over
  the other seven translations therefore gives
  \begin{equation*}
    \LL^8(\{t\in B:|L_{\iii,\jjj}(t)|<\eps\})
      \le\frac{2\eps}{\kappa b^{n-1}}.
  \end{equation*}
  We now sum this estimate over all word pairs and lengths. Set $C=\kappa/8$ and $c=b/128$,
  and let $E_n\subset B$ consist of the translations for
  which $|L_{\iii,\jjj}(t)|<Cc^n$ for at least one pair of distinct words of length $n$.
  There are fewer than $8^{2n}=64^n$ such pairs, so
  \begin{equation*}
    \LL^8(E_n)\le64^n\frac{2Cc^n}{\kappa b^{n-1}}=\frac b4\,2^{-n}
    \quad\text{and}\quad
    \LL^8\biggl(\bigcup_{n\ge1}E_n\biggr)\le\frac b4<1=\LL^8(B).
  \end{equation*}
  We can consequently choose $t\in B\setminus\bigcup_{n\ge1}E_n$. For this fixed choice,
  every pair of distinct length-$n$ words satisfies
  $|g_\iii^t(0)-g_\jjj^t(0)|\ge Cc^n$. Since $0<C\le1$, we have
  $\Delta_n^\Psi\ge Cc^n\ge(Cc)^n$ for every $n\ge1$, giving projective exponential separation.
  To rule out exact overlaps, suppose that $g_\iii^t=g_\jjj^t$ for distinct finite words
  $\iii,\jjj$. If one word is a prefix of the other, cancelling the corresponding injective
  similarity makes a non-empty composition the identity, contradicting strict contraction.
  Otherwise, $\iii\jjj$ and $\jjj\iii$ are distinct words of equal length, yet
  $g_{\iii\jjj}^t=g_{\jjj\iii}^t$, contradicting the separation bound above.

  For $0\le u\le1$, the matrix whose $(i,j)$ entry sums the $u$th powers of the contractions
  from coordinate $j$ to coordinate $i$ is
  \begin{equation*}
    M(u)=2\begin{pmatrix}a^u&a^u\\b^u&b^u\end{pmatrix}.
  \end{equation*}
  Each product $A_\iii$ has two nonzero entries, whose absolute values are its singular values.
  The sum of their $u$th powers lies between $\fii^u(A_\iii)$ and $2\fii^u(A_\iii)$, so
  \begin{equation*}
    \sum_{|\iii|=n}\fii^u(A_\iii)
    \le(1,1)M(u)^n\begin{pmatrix}1\\1\end{pmatrix}
    \le2\sum_{|\iii|=n}\fii^u(A_\iii).
  \end{equation*}
  Taking logarithms, dividing by $n$ and letting $n\to\infty$ gives
  $P(\A,u)=\log(2a^u+2b^u)$, since $2a^u+2b^u$ is the only nonzero eigenvalue of $M(u)$.
  Since $2a^{s_1}+2b^{s_1}=1$, the affinity dimension is $s=s_1$. Let $\mu$ be the
  $s$-equilibrium state.

  Write $v=(p,1-p)$. Since $M(s)v=v$, the transition probabilities on the two orientation
  states of the $s$-Gibbs measure $\mu_1$ are $Q_{ij}=M(s)_{ij}v_j/v_i$, giving
  \begin{equation*}
    Q=\begin{pmatrix}p&1-p\\p&1-p\end{pmatrix}.
  \end{equation*}
  Its stationary distribution is $(p,1-p)$. Each transition has two equiprobable map choices,
  so from either state the four labelled transitions have probabilities
  $p/2,p/2,(1-p)/2,(1-p)/2$. By \cref{thm:entropy_eq_gibbs},
  \begin{equation*}
    h_0(\mu)=-p\log(p/2)-(1-p)\log((1-p)/2)=\log2+H(p).
  \end{equation*}
  The selected coordinate contracts by $a$ in the first state and by $b$ in the second,
  while the other coordinate has the complementary contraction. Averaging with weights
  $(p,1-p)$ and using \cref{thm:eq_state_is_sum} therefore gives
  \begin{align*}
    \lambda_1(\A,\mu)&=p\log a+(1-p)\log b=-\frac{\log2+H(p)}{s_1},\\
    \lambda_2(\A,\mu)&=p\log b+(1-p)\log a=-\frac{\log2+J(p)}{s_1}.
  \end{align*}
  Consequently,
  \begin{equation*}
    -\frac{h_0(\mu)}{\lambda_1(\A,\mu)}=s_1
    \quad\text{and}\quad
    -\frac{h_0(\mu)}{\lambda_2(\A,\mu)}
      =s_1\frac{\log2+H(p)}{\log2+J(p)}=s_2.
  \end{equation*}
  By \cref{thm:eq_state_proj_dim}, the coordinate projections $\nu_x$ and $\nu_y$ are
  double-dimensional, with
  \begin{equation*}
    \udimp(\nu_x)=\udimp(\nu_y)=s_1
    \quad\text{and}\quad
    \ldimh(\nu_x)=\ldimh(\nu_y)=s_2.
  \end{equation*}
  Each pointwise dimension value occurs on a set of positive measure in each projection.
\end{example}

\subsection{Conditional measures of the average measure} \label{sec:average_gibbs}
Fix $0<s<2$. Let $\mu_1$ and $\mu_2$ be the $s$-Gibbs measures of the normalised tuple
defined in \cref{sec:gibbs_measures}, and define
\begin{equation} \label{eq:average_gibbs}
  \overline{\mu} = \frac{\mu_1+\mu_2}{2}
\end{equation}
to be their average. In this section, we examine the conditional measures of $\overline{\mu}$
associated with $\overline{\xi}_1$ and determine conditions under which its entropy vanishes.

Suppose that $\A \in \GL^N$ is proximal and box-like with rotations, containing a diagonal matrix
whose diagonal entries have distinct absolute values, and let $\mu \in \MM_\sigma(\Sigma)$ be the
unique ergodic $s$-equilibrium state with $0<s<2$.
Recalling \cref{eq:simple_lyapunov}, let $\Sigma' \subseteq \Sigma$ be the set of full $\mu$-measure
as in \cref{sec:lyapunov_exponents}. Recall from \cref{eq:sigma_tau_omega} that
\begin{align*}
  \Sigma^\tau &= \{\iii \in \Sigma' : V_\iii^1 \text{ is the $y$-axis}\}, \\ 
  \Sigma^\omega &= \{\iii \in \Sigma' : V_\iii^1 \text{ is the $x$-axis}\}.
\end{align*}

\begin{lemma} \label{thm:gibbs_tau_omega_positive}
  If $X$ is a proximal box-like self-affine set with rotations whose tuple contains a proximal
  diagonal generator and $\mu_1, \mu_2 \in \MM_\sigma(\Sigma_Z)$ are the $s$-Gibbs measures with
  $0<s<2$, then
  \begin{align*}
    \mu_1(\tau(\Sigma^\tau)) + \mu_1(\omega(\Sigma^\omega)) = \mu_2(\omega(\Sigma^\tau)) +
    \mu_2(\tau(\Sigma^\omega)) &= 1, \\
    \mu_2(\tau(\Sigma^\tau)) = \mu_2(\omega(\Sigma^\omega)) = \mu_1(\omega(\Sigma^\tau)) =
    \mu_1(\tau(\Sigma^\omega)) &= 0.
  \end{align*}
  Furthermore, the four summands in the first line are positive:
  \begin{equation*}
    \min\{\mu_1(\tau(\Sigma^\tau)), \mu_1(\omega(\Sigma^\omega)), \mu_2(\omega(\Sigma^\tau)),
    \mu_2(\tau(\Sigma^\omega))\} > 0.
  \end{equation*}
\end{lemma}

\begin{proof}
  Let $\mu \in \MM_\sigma(\Sigma)$ be the $s$-equilibrium state given by \cref{thm:irred_eq_state}.
  By \cref{thm:eq_state_is_sum} and \cref{eq:lyapunov_exponent_gsystem}, 
  \begin{equation*}
    \lambda_1(\A,\mu) = \overline{\lambda}(\A,\mu_1) = \lim_{n \to \infty} \frac{1}{n} \log
    |g'_{i_{-n} \cdots i_{-1}}(0)|
  \end{equation*}
  for $\mu_1$-almost all $\iii \in \Sigma_Z$ and
  \begin{equation*}
    \lambda_2(\A,\mu) = \overline{\lambda}(\A,\mu_2) = \lim_{n \to \infty} \frac{1}{n} \log
    |g'_{i_{-n} \cdots i_{-1}}(0)|
  \end{equation*}
  for $\mu_2$-almost all $\iii \in \Sigma_Z$. By \cref{thm:lyapunov_exponents_via_g_system}, we have
  \begin{equation*}
    \lim_{n \to \infty} \frac{1}{n} \log |g'_{i_{-n} \cdots i_{-1}}(0)| = 
    \begin{cases}
      \lambda_1(\A,\mu), &\text{if } \iii \in \tau(\Sigma^\tau) \cup \omega(\Sigma^\omega), \\ 
      \lambda_2(\A,\mu), &\text{if } \iii \in \tau(\Sigma^\omega) \cup \omega(\Sigma^\tau),
    \end{cases}
  \end{equation*}
  for all $\iii\in\tau(\Sigma')\cup\omega(\Sigma')$. This union has full measure for both
  $\mu_1$ and $\mu_2$, since $\mu(\Sigma')=1$ and $q_*\mu_k=\mu$ for $k\in\{1,2\}$, where
  $q$ is the orientation-forgetting map from \cref{sec:entropy_eq_state}.
  Since $\mu = \overline{\tau}_*\mu_1 +
  \overline{\omega}_*\mu_1 = \overline{\tau}_*\mu_2 + \overline{\omega}_*\mu_2$ by
  \cref{thm:eq_state_is_sum}, the first claim follows by recalling \cref{eq:tau_omega_partition} and
  \cref{eq:tau_omega_partition2} and combining the above identities.

  Since $\mu(\Sigma^\tau) > 0$ and $\mu(\Sigma^\omega) > 0$ by
  \cref{thm:sigma_tau_omega_positive_eq_measure} and $\mu = \overline{\tau}_*\mu_1 +
  \overline{\tau}_*\mu_2 = \overline{\omega}_*\mu_1 + \overline{\omega}_*\mu_2$ by
  \cref{thm:eq_state_is_sum}, the second claim follows by combining this information with the first
  claim.
\end{proof}

Let $T\in\{\tau,\omega\}$ and $E\in\{\Sigma^\tau,\Sigma^\omega\}$. By
\cref{thm:eq_state_is_sum,eq:inverse-coding-restrictions}, the unnormalised restrictions satisfy
\begin{equation} \label{eq:orientation-restrictions}
  \overline\mu|_{T(E)}=\frac12 T_*(\mu|_E) \quad\text{and}\quad \overline\mu(T(E))=\frac12\mu(E)>0.
\end{equation}
The positivity follows from \cref{thm:sigma_tau_omega_positive_eq_measure}. Normalising gives
\begin{equation} \label{eq:normalised-orientation-restrictions}
  \overline\mu_{T(E)}=T_*\mu_E.
\end{equation}
By \cref{thm:gibbs_tau_omega_positive}, this normalised restriction also equals
$(\mu_1)_{T(E)}$ for $(T,E)=(\tau,\Sigma^\tau)$ or $(\omega,\Sigma^\omega)$, and equals
$(\mu_2)_{T(E)}$ in the other two cases.

\begin{lemma} \label{thm:conditional_two_ways}
  Let $X$ be a proximal box-like self-affine set with rotations whose tuple contains a proximal
  diagonal generator, $\mu \in \MM_\sigma(\Sigma)$ be the unique ergodic $s$-equilibrium state with
  $0<s<2$, $\mu_1, \mu_2 \in \MM_\sigma(\Sigma_Z)$ be the $s$-Gibbs measures defined in
  \cref{sec:gibbs_measures}, and $\overline{\mu} \in \MM_\sigma(\Sigma_Z)$ be their average defined
  in \cref{eq:average_gibbs}. Then
  \begin{equation*}
    \overline{\mu}_\iii^{\overline{\xi}_1} =
    \begin{cases}
      (\mu_1)_\iii^{\overline{\xi}_1} = \tau_*\mu_{\overline{\tau}(\iii)}^{\xi_1}, 
      &\text{if } \iii
      \in \tau(\Sigma^\tau), \\[3pt] 
      (\mu_2)_\iii^{\overline{\xi}_1} = \omega_*\mu_{\overline{\omega}(\iii)}^{\xi^\bot_1},
      &\text{if } \iii \in \omega(\Sigma^\tau), \\[3pt]
      (\mu_1)_\iii^{\overline{\xi}_1} = \omega_*\mu_{\overline{\omega}(\iii)}^{\xi_1}, 
      &\text{if }
      \iii \in \omega(\Sigma^\omega), \\[3pt] 
      (\mu_2)_\iii^{\overline{\xi}_1} = \tau_*\mu_{\overline{\tau}(\iii)}^{\xi^\bot_1}, 
      &\text{if }
      \iii \in \tau(\Sigma^\omega), 
    \end{cases}
  \end{equation*}
  for $\overline{\mu}$-almost all $\iii \in \Sigma_Z$.
\end{lemma}

\begin{proof}
  Fix one of the four cases in the statement. Write $T$ for its orientation map, $E$ for the
  corresponding set $\Sigma^\tau$ or $\Sigma^\omega$, and $\zeta$ for its planar partition
  $\xi_1$ or $\xi_1^\bot$. Set $F=T(E)$, and let $k\in\{1,2\}$ be the Gibbs index in that
  case. By \cref{eq:normalised-orientation-restrictions} and the observation following it,
  \begin{equation*}
    \overline\mu_F=T_*\mu_E=(\mu_k)_F.
  \end{equation*}
  By \cref{thm:partitions_xi1}, $E$ and $F$ are unions of the corresponding partition atoms,
  and $T$ maps $\zeta|_E$ to $\overline\xi_1|_F$. Pushing forward the conditional measures
  of $\mu_E$ therefore gives a disintegration of $\overline\mu_F$. By uniqueness,
  \begin{equation*}
    (\overline\mu_F)_\jjj^{\overline\xi_1}
    =((\mu_k)_F)_\jjj^{\overline\xi_1}
    =T_*(\mu_E)_{T^{-1}(\jjj)}^\zeta
  \end{equation*}
  for $\overline\mu$-almost every $\jjj\in F$. Applying
  \cref{thm:conditional_measures_of_a_conditional_measure} removes the restrictions to $E$
  and $F$, proving all four claims.
\end{proof}

\subsection{Pointwise one-dimensional entropy of the average measure} \label{sec:average_gibbs_entropy}
Fix $0<s<2$, let $\mu_1$ and $\mu_2$ be the $s$-Gibbs measures of the normalised tuple
defined in \cref{sec:gibbs_measures}, and let
\begin{equation*}
  \overline{\mu} = \frac{\mu_1+\mu_2}{2}
\end{equation*}
be their average as introduced in \cref{eq:average_gibbs}. In this section, we study the pointwise
$1$-dimensional entropy of $\overline{\mu}$. Our first lemma shows that, under the projective
exponential separation, $\overline{\mu}$ has zero pointwise $1$-dimensional entropy,
$\overline{h}_1(\overline{\mu},\iii)=0$, for $\overline{\mu}$-almost all $\iii \in \Sigma_Z$ when
the Kolmogorov--Sinai entropy of the $s$-equilibrium state satisfies $h_0(\mu) \le
-\lambda_1(\A,\mu)$. 

\begin{lemma} \label{thm:overline_h1_zero}
  If $X$ is a proximal box-like self-affine set with rotations satisfying the projective exponential
  separation condition, $\A$ contains a proximal diagonal generator, and $\mu \in
  \MM_\sigma(\Sigma)$ is the unique ergodic $s$-equilibrium state with $0<s<2$ such that $h_0(\mu)
  \le -\lambda_1(\A,\mu)$, then the measure $\overline{\mu} \in \MM_\sigma(\Sigma_Z)$ defined in
  \cref{eq:average_gibbs} satisfies
  \begin{equation*}
    \overline{h}_1(\overline{\mu},\iii)=0
  \end{equation*}
  for $\overline{\mu}$-almost all $\iii \in \Sigma_Z$.
\end{lemma}

\begin{proof}
  By Feng--Hu \cite[Theorem 2.8]{FengHu2009}, we see that the canonical projections
  $\overline{\pi}_*\mu_1$ and $\overline{\pi}_*\mu_2$ of the $s$-Gibbs measures $\mu_1, \mu_2 \in
  \MM_\sigma(\Sigma_Z)$ defined in \cref{sec:gibbs_measures} are exact-dimensional. Therefore, since
  $X$ satisfies the projective exponential separation condition, Jordan--Rapaport \cite[Theorem
  1.1]{JordanRapaport2021} together with \cref{thm:entropy_eq_gibbs,thm:eq_state_is_sum} show that
  \begin{align*}
    \dim(\overline{\pi}_*\mu_1) 
    &= \min\biggl\{ 1,
    -\frac{\overline{h}_0(\mu_1)}{\overline{\lambda}(\A,\mu_1)} \biggr\} =
    -\frac{h_0(\mu)}{\lambda_1(\A,\mu)} \le 1, \\
    \dim(\overline{\pi}_*\mu_2) 
    &= \min\biggl\{ 1,
    -\frac{\overline{h}_0(\mu_2)}{\overline{\lambda}(\A,\mu_2)} \biggr\} =
    -\frac{h_0(\mu)}{\lambda_2(\A,\mu)} \le -\frac{h_0(\mu)}{\lambda_1(\A,\mu)},
  \end{align*}
  where the last inequality follows from \cref{eq:simple_lyapunov}. Since, by \cref{thm:feng1},
  \begin{align*}
    \dim(\overline{\pi}_*\mu_1) 
	&=
	\frac{\overline{h}_1(\mu_1)-\overline{h}_0(\mu_1)}{\overline{\lambda}(\A,\mu_1)} =
	\frac{\overline{h}_1(\mu_1)-h_0(\mu)}{\lambda_1(\A,\mu)}, \\ 
	\dim(\overline{\pi}_*\mu_2) 
	&=
	\frac{\overline{h}_1(\mu_2)-\overline{h}_0(\mu_2)}{\overline{\lambda}(\A,\mu_2)} =
	\frac{\overline{h}_1(\mu_2)-h_0(\mu)}{\lambda_2(\A,\mu)},
  \end{align*}
  we see that
  \begin{equation} \label{eq:h1_gibbs_entropy_zero}
    \overline{h}_1(\mu_1)=0 \quad\text{and}\quad \overline{h}_1(\mu_2)=0.
  \end{equation}
  Since $\mu_1,\mu_2 \in \MM_\sigma(\Sigma_Z)$ are ergodic, we conclude from
  \cref{thm:feng_entropy_as_limit,thm:conditional_two_ways,thm:gibbs_tau_omega_positive} that
  \begin{align*}
    \overline{h}_1(\overline{\mu},\iii)
	&= -\lim_{n \to \infty}
	\frac{1}{n}\log\overline{\mu}_\iii^{\overline{\xi}_1}(\overline{\cP}_0^{n-1}(\iii)) \\
	&=
	\begin{cases}
	  -\lim_{n \to \infty} \frac{1}{n}\log(\mu_1)_\iii^{\overline{\xi}_1}(\overline{\cP}_0^{n-1}(\iii)) =
	  \overline{h}_1(\mu_1), &\text{if } \iii \in \tau(\Sigma^\tau) \cup \omega(\Sigma^\omega), \\
	  -\lim_{n \to \infty} \frac{1}{n}\log(\mu_2)_\iii^{\overline{\xi}_1}(\overline{\cP}_0^{n-1}(\iii)) =
	  \overline{h}_1(\mu_2), &\text{if } \iii \in \omega(\Sigma^\tau) \cup \tau(\Sigma^\omega),
	\end{cases} \\ 
	&= 0
  \end{align*}
  by \cref{eq:h1_gibbs_entropy_zero} for $\overline{\mu}$-almost all $\iii \in \Sigma_Z$ as claimed.
\end{proof}

The following lemma shows that the one-dimensional entropy of the $s$-Gibbs measure $\mu_2$ vanishes
when the Kolmogorov--Sinai entropy of the $s$-equilibrium state satisfies $h_0(\mu) \le
-\lambda_2(\A,\mu)$.

\begin{lemma} \label{thm:overline_h1_zero_large}
  If $X$ is a proximal box-like self-affine set with rotations satisfying the projective exponential
  separation condition, $\A$ contains a proximal diagonal generator, and $\mu \in
  \MM_\sigma(\Sigma)$ is the unique ergodic $s$-equilibrium state with $0<s<2$ such that
  $h_0(\mu) \le -\lambda_2(\A,\mu)$, then $\overline{h}_1(\mu_2)=0$ and the
  measure $\overline{\mu} \in \MM_\sigma(\Sigma_Z)$ defined in \cref{eq:average_gibbs} satisfies
  \begin{equation*}
    \overline{\mu}(\{\iii \in \omega(\Sigma^\tau) \cup
    \tau(\Sigma^\omega) : \overline{\mu}_\iii^{\overline{\xi}_1}=\delta_\iii\}) =
    \overline{\mu}(\omega(\Sigma^\tau) \cup \tau(\Sigma^\omega)).
  \end{equation*}
\end{lemma}

\begin{proof}
  By \cref{thm:entropy_eq_gibbs} and \cref{thm:eq_state_is_sum}, the assumption $h_0(\mu) \le
  -\lambda_2(\A,\mu)$ implies
  \begin{equation*}
    -\frac{\overline{h}_0(\mu_2)}{\overline{\lambda}(\A,\mu_2)} =
    -\frac{h_0(\mu)}{\lambda_2(\A,\mu)} \le 1.
  \end{equation*}
  Since $X$ satisfies the projective exponential separation condition, Jordan--Rapaport \cite[Theorem
  1.1]{JordanRapaport2021} thus gives
  \begin{equation} \label{eq:overline_h1_zero_large1}
    \dim(\overline{\pi}_*\mu_2) = -\frac{\overline{h}_0(\mu_2)}{\overline{\lambda}(\A,\mu_2)}.
  \end{equation}
  By \cref{thm:feng1}, we have
  \begin{equation} \label{eq:overline_h1_zero_large2}
    \dim(\overline{\pi}_*\mu_2) =
    \frac{\overline{h}_1(\mu_2)-\overline{h}_0(\mu_2)}{\overline{\lambda}(\A,\mu_2)}.
  \end{equation}
  Combining \cref{eq:overline_h1_zero_large1} and \cref{eq:overline_h1_zero_large2}, we obtain
  $\overline{h}_1(\mu_2)=0$ which, by \cref{thm:feng_entropy_zero}, means that
  \begin{equation*}
    \mu_2(\{\iii \in \Sigma_Z : (\mu_2)_\iii^{\overline{\xi}_1}=\delta_\iii\}) = 1.
  \end{equation*}
  Recalling \cref{thm:conditional_two_ways}, we thus see that
  $\overline{\mu}_\iii^{\overline{\xi}_1} = \delta_\iii$ for $\mu_2$-almost all $\iii \in
  \omega(\Sigma^\tau) \cup \tau(\Sigma^\omega)$. Since $\mu_1(\omega(\Sigma^\tau)) =
  \mu_1(\tau(\Sigma^\omega)) = 0$ by \cref{thm:gibbs_tau_omega_positive}, we have finished the
  proof.
\end{proof}

The preceding lemma only gives vanishing on the part carrying $\mu_2$. On the part carrying
$\mu_1$, the pointwise entropy of $\overline\mu$ agrees with $h_1(\mu)$, as we now show.
This comparison does not require projective separation.

\begin{lemma} \label{thm:overline_h1_large}
  If $X$ is a proximal box-like self-affine set with rotations, $\A$ contains a proximal diagonal
  generator, and $\mu\in\MM_\sigma(\Sigma)$ is the unique ergodic $s$-equilibrium state with
  $0<s<2$, then the measure
  $\overline{\mu} \in \MM_\sigma(\Sigma_Z)$ defined in \cref{eq:average_gibbs} satisfies
  \begin{equation*}
    \overline{h}_1(\overline{\mu},\iii) = \overline{h}_1(\mu_1) = h_1(\mu)
  \end{equation*}
  for $\overline{\mu}$-almost all $\iii \in \tau(\Sigma^\tau) \cup \omega(\Sigma^\omega)$. 
\end{lemma}

\begin{proof}
  Consider $(T,E)=(\tau,\Sigma^\tau)$ or $(\omega,\Sigma^\omega)$. On $T(E)$, we have
  $\overline\mu|_{T(E)}=\tfrac12\mu_1|_{T(E)}=\tfrac12 T_*(\mu|_E)$ by
  \cref{thm:gibbs_tau_omega_positive,eq:orientation-restrictions}.
  By \cref{thm:conditional_two_ways} and the preimage identity for cylinders in the proof of
  \cref{thm:partition_tau_conj}, we have
  \begin{equation*}
    \overline\mu_{T(\iii)}^{\overline\xi_1}(\overline\cP_0^{n-1}(T(\iii))) = (\mu_1)_{T(\iii)}^{\overline\xi_1}(\overline\cP_0^{n-1}(T(\iii))) = \mu_\iii^{\xi_1}(\cP_0^{n-1}(\iii))
  \end{equation*}
  for $\mu$-almost every $\iii\in E$ and every $n\ge1$.
  Taking logarithms, dividing by $-n$ and applying \cref{thm:feng_entropy_as_limit} gives
  \begin{equation*}
    \overline h_1(\overline\mu,T(\iii))
    =\overline h_1(\mu_1,T(\iii))=h_1(\mu,\iii).
  \end{equation*}
  Since $\mu_1$ and $\mu$ are ergodic, the last two quantities equal
  $\overline h_1(\mu_1)$ and $h_1(\mu)$ almost everywhere, proving the assertion.
\end{proof}

\subsection{Higher dimensional entropies of the equilibrium state} \label{sec:higher_dimensional_entropies}
As outlined in \cref{sec:ledrappier_young_formula}, the Ledrappier--Young formula provides a method
to calculate the Hausdorff dimension of the self-affine set, requiring knowledge of the entropies
$h_i(\mu)$ of the $s$-equilibrium state $\mu$. Our previous analysis in \cref{sec:entropy_eq_state}
established that the Kolmogorov--Sinai entropy of $\mu$ equals
$h_0(\mu)=\overline{h}_0(\mu_1)=\overline{h}_0(\mu_2)$, where $\mu_1$ and $\mu_2$ are the $s$-Gibbs
measures defined in \cref{sec:gibbs_measures}. 

In this section, we study the higher dimensional entropies of the $s$-equilibrium state. In
particular, we consider horizontal and vertical slices of $\mu$ and, by leveraging the projective
exponential separation defined in \cref{sec:exponential_separation}, we demonstrate that $h_1(\mu)$
vanishes when $h_0(\mu) \le -\lambda_1(\A,\mu)$ and equals $h_0(\mu)+\lambda_1(\A,\mu)$ when the
Kolmogorov--Sinai entropy of the $s$-equilibrium state lies between the negative of the Lyapunov
exponents. We also show that $h_2(\mu)$ vanishes when $h_0(\mu) \le -\lambda_2(\A,\mu)$.

Let $\xi_i$ be the measurable partition of $\Sigma'$ and $\widehat{\xi}_i$ the $\sigma$-algebra
generated by $\xi_i$ as in \cref{sec:entropies}. Relying on the analysis done with the average of
the $s$-Gibbs measures in \cref{sec:average_gibbs}, we first prove that the $s$-equilibrium state
$\mu$ has zero $1$-dimensional entropy, $h_1(\mu)=0$, when $h_0(\mu) \le -\lambda_1(\A,\mu)$.

\begin{proposition} \label{thm:eq_state_h1_zero}
  If $X$ is a proximal box-like self-affine set with rotations satisfying the projective exponential
  separation condition, $\A$ contains a proximal diagonal generator, and $\mu \in
  \MM_\sigma(\Sigma)$ is the $s$-equilibrium state with $0<s<2$ such that $h_0(\mu) \le
  -\lambda_1(\A,\mu)$, then
  \begin{equation*}
    h_1(\mu)=0 \quad\text{and}\quad \mu(\{\iii \in \Sigma' : \mu_\iii^{\xi_1} = \delta_\iii\}) = 1.
  \end{equation*}
\end{proposition}

\begin{proof}
  By \cref{thm:overline_h1_zero,thm:feng_entropy_zero}, we have
  $\overline\mu_\jjj^{\overline\xi_1}=\delta_\jjj$ for $\overline\mu$-almost every $\jjj$.
  For $(T,E)=(\tau,\Sigma^\tau)$ or $(\omega,\Sigma^\omega)$,
  \cref{thm:conditional_two_ways,eq:orientation-restrictions} therefore gives
  \begin{equation*}
    T_*\mu_\iii^{\xi_1}=\delta_{T(\iii)}
  \end{equation*}
  for $\mu$-almost every $\iii\in E$. 
  Since $T$ is injective, $\mu_\iii^{\xi_1}=\delta_\iii$ at these points.
  The sets $\Sigma^\tau$ and $\Sigma^\omega$ partition $\Sigma'$, so this identity holds
  $\mu$-almost everywhere. Now \cref{thm:feng_entropy_zero} gives $h_1(\mu)=0$.
\end{proof}

The next proposition calculates the $1$-dimensional entropy of the $s$-equilibrium state $\mu$ when
the Kolmogorov--Sinai entropy of $\mu$ is large.

\begin{proposition} \label{thm:eq_state_h1_large}
  If $X$ is a proximal box-like self-affine set with rotations satisfying the projective exponential
  separation condition, $\A$ contains a proximal diagonal generator, and $\mu \in
  \MM_\sigma(\Sigma)$ is the $s$-equilibrium state with $0<s<2$ such that $h_0(\mu) \ge
  -\lambda_1(\A,\mu)$, then
  \begin{equation*}
    h_1(\mu)=h_0(\mu)+\lambda_1(\A,\mu).
  \end{equation*}
\end{proposition}

\begin{proof}
  Let $\mu_1$ be the $s$-Gibbs measure defined in \cref{sec:gibbs_measures}. Since
  $-\lambda_1(\A,\mu) \le h_0(\mu)$ and $X$ satisfies the projective exponential separation
  condition, Jordan--Rapaport \cite[Theorem 1.1]{JordanRapaport2021} together with
  \cref{thm:entropy_eq_gibbs,thm:eq_state_is_sum} show that
  \begin{equation} \label{eq:eq_state_h1_large1}
    \dim(\overline{\pi}_*\mu_1) = \min\biggl\{ 1,
    -\frac{\overline{h}_0(\mu_1)}{\overline{\lambda}(\A,\mu_1)} \biggr\} = \min\biggl\{ 1,
  -\frac{h_0(\mu)}{\lambda_1(\A,\mu)} \biggr\} = 1.
\end{equation}
By \cref{thm:feng1} and \cref{thm:overline_h1_large}, we have
\begin{equation} \label{eq:eq_state_h1_large2}
  \dim(\overline{\pi}_*\mu_1) =
  \frac{\overline{h}_1(\mu_1)-\overline{h}_0(\mu_1)}{\overline{\lambda}(\A,\mu_1)} =
  \frac{h_1(\mu)-h_0(\mu)}{\lambda_1(\A,\mu)}.
\end{equation}
Combining \cref{eq:eq_state_h1_large1} and \cref{eq:eq_state_h1_large2}, we obtain
$h_1(\mu)=h_0(\mu)+\lambda_1(\A,\mu)$.
\end{proof}

To finish this section, let us show in the following three propositions that $h_2(\mu)$ vanishes
when $h_0(\mu) \le -\lambda_2(\A,\mu)$ and also if the strong open set condition is satisfied. The
first proposition covers the case where the affinity dimension is strictly smaller than one.

\begin{proposition} \label{thm:eq_state_h2_zero}
  If $X$ is a proximal box-like self-affine set with rotations satisfying the projective exponential
  separation condition, $\A$ contains a proximal diagonal generator, and $\mu \in
  \MM_\sigma(\Sigma)$ is the $s$-equilibrium state with $0<s<2$ such that $h_0(\mu) \le
  -\lambda_1(\A,\mu)$, then
  \begin{equation*}
    h_2(\mu)=0.
  \end{equation*}
\end{proposition}

\begin{proof}
  Write $U_i = \{\iii \in \Sigma : \mu_\iii^{\xi_i} = \delta_\iii\}$ and notice that, by
  \cref{thm:eq_state_h1_zero}, $\mu(U_1)=1$. Since $\xi_2(\iii) \subseteq \xi_1(\iii)$ for all $\iii
  \in \Sigma'$, \cref{thm:finer_partition_dirac_mass} shows that $\mu(U_2)=1$. The proof is finished
  by \cref{thm:feng_entropy_zero}.
\end{proof}

The second proposition covers the case where the Kolmogorov--Sinai entropy of the $s$-equilibrium
state satisfies $h_0(\mu) \le -\lambda_2(\A,\mu)$.

\begin{proposition} \label{thm:orto_partition_dirac}
  If $X$ is a proximal box-like self-affine set with rotations satisfying the projective exponential
  separation condition, $\A$ contains a proximal diagonal generator, and $\mu \in
  \MM_\sigma(\Sigma)$ is the $s$-equilibrium state with $0<s<2$ such that
  $h_0(\mu) \le -\lambda_2(\A,\mu)$, then
  \begin{equation*}
    h_2(\mu) = 0.
  \end{equation*}
\end{proposition}

\begin{proof}
  By \cref{thm:overline_h1_zero_large}, the conditional measures of $\overline\mu$ are Dirac
  for $\overline\mu$-almost every point of $\omega(\Sigma^\tau)\cup\tau(\Sigma^\omega)$.
  Apply \cref{thm:conditional_two_ways,eq:orientation-restrictions} with
  $(T,E)=(\omega,\Sigma^\tau)$ or $(\tau,\Sigma^\omega)$ to obtain
  \begin{equation*}
    T_*\mu_\iii^{\xi_1^\bot}=\delta_{T(\iii)}
  \end{equation*}
  for $\mu$-almost every $\iii\in E$.
  Injectivity of $T$ gives $\mu_\iii^{\xi_1^\bot}=\delta_\iii$ for $\mu$-almost every
  $\iii\in\Sigma'$. Since $\xi_2$ refines $\xi_1^\bot$,
  \cref{thm:finer_partition_dirac_mass} gives $\mu_\iii^{\xi_2}=\delta_\iii$ almost everywhere.
  Applying \cref{thm:feng_entropy_zero}, we obtain $h_2(\mu)=0$.
\end{proof}

Finally, we conclude that under the strong open set condition, defined in
\cref{sec:box_like_self_affine}, the $2$-dimensional entropy $h_2$ vanishes.

\begin{proposition} \label{thm:SOSC_h2_zero}
  If $X$ is a proximal box-like self-affine set with rotations satisfying the strong open set
  condition and $\mu \in \MM_\sigma(\Sigma)$ is the $s$-equilibrium state with $0<s<2$, then
  \begin{equation*}
    h_2(\mu)=0.
  \end{equation*}
\end{proposition}

\begin{proof}
  By \cref{thm:irred_eq_state}, $\mu$ is ergodic and fully supported. Let $U$ satisfy the strong
  open set condition. Iterating any point of $U$ shows that $X\subseteq\overline U$.
  Since $X\cap U\ne\emptyset$ and cylinder diameters tend to zero, there is a finite word $w$
  with $\fii_w(X)\subseteq U$. Ergodicity and $\mu([w])>0$ imply that almost every sequence
  contains $w$ infinitely often in its future. For such a sequence $\iii$, every tail
  $\pi(\sigma^k(\iii))$ lies in $U$, by choosing a later occurrence of $w$ and using forward
  invariance of $U$.

  If another one-sided address $\jjj\ne\iii^+$ codes $\pi(\iii)$, let
  $k=\min\{m\ge0:i_m\ne j_m\}$. Since $i_m=j_m$ for $m<k$ and each $\fii_i$ is injective,
  we have $\pi(\sigma^k(\iii))=\pi^+(\sigma^k(\jjj))$. Hence
  \begin{equation*}
    \pi(\sigma^k(\iii))\in\fii_{i_k}(U)\cap\fii_{j_k}(X)
    \subseteq\fii_{i_k}(U)\cap\overline{\fii_{j_k}(U)}=\emptyset,
  \end{equation*}
  since the open images of $U$ are pairwise disjoint. Thus the future address is unique.
  By \cref{eq:simple_lyapunov}, $\mu$ has simple Lyapunov spectrum, and $\xi_2$ fixes both the
  past and the coded point. Hence $\xi_2(\iii)=\{\iii\}$ and $\mu_\iii^{\xi_2}=\delta_\iii$
  almost everywhere. The conclusion follows from \cref{thm:feng_entropy_zero}.
\end{proof}

\subsection{Dimension of box-like self-affine set} \label{sec:dimension_result_proof}
In \cref{thm:main_extended}, we first compute $\dim(\pi_*\mu)$ for the equilibrium state $\mu$ under
projective exponential separation, up to the threshold and under the strong open set condition,
which gives the last claim of \cref{thm:main2}. We then prove the dimension formula of
\cref{thm:main2} using finite diagonal subsystems. In its cases \cref{it:main2_threshold} and
\cref{it:main2_sosc}, the only separation assumption on the projective system is the absence of
exact overlaps.

\begin{theorem} \label{thm:main_extended}
  Let $X$ be a proximal box-like self-affine set with rotations, and let
  $\mu \in \MM_\sigma(\Sigma)$ be the unique $s$-equilibrium state for $s = \dimaff(\A)$.
  For $s<2$, assume additionally that $X$ satisfies the projective exponential separation condition.
  Then, if $h_0(\mu) \le -\lambda_2(\A,\mu)$ or $X$ satisfies the strong open set condition, we have
  \begin{equation*}
    \dim(\pi_*\mu) = \dimaff(\A) \le 2.
  \end{equation*}
\end{theorem}

\begin{proof}
  Since $N\ge2$, we have $P(\A,0)=\log N>0$, so $s>0$. We first consider the case $s\ge2$, which
  does not require the additional hypothesis of the projective exponential separation condition.
  By \cref{eq:eq_state_def,eq:energy_def},
  \begin{equation*}
    h_0(\mu)=-\frac{s}{2}\bigl(\lambda_1(\A,\mu)+\lambda_2(\A,\mu)\bigr)
    \ge-\lambda_1(\A,\mu)-\lambda_2(\A,\mu)>-\lambda_2(\A,\mu).
  \end{equation*}
  Thus the entropy hypothesis cannot hold, and $X$ must satisfy the strong open set condition.
  Choose an open set $U$ as in its definition. Since all the maps are contractions, there is
  a large open ball $B\supset X$ such that $\fii_i(B)\subseteq B$ for every $i$. Then
  $V=U\cap B$ is a nonempty bounded open set whose images $\fii_i(V)$ are pairwise disjoint
  subsets of $V$.

  Let $\LL^2$ denote planar Lebesgue measure. By change of variables,
  \begin{equation*}
    \LL^2(V)\sum_{i=1}^N |\det A_i| = \sum_{i=1}^N\LL^2(\fii_i(V)) \le \LL^2(V).
  \end{equation*}
  Since $0<\LL^2(V)<\infty$, we obtain $\sum_i |\det A_i|\le1$, and hence 
  $P(\A,2)=\log\sum_i |\det A_i|\le0$. The pressure is strictly decreasing and $P(\A,s)=0$, 
  so $s\ge2$ forces $s=2$ and $\sum_i |\det A_i|=1$. Consequently, the sets $\fii_i(V)$ cover 
  $V$ up to a Lebesgue-null set. Define the probability measure
  $\nu=\LL^2|_V/\LL^2(V)$. For every Borel set $E\subseteq\R^2$, we have
  \begin{equation*}
    \sum_{i=1}^N |\det A_i|(\fii_i)_*\nu(E) = \frac{1}{\LL^2(V)}\sum_{i=1}^N\LL^2(E\cap\fii_i(V)) = \nu(E).
  \end{equation*}
  By \cref{thm:bernoulli_equilibrium}, at $s=2$ the equilibrium state $\mu$ is Bernoulli with weights
  $|\det A_i|$. Its projection $\pi_*\mu$ therefore satisfies the same invariance equation.
  Uniqueness of the invariant probability measure for a contractive iterated function system,
  proved by Hutchinson~\cite{Hutchinson1981}, gives $\pi_*\mu=\nu$. Since $V$ is open, $\nu$
  has local dimension two at every point of $V$, and hence is exact-dimensional of dimension
  two.

  We may now assume $0<s<2$ and, by \cref{thm:two-block-normalisation}, that $\A$ contains 
  a proximal diagonal generator. By \cref{thm:irred_eq_state}, the equilibrium state $\mu$ is ergodic and satisfies $\dimaff(\A) = \diml(\A,\mu) = s < 2$. It follows from \cref{thm:feng1} that
  \begin{equation} \label{eq:feng_ledrappier_young}
    \begin{split}
      \dim(\pi_*\mu) 
      &= \frac{h_1(\mu)-h_0(\mu)}{\lambda_1(\A,\mu)} +
      \frac{h_2(\mu)-h_1(\mu)}{\lambda_2(\A,\mu)} \\ 
      &= \frac{h_2(\mu)-h_0(\mu)}{\lambda_2(\A,\mu)} + \biggl(
      1-\frac{\lambda_1(\A,\mu)}{\lambda_2(\A,\mu)} \biggr)
      \frac{h_1(\mu)-h_0(\mu)}{\lambda_1(\A,\mu)}. 
    \end{split}
  \end{equation}
  Recall that, by \cref{eq:simple_lyapunov}, we have $0 < -\lambda_1(\A,\mu) < -\lambda_2(\A,\mu)$.
  Let us first assume that $h_0(\mu) \le -\lambda_1(\A,\mu)$. Since $X$ satisfies the projective
  exponential separation condition, \cref{thm:eq_state_h2_zero,thm:eq_state_h1_zero} show that
  $h_1(\mu)=0$ and $h_2(\mu)=0$. By \cref{eq:feng_ledrappier_young,eq:diml_min2}, we obtain
  \begin{equation} \label{eq:main_extended1}
    \dim(\pi_*\mu) = -\frac{h_0(\mu)}{\lambda_1(\A,\mu)} = \diml(\A,\mu).
  \end{equation}
  Suppose next that $-\lambda_1(\A,\mu) < h_0(\mu) \le -\lambda_2(\A,\mu)$.
  In this case, since $X$ satisfies the projective exponential separation condition,
  \cref{thm:eq_state_h1_large,thm:orto_partition_dirac} give $h_1(\mu) = h_0(\mu)+\lambda_1(\A,\mu)$
  and $h_2(\mu)=0$. By \cref{eq:feng_ledrappier_young,eq:diml_min2}, we obtain
  \begin{equation} \label{eq:main_extended2}
    \dim(\pi_*\mu) = -\frac{h_0(\mu)}{\lambda_2(\A,\mu)} + \biggl( 1-\frac{\lambda_1(\A,\mu)}{\lambda_2(\A,\mu)} \biggr) = \diml(\A,\mu),
  \end{equation}
  finishing the proof under the assumption $h_0(\mu) \le -\lambda_2(\A,\mu)$.

  If $X$ satisfies the strong open set condition, then \cref{thm:SOSC_h2_zero} shows that 
  $h_2(\mu)=0$. If $h_0(\mu) \le -\lambda_2(\A,\mu)$, then the preceding cases apply. We may 
  therefore assume that $h_0(\mu) > -\lambda_2(\A,\mu) > -\lambda_1(\A,\mu)$. Then \cref{thm:eq_state_h1_large} gives $h_1(\mu)=h_0(\mu)+\lambda_1(\A,\mu)$, and since $1<s<2$ by \cref{eq:diml_min2}, we again obtain \cref{eq:main_extended2}.
\end{proof}

To prove \cref{thm:main2}, we approximate the entropy and Lyapunov exponents of equilibrium states
by those of Bernoulli measures on finite diagonal subsystems. This allows us to apply
Kittle--Kogler~\cite[Theorem~1.1]{KittleKogler2026-preprint} to their coordinate projections and, in
case \cref{it:main2_pesc}, \cref{thm:rapaport} to the subsystems themselves.

\begin{proof}[Proof of \cref{thm:main2}]
  Write $s=\dimaff(\A)>0$, and consider first the case $s\ge2$.
  Then $2-\lambda_1(\A,\mu)/\lambda_2(\A,\mu)<2\le s$ by the negativity of both Lyapunov exponents
  (see \cref{sec:lyapunov_exponents}), so \cref{it:main2_threshold} fails. If
  \cref{it:main2_sosc} holds, \cref{thm:main_extended} gives $s=2$ and $\dim(\pi_*\mu)=2$.
  Since $\pi_*\mu$ is supported on $X\subseteq\R^2$, we obtain $\dimh(X)=2=s$.
  Under \cref{it:main2_pesc}, it suffices to prove $\dimh(X)\ge t$ for every $0<t<2$.

  For the remaining cases, set $t=s$ if $s<2$, and otherwise fix $0<t<2$.
  Passing to the second iterate if necessary, we may assume by \cref{thm:two-block-normalisation}
  that $\A$ contains a proximal diagonal generator. This preserves the attractor, the affinity
  dimension and the ratio of the Lyapunov exponents, as well as the strong open set and
  projective exponential separation conditions.
  The absence of exact overlaps also passes to the second iterate: each of its
  projective generators corresponds to a distinct two-letter word in $\Psi$, as in the proof of
  that lemma, and concatenation of these words is injective.
  Let $\mu_t$ be the $t$-equilibrium state, and write $h=h_0(\mu_t)$ and
  $\lambda_i=\lambda_i(\A,\mu_t)$ for $i\in\{1,2\}$. By \cref{eq:eq_state_def} and the
  monotonicity of the pressure,
  \begin{equation*}
    h+\Lambda(\A,\mu_t,t)=P(\A,t)\ge P(\A,s)=0.
  \end{equation*}
  Hence $\diml(\A,\mu_t)\ge t$, with equality when $t=s$.
  By \cref{eq:simple_lyapunov}, we also have $\lambda_2<\lambda_1<0$.

  Let $\mu_1$ be the Gibbs measure on $\Sigma_Z$ associated to $f_{1,t}$. By
  \cref{thm:entropy_eq_gibbs,thm:eq_state_is_sum}, its entropy is $h$, and
  \begin{equation*}
    \int f_{1,1}\dd\mu_1=\lambda_1 \quad\text{and}\quad \int f_{2,1}\dd\mu_1=\lambda_2.
  \end{equation*}
  Since $f_{1,t}$ depends only on the current symbol and $Z$ is mixing, $\mu_1$ is a mixing
  Markov measure with full support. Mixing gives
  \begin{equation*}
    \mu_1\bigl(\tau(\Sigma)\cap\sigma^{-n}\tau(\Sigma)\bigr)
    \to\mu_1(\tau(\Sigma))^2>0.
  \end{equation*}
  For $0<2\eps<\lambda_1-\lambda_2$ and $n\ge1$, define
  \begin{align*}
    W_{n,\eps} = \{\iii|_n : \;&\iii\in\tau(\Sigma)\cap\sigma^{-n}\tau(\Sigma), \\ 
    &e^{-n(h+\eps/2)}\le\mu_1([\iii|_n])\le e^{-n(h-\eps/2)}, \\
    &\text{and }\max_{k\in\{1,2\}}|n^{-1}S_n f_{k,1}(\iii)-\lambda_k|<\eps\}.
  \end{align*}
  The first condition implies that we only take words that return to the same orientation, so their
  associated matrices are diagonal. For each fixed $\eps>0$ in this range, the
  Shannon--McMillan--Breiman and Birkhoff theorems, together with the mixing statement above, give
  \begin{equation*}
    \sum_{\iii\in W_{n,\eps}}\mu_1([\iii])\to\mu_1(\tau(\Sigma))^2>0.
  \end{equation*}
  Consequently, for all sufficiently large $n$,
  \begin{equation*}
    \exp(n(h-\eps))\le\#W_{n,\eps}\le\exp(n(h+\eps)).
  \end{equation*}

  Write $\mathbf{k}=(k_1,\ldots,k_{2N})\in\{0,\ldots,n\}^{2N}$ for the vector of digit frequencies and
  set
  \begin{equation*}
    W_{n,\eps}(\mathbf{k})=\{\iii\in W_{n,\eps}\colon
    \#\{p \in \{0,\ldots,n-1\} : i_p=j\}=k_j\text{ for }j\in\{1,\ldots,2N\}\}.
  \end{equation*}
  By pigeonholing, we may choose $\mathbf{k}$ such that $\#W_{n,\eps}(\mathbf{k}) \ge
  (n+1)^{-2N}\exp(n(h-\eps))$.
  By construction, for all words in $W_{n,\eps}(\mathbf{k})$ both diagonal entries are products of the
  same scalar factors with the same multiplicities, so the associated planar matrices coincide.
  Letting $\eps$ tend to zero and choosing sufficiently large $n$, define
  \begin{equation*}
    W_n=q\bigl(W_{n,\eps}(\mathbf{k})\bigr)\subseteq\Sigma_n^+,
  \end{equation*}
  where the map $q$, defined in \cref{sec:entropy_eq_state}, allows us to pass from the doubled
  alphabet to $\{1,\ldots,N\}$.
  Since all selected words start in the same horizontal orientation, $q$ is injective on
  $W_{n,\eps}(\mathbf{k})$ and thus $\#W_n=\#W_{n,\eps}(\mathbf{k})$. We may write
  $A_w=\operatorname{diag}(a_n,b_n)$ for all $w\in W_n$, with $0<|b_n|<|a_n|<1$ by contraction
  and the choice of $\eps$.
  The preceding estimates give
  \begin{equation} \label{eq:diagonal-subsystem-limits}
    \frac{\log\#W_n}{n}\to h,\quad
    \frac{\log|a_n|}{n}\to\lambda_1,\quad\text{and}\quad
    \frac{\log|b_n|}{n}\to\lambda_2,
  \end{equation}
  as $n\to\infty$.

  Let $\nu_n$ be the Bernoulli measure on $W_n^\Z$ assigning weight $(\#W_n)^{-1}$ to each
  symbol, and let $\pi_n\colon W_n^\Z\to\R^2$ be the coding map for the subsystem
  $(\fii_w)_{w\in W_n}$. Then $h_0(\nu_n)=\log\#W_n$, and the Lyapunov exponents of $\nu_n$ are
  $\log|a_n|$ and $\log|b_n|$.
  The systems induced on the coordinate axes have no exact overlaps: their compositions
  correspond to distinct equal-length words in $\Psi$. Hence
  Kittle--Kogler~\cite[Theorem~1.1]{KittleKogler2026-preprint} gives
  \begin{align*}
    \dim((\proj_x\circ\pi_n)_*\nu_n)&=\min\biggl\{1,\frac{\log\#W_n}{-\log|a_n|}\biggr\},\\
    \dim((\proj_y\circ\pi_n)_*\nu_n)&=\min\biggl\{1,\frac{\log\#W_n}{-\log|b_n|}\biggr\}.
  \end{align*}
  When $s<2$, we have $t=s$ and
  $\log\#W_n/(-\log|a_n|)\to-h/\lambda_1\ge\min\{1,s\}$ by \cref{eq:diml_min2}.
  Since $(\proj_x\circ\pi_n)_*\nu_n$ is supported on $\proj_x(X)$, it follows that
  $\dimh(\proj_x(X))\ge\min\{1,s\}$. An anti-diagonal generator maps $\proj_x(X)$ bi-Lipschitzly into
  $\proj_y(X)$, and so both coordinate
  projections of $X$ have Hausdorff (and Minkowski) dimension $\min\{1,s\}$.
  Thus, if $s\le1$, we have $s=\dimh(\proj_x(X))\le\dimh(X)\le s$ by
  \cref{eq:hausdorff_affinity_upper_bound}, proving \cref{thm:main2} for $s\in(0,1]$.

  We now prove $\dimh(X)\ge t$ under assumption \cref{it:main2_pesc}. Let $t=s$ when $1<s<2$ and
  choose $t$ arbitrarily close to $2$ when $s\ge 2$.
  Let $X_n\subseteq X$ be the attractor of
  $(\fii_w)_{w\in W_n}$. Fix $0<c<1$ such that $\Delta_m^\Psi\ge c^m$ for infinitely many $m$.
  A composition of $k$ subsystem maps acts on either coordinate axis by a map $g_\iii$
  with $|\iii|=kn$. Distinct compositions give distinct words, and on a given axis the
  derivatives are equal. Given $m\ge n$ with $\Delta_m^\Psi\ge c^m$, put $k=\lfloor m/n\rfloor$,
  so that $kn\le m<2kn$. By the monotonicity of $(\Delta_j^\Psi)_{j\in\N}$ established in
  \cref{sec:exponential_separation}, any two such words $\iii\ne\jjj$ satisfy
  \begin{equation*}
    |g_\iii(0)-g_\jjj(0)|\ge\Delta_{kn}^\Psi\ge\Delta_m^\Psi\ge c^m\ge(c^{2n})^k.
  \end{equation*}
  Hence, both induced systems are exponentially separated.
  Since $|b_n|<|a_n|$, the subsystem is proximal, and \cref{thm:rapaport} gives
  $\dimh(X_n)=\min\{2,\dimaff(X_n)\}$.

  The pressure of this subsystem at $u\ge0$ is
  $\log\#W_n+\log\fii^u(\operatorname{diag}(a_n,b_n))$.
  By \cref{eq:diml_min2} and the three limits above,
  $\dimaff(X_n)\to\diml(\A,\mu_t)$. Thus
  \begin{equation*}
    \dimh(X)\ge\lim_{n \to \infty}\dimh(X_n)=\min\{2,\diml(\A,\mu_t)\}\ge t.
  \end{equation*}
  Since $t=s$, or $t$ is arbitrarily close to $2$, we obtained the required lower bound.
  The upper bound is \cref{eq:hausdorff_affinity_upper_bound} and we have proved \cref{thm:main2}
  under assumption \cref{it:main2_pesc}.

  It remains to treat \cref{it:main2_threshold,it:main2_sosc} for $1<s<2$.
  Since $t=s$, we have $\diml(\A,\mu_t)=s$, and
  \cref{eq:diml_min2,eq:simple_lyapunov} show that \cref{it:main2_threshold} is equivalent to
  $h\le-\lambda_2$.
  Since $\nu_n$ is Bernoulli and the Oseledets subspace $V^1$ is the $y$-axis,
  $h_1(\nu_n)$ is the conditional entropy of the first symbol given the horizontal coordinate
  of the coded point. The scalar Ledrappier--Young formula
  \cite[Theorems~1.2 and~1.3]{Feng2023} gives
  \begin{equation*}
    h_1(\nu_n)=\log\#W_n+\log|a_n|\dim((\proj_x\circ\pi_n)_*\nu_n).
  \end{equation*}
  We now modify $W_n$ to obtain $h_2(\nu_n)=0$ while preserving the limits in
  \cref{eq:diagonal-subsystem-limits} and will continue to write $\nu_n$ and $\pi_n$ for the
  uniform Bernoulli measure and coding map of the modified subsystem.

  Assuming \cref{it:main2_threshold}, we redefine $W_n$ by retaining any (arbitrary) words such that
  $\#W_n = \min\{\#W_n,\lfloor|b_n|^{-1}\rfloor\}$.
  Since $h\le-\lambda_2$, the limits in \cref{eq:diagonal-subsystem-limits} are preserved
  and $\#W_n\le|b_n|^{-1}$. Hence, the conditional entropy is
  \begin{equation*}
    H_{\nu_n}(\cP\,|\,\widehat{\xi_1^\bot})
    =\log\#W_n+\log|b_n|\dim((\proj_y\circ\pi_n)_*\nu_n)=0.
  \end{equation*}
  Further, since $\xi_2$ refines $\xi_1^\bot$, we obtain $h_2(\nu_n)=0$.

  Assuming instead \cref{it:main2_sosc}, we let $U$ be a feasible open set and choose a word $v$ 
  with $\fii_v(X)\subseteq U$, as in the proof of \cref{thm:SOSC_h2_zero}. By appending an anti-diagonal generator if necessary, we may assume that $A_v$ is diagonal. 
  We now replace $W_n$ by $\{vw : w\in W_n\}$, leaving the cardinality unchanged, while the lengths 
  and logarithms of the coordinate contraction ratios change by bounded amounts. 
  Thus the limits in \cref{eq:diagonal-subsystem-limits} hold and $|b_n|<|a_n|$ for all sufficiently large $n$. 
  The resulting attractor lies in $\fii_v(X)\subseteq U$, so its first-level images lie 
  in pairwise disjoint level-$(|v|+n)$ images of $U$. The subsystem therefore satisfies the strong 
  separation condition, which gives unique codings and hence $h_2(\nu_n)=0$.

  Using the projection formula and the three limits, \cref{thm:feng1,eq:diml_min2} now give
  \begin{equation*}
    \dim((\pi_n)_*\nu_n)
    =\frac{h_1(\nu_n)-\log\#W_n}{\log|a_n|}
     -\frac{h_1(\nu_n)}{\log|b_n|}
    \to s.
  \end{equation*}
  Since $(\pi_n)_*\nu_n$ is supported on $X$, we obtain $\dimh(X)\ge s$. Combining this with the
  upper bound \cref{eq:hausdorff_affinity_upper_bound} gives the desired conclusion.

  Finally, we prove the last assertion in \cref{thm:main2}. 
  Suppose that \cref{it:main2_pesc} holds together with \cref{it:main2_threshold} or
  \cref{it:main2_sosc}. Under \cref{it:main2_threshold}, we have $s<2$, and
  \cref{eq:diml_min2,eq:simple_lyapunov} give
  $h_0(\mu)\le-\lambda_2(\A,\mu)$. In either case, the hypotheses of
  \cref{thm:main_extended} are satisfied, and hence $\dim(\pi_*\mu)=\dimaff(X)\le2$.
\end{proof}

\begin{remark} \label{rem:threshold_bound}
  If $X$ is a proximal box-like self-affine set with rotations whose projective system has no exact
  overlaps, $\dimaff(X)<2$, and $\mu$ is the unique $s$-equilibrium state for $s=\dimaff(X)$, then
  the proof of \cref{thm:main2} gives
  \begin{equation*}
    \dimh(X) \ge \min\biggl\{\dimaff(X),2-\frac{\lambda_1(\A,\mu)}{\lambda_2(\A,\mu)}\biggr\}
  \end{equation*}
  without assuming the strong open set condition or projective exponential separation. In the
  notation of that proof, suppose that $h>-\lambda_2$ and reduce $W_n$ to $\lfloor|b_n|^{-1}\rfloor$
  words. Then $\#W_n\le|b_n|^{-1}$ gives $h_2(\nu_n)=0$ by the same argument.
  Since $n^{-1}\log\#W_n\to-\lambda_2>-\lambda_1$, we have
  $\dim((\proj_x\circ\pi_n)_*\nu_n)=1$ for all sufficiently large $n$, and hence
  \begin{equation*}
    \dim((\pi_n)_*\nu_n)=1+\frac{\log\#W_n+\log|a_n|}{-\log|b_n|}
    \to2-\frac{\lambda_1}{\lambda_2}.
  \end{equation*}
  The bound $\log\#W_n\le-\log|b_n|$, which gives $h_2(\nu_n)=0$, is the analogue of the hypothesis of
  B\'ar\'any--Rams--Simon~\cite[Theorem~A(2)]{BaranyRamsSimon2016}, who use the theorem of
  Hochman~\cite[Theorem~1.1]{Hochman2014} whereas we use that of Kittle--Kogler. Above the threshold
  and without the reduction, the conditional entropy given the vertical coordinate is
  $\log\#W_n+\log|b_n|>0$ for all sufficiently large $n$, which does not force $h_2(\nu_n)=0$.
  If \cref{thm:rapaport} held with exponential separation of the induced
  systems replaced by the absence of exact overlaps, then the argument for \cref{it:main2_pesc}
  would prove \cref{thm:main2} without any of the conditions
  \cref{it:main2_threshold}--\cref{it:main2_pesc}, see
  Rapaport~\cite[Remark~1.6]{Rapaport2023-preprint}.
\end{remark}

\bibliographystyle{abbrv}
\bibliography{Bibliography}

\end{document}